# Invariant measures for Burgers equation with stochastic forcing

By WEINAN E, K. KHANIN, A. MAZEL, and YA. SINAI

## 1. Introduction

In this paper we study the following Burgers equation

$$(1.1) \qquad \frac{\partial u}{\partial t} + \frac{\partial}{\partial x}\Big(\frac{u^2}{2}\Big) = \varepsilon \frac{\partial^2 u}{\partial x^2} + f(x,t)$$

where $f(x,t) = \frac{\partial F}{\partial x}(x,t)$ is a random forcing function, which is periodic in $x$ with period 1, and with white noise in $t$. The general form for the potentials of such forces is given by:

$$(1.2) \qquad F(x,t) = \sum_{k=1}^{\infty} F_k(x) \dot{B}_k(t)$$

where the $\{B_k(t), t \in (-\infty, \infty)\}$'s are independent standard Wiener processes defined on a probability space $(\Omega, \mathcal{F}, P)$ and the $F_k$'s are periodic with period 1. We will assume for some $r \geq 3$,

$$(1.3) \qquad f_k(x) = F_k'(x) \in C^r(\mathbb{S}^1), \ \|f_k\|_{C^r} \leq \frac{C}{k^2}.$$

Here $\mathbb{S}^1$ denotes the unit circle, and $C$, a generic constant. Without loss of generality, we can assume that for all $k$, $\int_0^1 F_k(x) dx = 0$. We will denote the elements in the probability space $\Omega$ by $\omega = (\dot{B}_1(\cdot), \dot{B}_2(\cdot), \ldots)$. Except in Section 8 where we study the convergence as $\varepsilon \to 0$, we will restrict our attention to the case when $\varepsilon = 0$:

$$(1.4) \qquad \frac{\partial u}{\partial t} + \frac{\partial}{\partial x}\Big(\frac{u^2}{2}\Big) = \frac{\partial F}{\partial x}(x,t).$$

Besides establishing existence and uniqueness of an invariant measure for the Markov process corresponding to (1.4), we will also give a detailed description of the structure and regularity properties for the solutions that live on the support of this measure.

The randomly forced Burgers equation (1.1) is a prototype for a very wide range of problems in nonequilibrium statistical physics where strong nonlinear effects are present. It arises in studies of various one-dimensional systems such as vortex lines in superconductors [BFGLV], charge density waves [F], directed polymers [KS], etc. (1.1) with its high-dimensional analog is the differentiated version of the well-known KPZ equation describing, among other things, kinetic roughening of growing surfaces [KS]. Most recently, (1.1) has received a great deal of interest as the testing ground for field-theoretic techniques in hydrodynamics [CY], [Pol], [GM], [BFKL], [GK2]. In fact, we expect that the randomly forced Burgers equation will play no lesser a role in the understanding of nonlinear non-equilibrium phenomena than that of the Burgers equation in the understanding of nonlinear waves.

Before proceeding further let us give an indication why an invariant measure is expected for (1.1) even when $\varepsilon = 0$. Since energy is continuously supplied to the system, a dissipation mechanism has to be present to maintain an invariant distribution. In the case when $\varepsilon > 0$, the viscous term provides the necessary energy dissipation, and the existence of an invariant measure has already been established in [S1], [S2]. When $\varepsilon = 0$, it is well-known that discontinuities are generally present in solutions of (1.4) in the form of shock waves [La]. These weak solutions are limits of solutions of (1.1) as $\varepsilon \to 0$, and satisfy an additional entropy condition: $u(x+,t) \leq u(x-,t)$, for all $(x,t)$. It turns out that this entropy condition enforces sufficient energy dissipation (in the shocks) for maintaining an invariant measure. We will always restrict our attention to weak solutions of (1.4) that satisfy the entropy condition.

The starting point of our analysis is the following variational characterization of solutions of (1.4) satisfying the entropy condition [Li]:

For any Lipschitz continuous curve $\xi\colon [t_1, t_2] \to \mathbb{S}^1$, define its action

$$(1.5) \qquad \mathcal{A}_{t_1,t_2}(\xi) = \int_{t_1}^{t_2} \left\{ \frac{1}{2} \dot{\xi}(s)^2 ds + \sum_k F_k(\xi(s)) dB_k(s) \right\}.$$

Then for $t > \tau$, solutions of (1.4) satisfy

$$(1.6) \qquad u(x,t) = \frac{\partial}{\partial x} \inf_{\xi(t)=x} \left\{ \mathcal{A}_{\tau,t}(\xi) + \int_0^{\xi(\tau)} u(z,\tau) dz \right\}$$

where the infimum is taken over all Lipschitz continuous curves on $[\tau, t]$ satisfying $\xi(t) = x$.

Here and below, we avoid in the notation explicit indication of the dependence on realization of the random force when there is no danger of confusion. Otherwise we indicate such dependence by a super- or subscript $\omega$. In addition, we will denote by $\theta^\tau$ the shift operator on $\Omega$ with increment $\tau$: $\theta^\tau \omega(t) = \omega(t + \tau)$, and by $S_\omega^\tau w$ the solution of (1.1) at time $t = \tau$ when the realization of the force is $\omega$ and the initial datum at time $t = 0$ is $w$. We will denote by $D$ the Skorohod space on $\mathbb{S}^1$ (see [B], [Pa]) consisting of functions having only discontinuities of the first kind; i.e., both left and right limits exist at each point, but they may not be equal.

It is easy to see that the dynamics of (1.4) conserves the quantity $\int_0^1 u(x,t)dx$. Therefore to look for unique invariant measure, we must restrict attention to the subspace

$$D_c = \{u \in D, \int_0^1 u(x)dx = c\}.$$

In this paper we will restrict most of our attention to the case when $c = 0$ but it is relatively easy to see that all of our results continue to hold for the case when $c \neq 0$. We will come back to this point at the end of this section. At the end of Section 3, we will outline the necessary changes for the case when $c \neq 0$.

Our basic strategy for the construction of an invariant measure is to show that the following "*one force, one solution*" principle holds for (1.4): For almost all $\omega$, there exists a unique solution of (1.4), $u^\omega$, defined on the time interval $(-\infty, +\infty)$. In other words, the random attractor consists of a single trajectory almost surely. Furthermore, if we denote the mapping between $\omega$ and $u^\omega$ by $\Phi$:

(1.7) $$u^\omega = \Phi(\omega),$$

then $\Phi$ is invariant in the sense that

(1.8) $$\Phi(\theta^\tau \omega) = S_\omega^\tau \Phi(\omega).$$

It is easy to see that if such a map exists, then the distribution of $\Phi_0 : \Omega \to D$:

$$\Phi_0(\omega)(x) = u^\omega(x, 0),$$

is an invariant measure for (1.4). Moreover, this invariant measure is necessarily unique.

This approach of constructing the invariant measure has the advantage that many statistical properties of the forces, such as ergodicity and mixing, carry over automatically to the invariant measure. More importantly, it facilitates the study of solutions supported by the invariant measure, i.e. the associated stationary Markov process. This study will be carried out in the second half of the present paper.

The construction of $u^\omega$ will be accomplished in Section 3. The variational principle (1.6) allows us to restrict our attention to $t = 0$.

Our construction of $\Phi$ relies heavily on the notion of one-sided minimizer. A curve $\xi$ $(-\infty, 0] \to \mathbb{S}^1$ is called a one-sided minimizer if it minimizes the action (1.5) with respect to all compact perturbations. More precisely, we introduce:

*Definition* 1.1. A piecewise $C^1$-curve $\{\xi(t),\, t \leq 0\}$ is a *one-sided minimizer* if for any Lipschitz continuous $\tilde{\xi}$ defined on $(-\infty, 0]$ such that $\tilde{\xi}(0) = \xi(0)$ and $\tilde{\xi} = \xi$ on $(-\infty, \tau]$ for some $\tau < 0$,
$$\mathcal{A}_{s,0}(\xi) \leq \mathcal{A}_{s,0}(\tilde{\xi})$$
for all $s \leq \tau$.

It is important to emphasize that the curves are viewed on the cylinder $\mathbb{R}^1 \times \mathbb{S}^1$. Similarly, we define one-sided minimizers on $(-\infty, t]$, for $t \in \mathbb{R}^1$.

The interest of this notion lies in the fact that we are considering an infinite interval. It is closely related to the notion of geodesics of type A introduced and studied by Morse [Mo] and Hedlund [H] and the notion of global minimal orbits in Aubry-Mather theory [A], [M]. In the geometric context, it has been studied by Bangert (see [Ba]) as geodesic rays. A somewhat surprising result is that, in the random case, one-sided minimizers are almost unique. More precisely, we have:

THEOREM 1.1. *With probability* 1, *except for a countable set of $x$ values, there exists a unique one-sided minimizer $\xi$, such that $\xi(0) = x$.*

This theorem states that one-sided minimizers are intrinsic objects to $(x, \omega)$. It allows us to construct $\Phi_0(\omega)$ by patching together all one-sided minimizers:

(1.9) $$\Phi_0\{\omega(\tau),\ \tau < 0\}(x) = u^\omega(x, 0) = \dot{\xi}(0)$$

where $\xi$ is the unique one-sided minimizer such that $\xi(0) = x$. In (1.9) we emphasized the fact that $\Phi_0$ depends only on the realization of $\omega$ in the past $\tau < 0$. Now (1.9) defines $u^\omega(\cdot, 0)$ except on a countable subset of $\mathbb{S}^1$. Similarly we construct $u^\omega(\cdot, t)$ for other values of $t \in \mathbb{R}^1$. It is easy to verify that this construction is self-consistent and satisfies the invariance condition (1.8), as a consequence of the variational principle (1.6).

The existence part of Theorem 1.1 is proved by studying limits of minimizers on finite intervals $[-k, 0]$ as $k \to +\infty$. The uniqueness part of Theorem 1.1 is proved by studying the intersection properties of one-sided minimizers. The absence of two intersections of two different minimizers is a general fact in calculus of variations. However, we will prove the absence of even one intersection which is a consequence of randomness.

We are now ready to define formally the invariant measure. There are two alternative approaches. Either we can define the invariant measure on the product space $(\Omega \times D_0, \mathcal{F} \times \mathcal{D})$ with a skew-product structure, or we can define it as an invariant distribution of the Markov process on $(D_0, \mathcal{D})$ defined by (1.4), where $\mathcal{D}$ is the $\sigma$-algebra generated by Borel sets on $D_0$. The skew-product structure is best suited for the exploration of the "one force, one solution" principle.

*Definition* 1.2. A measure $\mu(du, d\omega)$ on $(\Omega \times D_0, \mathcal{F} \times \mathcal{D})$ is called an *invariant measure* if it is preserved under the skew-product transformation $F^t \colon \Omega \times D_0 \to \Omega \times D_0$,

$$(1.10) \qquad F^t(\omega, u_0) = (\theta^t \omega, S^t_\omega u_0),$$

and if its projection to $\Omega$ is equal to $P$.

Alternatively we may consider a homogeneous Markov process on $D_0$ with the transition probability

$$(1.11) \qquad P_t(u, A) = \int_\Omega \chi_A(u, \omega) P(d\omega)$$

where $u \in D_0$, $A \in \mathcal{D}$, and

$$(1.12) \qquad \chi_A(u, \omega) = \begin{cases} 1 & \text{if } S^t_\omega u \in A \\ 0 & \text{otherwise}. \end{cases}$$

*Definition* 1.3. An *invariant measure* $\kappa(du)$ of the Markov process (1.11) is a measure on $(D_0, \mathcal{D})$ satisfying

$$(1.13) \qquad \kappa(A) = \int_{D_0} P_t(u, A) \kappa(du)$$

for any Borel set $A \in \mathcal{D}$ and any $t > 0$.

Let $\delta^\omega(du)$ be the atomic measure on $(D_0, \mathcal{D})$ concentrated at $\Phi_0(\omega) = u^\omega(\cdot, 0)$, and let $\mu(du, d\omega) = \delta^\omega(du) P(d\omega)$; we then have:

THEOREM 1.2. *If $\mu$ is an invariant measure for the skew-product transformation $F^t$, it is the unique invariant measure on $(\Omega \times D_0, \mathcal{F} \times \mathcal{D})$ with the given projection $P(d\omega)$ on $(\Omega, \mathcal{F})$.*

THEOREM 1.3. *For the Markov process (1.11), $\kappa(du) = \int_\Omega \mu(du, d\omega)$ is the unique invariant measure.*

The uniqueness result is closely related to the uniqueness of one-sided minimizers and reflects the lack of memory in the dynamics of (1.4): Consider solutions of (1.4) with initial data $u(x, -T) = u_0(x)$. Then for almost all $\omega \in \Omega$ and any $t \in \mathbb{R}^1$, $\lim_{T \to +\infty} u(\cdot, t)$ exists and does not depend on $u_0$. The key step in the proof of uniqueness is to prove a strengthened version of this statement.

In the second half of this paper, we study in detail the properties of solutions supported by the invariant measure. The central object is the two-sided minimizer which is defined similarly to the one-sided minimizer but for the interval $(-\infty, +\infty) = \mathbb{R}^1$. Under very weak nondegeneracy conditions, we prove that almost surely, the two-sided minimizer exists and is unique. In Section 6, we show that the two-sided minimizer is a hyperbolic trajectory of the dynamical system corresponding to the characteristics of (1.4):

$$\frac{dx}{dt} = u, \ \frac{du}{dt} = \frac{\partial F}{\partial x}(x,t).$$

We can therefore consider the stable and unstable manifolds of the two-sided minimizer using Pesin theory [Pes]. As a consequence, we show:

THEOREM 1.4. *With probability* 1, *the graph of* $\Phi_0(\omega)$ *is a subset of the unstable manifold* (*at* $t = 0$) *of the two-sided minimizer.*

We use this statement to show that, almost surely, $u^\varepsilon(\cdot, 0)$ is piecewise smooth and has a finite number of discontinuities. This is done in Section 7.

Dual to the two-sided minimizer is an object called the main shock which is a continuous shock curve $x^\omega \colon \mathbb{R}^1 \to \mathbb{S}^1$ defined on the whole line $-\infty < t < \infty$. The main shock is also unique. Roughly speaking, the main shock plays the role of an attractor for the one-sided minimizers while the two-sided minimizer plays the role of a repeller.

Finally in Section 8, we show that as $\varepsilon \to 0$, the invariant measures of (1.1) constructed in [S1], [S2] converge to the invariant measure of (1.4).

The results of this paper have been used to analyze the asymptotic behavior of tail probabilities for the gradients and increments of $u$ (see [EKMS]). It also provides the starting point for the work in [EV] on statistical theory of the solutions. These results are of direct interest to physicists since they can be compared with predictions based on field-theoretic methods (see [Pol], [GM], [GK2], [CY]).

Our theory is closely related to the Aubry-Mather theory [A], [M] which is concerned with special invariant sets of twist maps obtained from minimizing the action

$$(1.14) \qquad \tfrac{1}{2}\sum_i (x_i - x_{i-1} - \gamma)^2 + \lambda \sum_i V(x_i)$$

where $\gamma$ is a parameter and $V$ is a periodic function. The continuous version of (1.14) is

$$(1.15) \qquad \int \{\tfrac{1}{2}(\dot\xi(t) - a)^2 + F(\xi(t), t)\} dt$$

where $F$ is a periodic function in $x$ and $t$ [Mo]. The main result of the Aubry-Mather theory is the existence of invariant sets with arbitrary rotation number, with a suitable $a$. Such invariant sets are made from the velocities of the two-sided minimizers defined earlier. It can be proved that such an invariant set lies on the graph of the periodic solutions of (1.4) [E], [JKM], [So]. In this connection, the results of this paper apply to the random version of (1.15):

$$\text{(1.16)} \qquad \int \{\tfrac{1}{2}(\dot{\xi}(t) - a)^2 dt + \sum_k F_k(\xi(t))dB_k(t)\}.$$

Although only $a = 0$ is considered in this paper, extension to arbitrary $a$ is straightforward and the results are basically the same for different values of $a$. This is because, over a large interval of duration $T$, the contribution of the kinetic energy is of order $T$, and the contribution from the potential is typically of order $\sqrt{T}$ for the random case but of order $T$ for the periodic case. This gives rise to subtle balances between kinetic and potential energies in the latter case. Consequently the conclusions for the random case become much simpler. While in the deterministic case, there are usually many different two-sided minimizers in the invariant set and they are not necessarily hyperbolic, there is only one two-sided minimizer in the random case and it is always hyperbolic.

The value of $a$ is closely related to the value of $c$ discussed earlier. In the setting of Aubry-Mather theory, $a$ is the average speed of the global minimizers and is related to $c$ through the Legendre transform of the homogenized Hamiltonian. In the random case, $a = c$ for the reason given in the last paragraph.

## 2. The variational principle

Let us first define in the probabilistic context the notion of weak solutions of (1.4) with (deterministic) initial data $u(x, t_0) = u_0(x)$. We will always assume $u_0 \in L^\infty(\mathbb{S}^1)$.

*Definition* 2.1. Let $u_\omega$ be a random field parametrized by $(x, t) \in \mathbb{S}^1 \times [t_0, +\infty)$ such that for almost all $\omega \in \Omega$, $u_\omega(\cdot, t) \in D$ for all $t \in (t_0, \infty)$. Then $u_\omega$ is a *weak solution* of (1.4) if:

(i) For all $t > t_0$, $u_\omega(\cdot, t)$ is measurable with respect to the $\sigma$-algebra $\mathcal{F}_{t_0}^t$ generated by all $\dot{B}_k(s)$, $t_0 \leq s \leq t$.

(ii) $u_\omega \in L^1_{\text{loc}}(\mathbb{S}^1 \times [t_0, \infty))$ almost surely.

(iii) With probability 1, the following holds for all $\varphi \in C^2(\mathbb{S}^1 \times \mathbb{R}^1)$ with compact support:

$$\int_0^1 u_0(x)\varphi(x, t_0)dx + \int_{t_0}^\infty \int_0^1 \frac{\partial \varphi}{\partial t} u_\omega(x, t)\, dx\, dt + \frac{1}{2} \int_{t_0}^\infty \int_0^1 \frac{\partial \varphi}{\partial x} u_\omega^2(x, t) dx\, dt$$
$$= -\int_0^1 \sum_k \left\{ F_k(x) \int_{t_0}^\infty \frac{\partial^2 \varphi}{\partial x \partial t}(x, t)(B_k(t) - B_k(t_0)) dt \right\} dx.$$

Also, $u_\omega$ is an entropy-weak solution if, for almost all $\omega \in \Omega$,

$$u_\omega(x+, t) \leq u_\omega(x-, t)$$

for all $(x, t) \in \mathbb{S}^1 \times (t_0, \infty)$.

Our analysis is based on a variational principle characterizing entropy weak solutions of (1.4). To formulate this variational principle, we redefine the action in order to avoid using stochastic integrals. Given $\omega \in \Omega$, for any Lipschitz continuous curve $\xi: [t_1, t_2] \to \mathbb{S}^1$, define

(2.1)
$$\mathcal{A}_{t_1,t_2}(\xi) = \int_{t_1}^{t_2} \left\{ \tfrac{1}{2}\dot{\xi}(s)^2 - \sum_k f_k(\xi(s))\dot{\xi}(s)(B_k(s) - B_k(t_1)) \right\} ds$$
$$+ \sum_k F_k(\xi(t_2))(B_k(t_2) - B_k(t_1)).$$

(2.1) can be formally obtained from (1.5) with an integration by parts. It has the advantage that the integral in (2.1) can be understood in the Lebesgue sense instead of the Ito sense, for example.

LEMMA 2.1. *Let $u_0(x) \in D$. For almost all $\omega \in \Omega$, there exists a unique weak solution of* (1.4) *satisfying the entropy condition, such that $u(x, t_0) = u_0(x)$. For $t \geq t_0$, this solution is given by*:

(2.2)
$$u(x, t) = \frac{\partial}{\partial x} \inf_{\xi(t)=x} \left\{ \mathcal{A}_{t_0, t}(\xi) + \int_0^{\xi(t_0)} u_0(z) dz \right\}$$

*and $u(\,\cdot\,, t) \in D$.*

This type of result was obtained for the first time in [Ho], [La] and [Ol] for scalar conservation laws. The generalization to multi-dimensional Hamilton-Jacobi equations is given in [Li]. Extension to the random case is straightforward, but requires some additional arguments which we present in Appendix A.

Any action minimizer $\gamma$ satisfies the following Euler-Lagrange equation:

(2.3)
$$\dot{\gamma}(s) = v(s), \quad dv(s) = \sum_{k=1}^{\infty} f_k(\gamma(s)) dB_k(s).$$

Under the assumptions in (1.3), the stochastic differential equation (2.3) has a unique solution starting at any point $x$. It is nothing but the equation of characteristics for (1.4). Therefore the variational principle (2.2) can be viewed as the generalization of the method of characteristics to weak solutions. In general, characteristics intersect each other forward in time, resulting in the formation of shocks. Given initial data at time $t_0$: $u(x, t_0) = u_0(x)$, to find the solution at $(x, t)$, consider all characteristics $\gamma$ that arrive at $x$ at time $t$ and choose among them the ones that minimize $\mathcal{A}_{t_0,t}(\gamma) + \int_0^{\gamma(t_0)} u_0(z)dz$. If such a minimizing characteristic is unique, say $\gamma(\cdot)$, then $u(x,t) = \dot{\gamma}(t)$. In the case when there are several such minimizing characteristics, $\{\gamma_\alpha(\cdot)\}$, the solution $u(\cdot, t)$ has a jump discontinuity at $x$, with $u(x-,t) = \sup_\alpha \dot{\gamma}_\alpha(t)$ and $u(x+,t) = \inf_\alpha \dot{\gamma}_\alpha(t)$.

This characterization is closely related to the notion of backward characteristics developed systematically by Dafermos (see [D]).

Our task of finding the invariant measure for (1.4) is different from what is usually asked about (1.4). Instead of solving (1.4) with given initial data, we look for a special distribution of the initial data that has the invariance property. Translated into the language of the variational principle, we will look for special minimizers or characteristics.

## 3. One-sided minimizers

A fundamental object needed for the construction of invariant measures for (1.4) is the one-sided minimizer. These are curves that minimize the action (2.1) over the semi-infinite interval $(-\infty, t]$.

In the following we will study the existence and intersection properties of one-sided minimizers. Before doing this, we formulate some basic facts concerning the effect on the action as a result of reconnecting and smoothing of curves.

*Fact* 1. Let $\xi_1, \xi_2$ be two $C^1$-curves on $[t_1, t_2]$ with values in $\mathbb{S}^1$. Then one can find a reconnection of the two curves, $\xi_r$, such that $\xi_r(t_1) = \xi_1(t_1), \xi_r(t_2) = \xi_2(t_2)$ and

$$(3.1) \qquad |\mathcal{A}^\omega_{t_1,t_2}(\xi_1) - \mathcal{A}^\omega_{t_1,t_2}(\xi_r)|, \quad |\mathcal{A}^\omega_{t_1,t_2}(\xi_2) - \mathcal{A}^\omega_{t_1,t_2}(\xi_r)|$$

$$\leq C\{\omega(\tau), \tau \in [t_1, t_2]\} \|\xi_1(t) - \xi_2(t)\|_{C^1}(1 + |t_2 - t_1|)\left(1 + \max_{t \in [t_1, t_2]}(|\dot{\xi}_1(t)|, |\dot{\xi}_2(t)|)\right).$$

Here and in the following we will use norms such as $\|\cdot\|_{C^1}$ for functions that take values on $\mathbb{S}^1$. These will always be understood as the norms of a particular representation of the functions on $\mathbb{R}^1$. The choice of the representation will either be immaterial or obvious from the context.

*Fact* 2. If $\xi$ is a curve containing corners, i.e. jump discontinuities of $\dot{\xi}$, smoothing out a corner in a sufficiently small neighborhood strictly decreases the action.

Both facts are classical and are more or less obvious.

The following lemma provides a bound on the velocities of minimizers over a large enough time interval.

LEMMA 3.1. *For almost all $\omega \in \Omega$ and any $t \in (-\infty, \infty)$ there exist random constants $T(\omega, t)$ and $C(\omega, t)$ such that if $\gamma$ minimizes $\mathcal{A}^{\omega}_{t_1, t}(\cdot)$ and $t_1 < t - T(\omega, t)$, then*

$$|\dot{\gamma}(t)| \leq C(\omega, t). \tag{3.2}$$

*Proof.* Denote

$$C_1(\omega, t) = \tfrac{1}{4} + \max_{t-1 \leq s \leq t} \sum_{k=1}^{\infty} \|F_k(x)\|_{C^2} |B_k(s) - B_k(t)| \tag{3.3}$$

and set $C(\omega, t) = 20 C_1(\omega, t)$, $T(\omega, t) = (4 C_1(\omega, t))^{-1}$. Clearly $T(\omega, t) < 1$. If $|\dot{\gamma}(t)| \leq 16 C_1$ then (3.2) is true with $C = 16 C_1$.

If $|\dot{\gamma}(t)| > 16 C_1$, we first show that the velocity $\dot{\gamma}(s)$ cannot be too large inside the interval $[t - T, t]$. Denote

$$v_0 = |\dot{\gamma}(t)| \quad \text{and} \quad v = \max_{t-T \leq s \leq t} |\dot{\gamma}(s)|. \tag{3.4}$$

Integrating by parts from (2.3), one gets for $s \in [t - T, t]$

$$\begin{aligned}
|\dot{\gamma}(s)| &= \left| \dot{\gamma}(t) - \int_s^t \sum_{k=1}^{\infty} f_k(\gamma(r)) dB_k(r) \right| \\
&\leq v_0 + \left| \sum_{k=1}^{\infty} f_k(\gamma(s))(B_k(s) - B_k(t)) \right| \\
&\quad + \left| \int_s^t \dot{\gamma}(r) \sum_{k=1}^{\infty} f'_k(\gamma(r))(B_k(r) - B_k(t)) dr \right| \\
&\leq v_0 + C_1 + C_1 v T \\
&= v_0 + C_1 + \frac{1}{4} v.
\end{aligned} \tag{3.5}$$

Hence

$$v \leq v_0 + C_1 + \frac{v}{4}, \tag{3.6}$$

implying

$$v \leq \frac{4}{3}(v_0 + C_1) \leq \frac{3}{2} v_0 \tag{3.7}$$

since $v_0 > 16 C_1$.

Next we check that $|\dot{\gamma}(s)|$ remains of order $v_0$, i.e. sufficiently large, for $s \in [t-T, t]$. As before, we have

$$
(3.8) \quad |\dot{\gamma}(s) - \dot{\gamma}(t)| = \left| \int_s^t \sum_{k=1}^{\infty} f_k(\gamma(r)) dB_k(r) \right|
$$
$$
\leq C_1 + C_1 vT
$$
$$
\leq C_1 + \frac{3}{8} v_0
$$
$$
\leq \frac{1}{2} v_0 .
$$

The last step is to show that (3.8) contradicts the minimization property of $\gamma(s)$ if $v_0 > 20C_1$. Consider a straight line $\gamma_1(s)$ joining $\gamma(t)$ and $\gamma(t-T)$. Clearly $|\gamma(t) - \gamma(t-T)| \leq 1$ since $\gamma(t), \gamma(t-T) \in \mathbb{S}^1$. Then

$$
(3.9) \quad \mathcal{A}^\omega_{t-T,t}(\gamma_1) \leq \frac{1}{2} \left( \frac{\gamma(t) - \gamma(t-T)}{T} \right)^2 T + C_1 + C_1 \left| \frac{\gamma(t) - \gamma(t-T)}{T} \right| T \leq \frac{1}{2T} + 2C_1
$$

while

$$
(3.10) \quad \mathcal{A}^\omega_{t-T,t}(\gamma) \geq \frac{1}{2} \left( \frac{v_0}{2} \right)^2 T - C_1 - \frac{3}{2} v_0 C_1 T .
$$

It is easy to see that $\frac{1}{2} \frac{v_0^2}{4} T - \frac{3}{2} v_0 C_1 T > \frac{1}{2T} + 3C_1$ for $v_0 > 20C_1$; i.e.,

$$
(3.11) \quad \mathcal{A}^\omega_{t-T,t}(\gamma_1) < \mathcal{A}^\omega_{t-T,t}(\gamma) .
$$

This contradicts the minimization property of $\gamma$. Hence $v_0 \leq 20C_1$. □

Now we are ready to prove the existence of one-sided minimizers that arrive at any given point $x \in \mathbb{S}^1$.

THEOREM 3.1. *With probablity 1, the following holds. For any $(x,t) \in \mathbb{S}^1 \times \mathbb{R}^1$, there exists at least one one-sided minimizer $\gamma \in C^1(-\infty, t]$, such that $\gamma(t) = x$.*

*Proof.* Given $\omega \in \Omega$, fix $(x,t) \in \mathbb{S}^1 \times \mathbb{R}^1$. Consider a family of minimizers $\{\gamma_\tau\}$ for $\tau < t - T(\omega, t)$, where $\gamma_\tau$ minimizes $\mathcal{A}^\omega_{\tau,t}(\xi)$ subject to the constraint that $\xi(t) = x$, $\xi(\tau) \in \mathbb{S}^1$. From Lemma 3.1, we know that $\{\dot{\gamma}_\tau(t)\}$ is uniformly bounded in $\tau$. Therefore, there exists a subsequence $\{\tau_j\}$, $\tau_j \to -\infty$, and $v \in \mathbb{R}^1$, such that

$$
\lim_{\tau_j \to -\infty} \dot{\gamma}_{\tau_j}(t) = v .
$$

Furthermore, if we define $\gamma$ to be a solution of (2.3) on $(-\infty, t]$ such that $\gamma(t) = x$, $\dot{\gamma}(t) = v$, then $\gamma_{\tau_j}$ converges to $\gamma$ uniformly, together with their derivatives, on compact subsets of $(-\infty, t]$. We will show that $\gamma$ is a one-sided minimizer.

Assume that there exists a compact perturbation $\gamma_1 \in C^1(-\infty, t]$, of $\gamma$ such that $\gamma_1(t) = x$, support $(\gamma_1 - \gamma) \subset [t_2, t_3]$, and

$$\mathcal{A}^\omega_{t_2,t_3}(\gamma) - \mathcal{A}^\omega_{t_2,t_3}(\gamma_1) = \varepsilon > 0.$$

Let $j$ be sufficiently large such that $\tau_j \leq t_2$ and

(3.12) $$|\mathcal{A}^\omega_{t_2,t}(\gamma) - \mathcal{A}^\omega_{t_2,t}(\gamma_{\tau_j})| \leq \frac{\varepsilon}{3}$$

and

(3.13) $$\|\gamma(s) - \gamma_{\tau_j}(s)\|_{C^1[t_2-1,t_2]} \leq \delta$$

($\delta$ will be chosen later). Define a new curve $\gamma_2$ by

(3.14) $$\gamma_2(s) = \begin{cases} \gamma_{\tau_j}(s), & \text{for } s \in [\tau_j, t_2 - 1]; \\ \gamma_r(s), & \text{for } s \in [t_2 - 1, t_2]; \\ \gamma_1(s), & \text{for } s \in [t_2, t], \end{cases}$$

where $\gamma_r$ is the reconnecting curve described in Fact 1. We have

(3.15) $$\begin{aligned}\mathcal{A}^\omega_{\tau_j,t}(\gamma_{\tau_j}) - \mathcal{A}^\omega_{\tau_j,t}(\gamma_2) &= \mathcal{A}^\omega_{t_2,t}(\gamma_{\tau_j}) - \mathcal{A}^\omega_{t_2,t}(\gamma) \\ &+ \mathcal{A}^\omega_{t_2,t}(\gamma) - \mathcal{A}^\omega_{t_2,t}(\gamma_1) \\ &+ \mathcal{A}^\omega_{t_2-1,t_2}(\gamma_{\tau_j}) - \mathcal{A}^\omega_{t_2-1,t_2}(\gamma_2) \\ &\geq -\frac{\varepsilon}{3} + \varepsilon - C\delta \\ &\geq \frac{\varepsilon}{3},\end{aligned}$$

if $\delta$ is small enough. Here the constant $C$ depends only on $\omega$ and $\gamma_1$. This contradicts the minimization property of $\gamma_{\tau_j}(s)$. $\square$

Now we study the intersection properties of one-sided minimizers. We use $C^1_x(-\infty, t]$ to denote the set of $C^1$ curves $\gamma$ on $(-\infty, t]$ such that $\gamma(t) = x$. We start with a general fact for minimizers (see [A], [M]).

LEMMA 3.2. *Two different one-sided minimizers* $\gamma_1 \in C^1(-\infty, t_1]$ *and* $\gamma_2 \in C^1(-\infty, t_2]$ *cannot intersect each other more than once.*

In other words, if two one-sided minimizers intersect more than once, they must coincide on their common interval of definition.

*Proof.* Suppose that $\gamma_1$ and $\gamma_2$ intersect each other twice at times $t_3$ and $t_4$, with $t_4 > t_3$. Assume without loss of generality

(3.16) $$\mathcal{A}^\omega_{t_3,t_4}(\gamma_1) \leq \mathcal{A}^\omega_{t_3,t_4}(\gamma_2).$$

Then for the curve

(3.17) $$\gamma_3(s) = \begin{cases} \gamma_2(s), & \text{for } s \in (-\infty, t_3] \cup [t_4, t_2]; \\ \gamma_1(s), & \text{for } s \in [t_3, t_4], \end{cases}$$

one has

$$\mathcal{A}^\omega_{t_3,t_4}(\gamma_3) \leq \mathcal{A}^\omega_{t_3,t_4}(\gamma_2), \tag{3.18}$$

where $\gamma_3$ has two corners at $t_3$ and $t_4$. Smoothing out these corners, we end up with a curve $\gamma^* \in C^1(-\infty, t_2]$ for which

$$\mathcal{A}^\omega_{t_3-\tau,t_2}(\gamma^*) - \mathcal{A}^\omega_{t_3-\tau,t_2}(\gamma_2) < 0 \tag{3.19}$$

for some $\tau > 0$. This contradicts the assumption that $\gamma_2(s)$ is a one-sided minimizer. □

Exploiting the random origin of the force $f$, we can prove a result which is much stronger than Lemma 3.2.

THEOREM 3.2. *The following holds for almost all $\omega$. Let $\gamma_1$, $\gamma_2$ be two distinct one-sided minimizers on the intervals $(-\infty, t_1]$ and $(-\infty, t_2]$, respectively. Assume that they intersect at the point $(x,t)$. Then $t_1 = t_2 = t$, and $\gamma_1(t_1) = \gamma_2(t_2) = x$.*

In other words, two one-sided minimizers do not intersect except for the following situation: they both come to the point $(x, t)$, having no intersections before and they both are terminated at that point as minimizers. Of course they can be continued beyond time $t$ as the solution of SDE (2.3) but they are no longer one-sided minimizers.

The proof of Theorem 3.2 resembles that of Lemma 3.2 with an additional observation that, because of the randomness of $f$, two minimizers always have an "effective intersection at $t = -\infty$." The precise formulation of this statement is given by:

LEMMA 3.3. *With probability 1, for any $\varepsilon > 0$ and any two one-sided minimizers $\gamma_1 \in C^1(-\infty, t_1]$ and $\gamma_2 \in C^1(-\infty, t_2]$, there exist a constant $T = T(\varepsilon)$ and an infinite sequence $t_n(\omega, \varepsilon) \to -\infty$ such that*

$$|\mathcal{A}^\omega_{t_n-T,t_n}(\gamma_1) - \mathcal{A}^\omega_{t_n-T,t_n}(\gamma_{1,2})|, \quad |\mathcal{A}^\omega_{t_n-T,t_n}(\gamma_2) - \mathcal{A}^\omega_{t_n-T,t_n}(\gamma_{1,2})|,$$
$$|\mathcal{A}^\omega_{t_n-T,t_n}(\gamma_1) - \mathcal{A}^\omega_{t_n-T,t_n}(\gamma_{2,1})|, \quad |\mathcal{A}^\omega_{t_n-T,t_n}(\gamma_2) - \mathcal{A}^\omega_{t_n-T,t_n}(\gamma_{2,1})| < \varepsilon, \tag{3.20}$$

*where $\gamma_{1,2}$ is the reconnecting curve defined in Fact 1 with*

$$\gamma_{1,2}(t_n - T) = \gamma_1(t_n - T), \gamma_{1,2}(t_n) = \gamma_2(t_n),$$

*and $\gamma_{2,1}$ is the reconnecting curve satisfying*

$$\gamma_{2,1}(t_n - T) = \gamma_2(t_n - T), \gamma_{2,1}(t_n) = \gamma_1(t_n).$$

*Proof.* Fix $T$ sufficiently large. With probability 1, there exists a sequence $t_n(\omega, \varepsilon) \to -\infty$ such that

$$
\text{(3.21)} \qquad \max_{s \in \bigcup_n [t_n - T, t_n]} \sum_{k=1}^{\infty} \|F_k(x)\|_{C^2} |B_k(s) - B_k(t_n)| \leq C_1 = \frac{1}{4T}.
$$

Repeating the proof of Lemma 3.1, one can check that for any $n$

$$
\text{(3.22)} \qquad \max_{t_n - T \leq s \leq t_n} (|\dot{\gamma}_1(s)|, |\dot{\gamma}_2(s)|) \leq \frac{4}{3}(20 C_1 + C_1) = \frac{7}{T}.
$$

Using (3.22), we can choose $\gamma_{1,2}$, $\gamma_{2,1}$ such that

$$
\text{(3.23)} \qquad \max_{t_n - T \leq s \leq t_n} (|\dot{\gamma}_{1,2}(s)|, |\dot{\gamma}_{2,1}(s)|) \leq \frac{7}{T} + \frac{1}{T} = \frac{8}{T}.
$$

We then have

(3.24)
$$
\begin{aligned}
&|\mathcal{A}^\omega_{t_n - T, t_n}(\gamma_1) - \mathcal{A}^\omega_{t_n - T, t_n}(\gamma_{1,2})| \\
&\leq \left| \sum_{k=1}^{\infty} (F_k(\gamma_1(t_n)) - F_k(\gamma_{1,2}(t_n)))(B_k(t_n) - B_k(t_n - T)) \right| \\
&+ \int_{t_n - T}^{t_n} \left| \left(\frac{1}{2}\dot{\gamma}_1(t)^2 - \frac{1}{2}\dot{\gamma}_{1,2}(t)^2\right) \right. \\
&- \sum_{k=1}^{\infty} (B_k(t) - B_k(t_n - T))\Big(f_k(\gamma_1(t))(\dot{\gamma}_1(t) - \dot{\gamma}_{1,2}(t)) \\
&+ \left. (f_k(\gamma_1(t)) - f_k(\gamma_{1,2}(t)))\dot{\gamma}_{1,2}(t)\Big) \right| dt \\
&\leq \frac{1}{4T} + T\left(\frac{1}{2}\left(\frac{7}{T}\right)^2 + \frac{1}{2}\left(\frac{8}{T}\right)^2 + C_1\left(\frac{7}{T} + \frac{8}{T}\right) + C_1 \frac{8}{T}\right) \\
&= \frac{125}{2T} \\
&\leq \varepsilon,
\end{aligned}
$$

if $T \geq \frac{125}{2\varepsilon}$. Similarly, one proves other inequalities in (3.20). $\square$

*Proof of Theorem* 3.2. We will use $\tau$ to denote a sufficiently large negative number. Suppose that $\gamma_1$ and $\gamma_2$ intersect each other at time $t < \max(t_1, t_2)$ and for definiteness let $t_1 > t$. Then the curve

$$
\text{(3.25)} \qquad \gamma_3(s) = \begin{cases} \gamma_2(s), & \text{for } s \in (-\infty, t]; \\ \gamma_1(s), & \text{for } s \in [t, t_1] \end{cases}
$$

has a corner at time $t$. This corner can be smoothed out according to Fact 2, and the resulting curve $\gamma^* \in C^1(-\infty, t_1]$ satisfies

$$
\text{(3.26)} \qquad \mathcal{A}^\omega_{\tau, t_1}(\gamma_3) - \mathcal{A}^\omega_{\tau, t_1}(\gamma^*) = \delta > 0.
$$

Set $\varepsilon = \delta/4$. Choose sufficiently negative $t_n(\omega, \varepsilon)$ defined in Lemma 3.3 such that $\gamma^*(s) = \gamma_2(s)$ for $s \in (-\infty, t_n]$.

Assume that

$$\mathcal{A}^\omega_{t_n,t}(\gamma_2) - \mathcal{A}^\omega_{t_n,t}(\gamma_1) > 2\varepsilon. \qquad (3.27)$$

Then in view of Lemma 3.3

$$\gamma_4(s) = \begin{cases} \gamma_2(s), & \text{for } s \in (-\infty, t_n - T]; \\ \gamma_{2,1}(s), & \text{for } s \in [t_n - T, t_n]; \\ \gamma_1(s), & \text{for } s \in [t_n, t], \end{cases} \qquad (3.28)$$

is a local perturbation of $\gamma_2 \in C^1(-\infty, t]$ with

$$\begin{aligned}
\mathcal{A}^\omega_{\tau,t}(\gamma_2) - \mathcal{A}^\omega_{\tau,t}(\gamma_4) &= \mathcal{A}^\omega_{t_n-T,t_n}(\gamma_2) - \mathcal{A}^\omega_{t_n-T,t_n}(\gamma_{2,1}) \\
&\quad + \mathcal{A}^\omega_{t_n,t}(\gamma_2) - \mathcal{A}^\omega_{t_n,t}(\gamma_1) \\
&> -\varepsilon + 2\varepsilon \\
&= \varepsilon.
\end{aligned} \qquad (3.29)$$

This contradicts the assumption that $\gamma_2$ is a one-sided minimizer. Thus

$$\mathcal{A}^\omega_{t_n,t}(\gamma_1) - \mathcal{A}^\omega_{t_n,t}(\gamma_2) \geq -2\varepsilon, \qquad (3.30)$$

and

$$\gamma_5(s) = \begin{cases} \gamma_1(s), & \text{for } s \in (-\infty, t_n - T]; \\ \gamma_{1,2}(s), & \text{for } s \in [t_n - T, t_n]; \\ \gamma^*(s), & \text{for } s \in [t_n, t_1] \end{cases} \qquad (3.31)$$

is a local perturbation of $\gamma_1 \in C^1(-\infty, t_1]$ with

$$\begin{aligned}
\mathcal{A}^\omega_{\tau,t_1}(\gamma_1) - \mathcal{A}^\omega_{\tau,t_1}(\gamma_5) &= \mathcal{A}^\omega_{t_n-T,t_n}(\gamma_1) - \mathcal{A}^\omega_{t_n-T,t_n}(\gamma_{1,2}) \\
&\quad + \mathcal{A}^\omega_{t_n,t}(\gamma_1) - \mathcal{A}^\omega_{t_n,t}(\gamma_2) \\
&\quad + \mathcal{A}^\omega_{\tau,t_1}(\gamma_3) - \mathcal{A}^\omega_{\tau,t_1}(\gamma^*) \\
&\geq -\varepsilon - 2\varepsilon + \delta \\
&= \varepsilon > 0.
\end{aligned} \qquad (3.32)$$

This contradicts the assumption that $\gamma_1$ is a one-sided minimizer and proves the theorem. □

Theorem 3.2 implies the following remarkable properties of one-sided minimizers. Given $\omega$ and $t$, denote by $J(\omega, t)$ the set of points $x \in \mathbb{S}^1$ with more than one one-sided minimizer coming to $(x, t)$.

LEMMA 3.4. *The following holds with probability* 1. *For any $t$, the set $J(\omega, t)$ is at most countable.*

*Proof.* Any $x \in J(\omega, t)$ corresponds to a segment $[\gamma_-(t-1), \gamma_+(t-1)]$, where $\gamma_-$ and $\gamma_+$ are two different one-sided minimizers coming to $(x, t)$ and $\gamma_+(s) > \gamma_-(s)$, for $s < t$. In view of Theorem 3.2, these segments are mutually disjoint. This implies the lemma. □

LEMMA 3.5. *Given $\omega$ and $t$, consider a sequence of one-sided minimizers $\gamma_n(s)$ defined on $(-\infty, t]$ such that $\gamma_n(t) \to x$ and $\dot{\gamma}_n(t) \to v$ as $n \to \infty$. Let $\gamma$ be the solution of the SDE (2.3) on $(-\infty, t]$ with the initial data $\gamma(t) = x$ and $\dot{\gamma}(t) = v$. Then $\gamma$ is a one-sided minimizer.*

*Proof.* Suppose that $\gamma^* \in C^1(-\infty, t]$ coincides with $\gamma$ outside an interval $[t_1, t_2] \subset (-\infty, t]$ and $\mathcal{A}^\omega_{t_1, t_2}(\gamma) - \mathcal{A}^\omega_{t_1, t_2}(\gamma^*) = \varepsilon > 0$. It is clear that by taking sufficiently large $n$ one can make $\|\gamma(s) - \gamma_n(s)\|_{C^1[t_1-1,t]}$ arbitrarily small. Let $\gamma_1$ be the reconnecting curve on $[t_1 - 1, t_1]$ between $\gamma_n(t_1 - 1)$ and $\gamma(t_1)$, and for some $\delta > 0$, let $\gamma_2$ be the reconnecting curve on $[t - \delta, t]$ between $\gamma(t - \delta)$ and $\gamma_n(t)$. Then the curve

$$(3.33) \qquad \gamma^{**}(s) = \begin{cases} \gamma_n(s), & \text{for } s \in (-\infty, t_1 - 1]; \\ \gamma_1(s), & \text{for } s \in [t_1 - 1, t_1]; \\ \gamma^*(s), & \text{for } s \in [t_1, t_2]; \\ \gamma(s), & \text{for } s \in [t_2, t - \delta]; \\ \gamma_2(s), & \text{for } s \in [t - \delta, t] \end{cases}$$

satisfies $\mathcal{A}^\omega_{-\infty,t}(\gamma_n) - \mathcal{A}^\omega_{-\infty,t}(\gamma^{**}) > 0$ if $\delta$ and $\|\gamma(s) - \gamma_n(s)\|_{C^1[t_1-1,t]}$ are small enough. This contradicts the assumption that $\gamma_n$ is a one-sided minimizer since $\gamma^{**}$ is a local perturbation of $\gamma_n$. Note that (3.33) cannot be used if $t_2 = t$. In this case in the segment $[t - \delta, t]$ one can directly reconnect $\gamma_n$ and $\gamma^*$ and it is not hard to check that for $\delta$ small enough, $|\mathcal{A}^\omega_{t-\delta,t}(\gamma^*) - \mathcal{A}^\omega_{t-\delta,t}(\gamma^{**})|$ can be made arbitrarily small. □

LEMMA 3.6. *With probability one, the following holds. Fix an arbitrary sequence $t_n \to -\infty$ and a sequence of functions $\{v_n\}$, $v_n \in D_0$, $\int_0^1 v_n(z)dz = 0$. Consider (1.4) on the time interval $[t_n, t]$ with the initial condition $u(x, t_n) = v_n(x)$. Take any $x \in \mathbb{S}^1$ and a sequence of characteristics $\gamma_n \in C^1[t_n, t]$, $\gamma_n(t) = x$ minimizing $\mathcal{A}^\omega_{t_n, t}(\xi) + \int_0^{\xi(t_n)} v_n(z)dz$. Suppose that $v$ is a limiting point of the set $\{\dot{\gamma}_n(t)\}$. Then the solution $\gamma$ of SDE (2.3) with initial data $\gamma(t) = x$ and $\dot{\gamma}(t) = v$ is a one-sided minimizer on $(-\infty, t]$.*

*Proof.* The proof of this lemma is the same as the final part of the proof of Theorem 3.1. □

Next we study the measurability issues. Fix a time $t$ and consider all integer times $-n \leq t$. Introduce

$$(3.34) \qquad A^\omega_{-n,t}(x) = \min_{\substack{\xi \in C^1[-n,t] \\ \xi(t)=x}} \mathcal{A}^\omega_{-n,t}(\xi).$$

LEMMA 3.7. *The following statement holds with probability* 1. *Suppose that $\gamma \in C^1_x(-\infty, t]$ is a one-sided minimizer. Then for any $\varepsilon > 0$ there exist an infinite number of integer times $-n \leq t$ such that*

$$(3.35) \qquad |\mathcal{A}^\omega_{-n,t}(\gamma) - A^\omega_{-n,t}(x)| \leq \varepsilon.$$

*Conversely, if a curve $\xi \in C^1_x(-\infty, t]$ has the property that for any $\varepsilon > 0$ there exist an infinite number of integer times $-n \leq t$ such that*

$$(3.36) \qquad |\mathcal{A}^\omega_{-n,t}(\xi) - A^\omega_{-n,t}(x)| \leq \varepsilon,$$

*then $\xi$ is a one-sided minimizer.*

*Proof.* Suppose that for some $\varepsilon > 0$ and $n_0$

$$(3.37) \qquad |\mathcal{A}^\omega_{-n,t}(\gamma) - A^\omega_{-n,t}(x)| > \varepsilon$$

for all $-n \leq -n_0$. Consider the curves $\xi_{-n} \in C^1_x[-n,t]$ such that $\mathcal{A}^\omega_{-n,t}(\xi_{-n}) = A^\omega_{-n,t}(x)$. Then, according to Lemma 3.3, there exist an interval $[-n_1, -n_2] \subset (-\infty, -n_0]$ and a reconnecting curve $\gamma_r$ with $\gamma_r(-n_1) = \gamma(-n_1), \gamma_r(-n_2) = \xi_{-n_1}(-n_2)$, such that

$$(3.38) \qquad |\mathcal{A}^\omega_{-n_1,-n_2}(\gamma_{-n_1}) - \mathcal{A}^\omega_{-n_1,-n_2}(\gamma_r)| \leq \frac{\varepsilon}{2}.$$

Then

$$(3.39) \qquad \gamma_1(s) = \begin{cases} \gamma(s), & \text{for } s \in (-\infty, -n_1]; \\ \gamma_r(s), & \text{for } s \in [-n_1, -n_2]; \\ \xi_{-n_1}(s), & \text{for } s \in [-n_2, t] \end{cases}$$

is a local perturbation of $\gamma$ which lowers the action by at least $\varepsilon/2$. This contradicts the assumption that $\gamma$ is a one-sided minimizer.

Note that formally Lemma 3.3 cannot be applied here since $\xi_{-n_1}$ is not a one-sided minimizer but Lemma 3.1 remains valid for all $\xi_{-n}$ with sufficiently negative $-n$. Thus the same argument as in the proof of Lemma 3.3 proves (3.38).

To prove the second statement, observe that if $\xi_1$ is a local perturbation of $\xi$ lowering $\mathcal{A}^\omega_{-n,t}(\xi)$ by some $\varepsilon > 0$, then $\mathcal{A}^\omega_{-n,t}(\xi) \geq A^\omega_{-n,t}(x) + \varepsilon$ for all sufficiently negative $-n$. This contradicts (3.36). □

Now we are ready to define the main object of this paper. We will denote by $\{\gamma_{x,t,\alpha}(s)\}$ the family of all one-sided minimizers coming to $(x,t)$ indexing them by $\alpha$.

*Definition* 3.1.

$$u_+^\omega(x,t) = \inf_\alpha \dot\gamma_{x,t,\alpha}(t), \tag{3.40}$$

$$u_-^\omega(x,t) = \sup_\alpha \dot\gamma_{x,t,\alpha}(t). \tag{3.41}$$

It is clear that $u_+^\omega(x,t) = u_-^\omega(x,t)$ for $x \notin J(\omega, t)$.

LEMMA 3.8. *With probability 1, for every $x \in \mathbb{S}^1$*

$$\lim_{y \uparrow x} u_+^\omega(y,t) = u_-^\omega(x,t), \tag{3.42}$$

$$\lim_{y \downarrow x} u_+^\omega(y,t) = u_+^\omega(x,t), \tag{3.43}$$

*and hence $u_+^\omega(\cdot, t) \in D$ for fixed $t$.*

*Proof.* We will prove (3.42). The proof of (3.43) is similar. It was shown in Lemma 3.1 that $|u_+^\omega(y,t)| \leq C(\omega, t)$. Suppose that there exists a sequence $y_n \uparrow x$ such that $u_+^\omega(y_n, t) \to v \neq u_-^\omega(x,t)$. Then, according to Lemma 3.5, the solution $\gamma$ of SDE (2.3) with the initial data $\gamma(t) = x$ and $\dot\gamma(t) = v$ is a one-sided minimizer. Theorem 3.2 implies that $\dot\gamma(t) > u_-^\omega(x,t)$ which contradicts the definition of $u_-^\omega(x,t)$. □

It follows immediately from the construction that on any finite time interval $[t_1, t_2]$, $u_+^\omega$ is a weak solution of (1.4) with initial data $u_0(x) = u_+^\omega(x, t_1)$. Moreover, the following statement holds:

LEMMA 3.9. *Given $t$, the mapping $u_+^\omega(\cdot, t): \Omega \mapsto D$ is measurable.*

*Proof.* Without loss of generality, let $t = 0$. Since $\mathcal{D}$ is generated by cylinder sets of the type $A(x_1, \ldots, x_n)$ with $x_i$ from a dense subset of $\mathbb{S}^1$, it is enough to show that $u_+^\omega(x, 0): \Omega \to \mathbb{R}^1$ are measurable for a dense set of $x$ values. For any positive integer $n$, denote by $u_{-n,+}^\omega$ the right continuous weak solution of (1.4) on the time interval $[-n, 0]$ with the initial data $u_{-n,+}^\omega(x, -n) \equiv 0$. For any $x \in \mathbb{S}^1$ and $v \in \mathbb{R}^1$ denote by $\xi_{x,v}^\omega(s)$, $s \in [-n, 0]$ the backward solution of (2.3) with the initial data $\xi_{x,v}^\omega(0) = x$ and $\dot\xi_{x,v}^\omega(0) = v$. The action $A_{-n,0}^\omega(x,v) = \mathcal{A}_{-n,0}^\omega(\xi_{x,v}^\omega)$ is a continuous function on $\Omega \times \mathbb{S}^1 \times \mathbb{R}^1$. Hence the set $M = \{(\omega, x, v): A_{-n,0}^\omega(x,v) = A_{-n,0}^\omega(x)\}$ is closed. Let $M_{\omega,x} = \{v \in \mathbb{R}^1: (\omega, x, v) \in M\}$. We conclude that $u_{-n,+}^\omega(x, 0) = \max_v M_{\omega,x}$ is a measurable function on $\Omega \times \mathbb{S}^1$ and $u_{-n,+}^\omega(\cdot, 0)$ is a measurable mapping $\Omega \mapsto D$.

As in the proof of Theorem 3.1, it is easy to check that for

$$v_+^\omega(x,0) = \limsup_n u_{-n,+}^\omega(x,0) \tag{3.44}$$

the corresponding curve $\gamma_x^\omega = \xi_{x,v_+^\omega(x,0)}^\omega$ is a one-sided minimizer.

For positive integers $k$ and $m$, introduce measurable subsets of $\Omega \times \mathbb{S}^1 \times \mathbb{R}^1$:

$$M_n(k,m) = \big\{(\omega,x,v) \in \Omega \times \mathbb{S}^1 \times \mathbb{R}^1 \colon |v - u_{-n,+}^\omega(x,0)| \geq \frac{1}{k}, \tag{3.45}$$

$$\mathcal{A}_{-n,0}^\omega(x,v) - A_{-n,0}^\omega(x) \leq \frac{1}{m}\big\}.$$

Let $\Pi$ be the projection of the measurable set

$$\cup_k \cap_m \big(\cap_l \cup_{n=l}^\infty M_n(k,m)\big) \tag{3.46}$$

on $\Omega \times \mathbb{S}^1$. The points $(\omega, x)$ in $\Pi$ are characterized by the following property. There exists a backward solution $\xi_{x,v}^\omega$ different from the one-sided minimizer $\gamma_x^\omega$ such that at infinitely many negative integer times $-n$, the action $\mathcal{A}_{-n,0}^\omega(\xi_{x,v}^\omega)$ is arbitrarily close to its minimal value $A_{-n,0}^\omega(x)$. In view of Lemma 3.7, $\Pi$ is precisely the set $(\omega, x)$ having at least two one-sided minimizers coming to $(x,0)$; i.e.,

$$\Pi = J = \{\omega, x) \in \Omega \times \mathbb{S}^1 \colon u_+^\omega(x,0) \neq u_-^\omega(x,0)\}. \tag{3.47}$$

Consider the sections $J_x = \{\omega \in \Omega \colon u_+^\omega(x,0) \neq u_-^\omega(x,0)\}$ of $J$ and define $I = \{x \in \mathbb{S}^1 \colon P(J_x) > 0\}$. Using the measurability of $J$, Lemma 3.4 and Fubini's theorem, we conclude that the Lebesgue measure of $I$ is 0.

Fix $x \in \mathbb{S}^1 \setminus I$. Then for almost all $\omega$ one has

$$u_+^\omega(x,0) = \lim_{n \to \infty} u_{n,+}^\omega(x,0). \tag{3.48}$$

The functions $u_{n,+}^\omega(x,0)\colon \Omega \mapsto \mathbb{R}^1$ are measurable. Hence $u_+^\omega(x,0)\colon \Omega \mapsto \mathbb{R}^1$ is also measurable. This proves the lemma, since $\mathbb{S}^1 \setminus I$ is dense in $\mathbb{S}^1$. □

Before ending this section, we formulate some corollaries of Lemmas 3.1 and 3.6 that will be useful later.

LEMMA 3.10. *The following estimate holds*:

$$\|u^\omega(\,\cdot\,,t)\|_{L^\infty(\mathbb{S}^1)} \leq \bar{C}(\{\omega(s), s \in [t-1,t)\}).$$

*The stationary random variable $\bar{C}(\{\omega(s),\ s \in [t-1,t]\})$ has finite moments of all orders.*

This estimate says that $|u^\omega(x,t)|$ can be bounded by a quantity that depends only on $\{\omega(s), s \in [t-1, t)\}$.

LEMMA 3.11. *Let $\{v_n\}$ be a sequence of functions in $D_0$. Let $\{t_n\}$ be a sequence such that $t_n \to -\infty$, and let $u_n$ be the solution of (1.4) with initial condition $u_n(x, t_n) = v_n(x)$. Then for almost all $\omega$, $\lim_{t_n \to -\infty} u_n(\,\cdot\,, 0) = u_+^\omega(\,\cdot\,, 0)$ almost everywhere.*

This result follows directly from Lemma 3.6.

So far we restricted our attention to the case when $\int_0^1 u(x,t)dx = 0$. This is reflected in the variational principle used in (2.1). In the case when $\int_0^1 u(x,t)dx = c \neq 0$, the relevant action is replaced by

$$\mathcal{A}_{t_1,t_2}(\xi) = \int_{t_1}^{t_2} \left\{ \tfrac{1}{2}(\dot\xi(s) - c)^2 - \sum_k f_k(\xi(s))\dot\xi(s)(B_k(s) - B_k(t_1)) \right\} ds$$
$$+ \sum_k F_k(\xi(t_2))(B_k(t_2) - B_k(t_1)).$$

In this more general case, we can also define one-sided minimizers in an analogous way, and all the results obtained so far hold in the general case. Instead of giving all the details, let us just comment on the most important aspect when $c \neq 0$. When $c = 0$, the one-sided minimizers stay roughly inside one period so that their asymptotic speed (which is the analog of rotation number in Aubry-Mather theory) is zero when lifted to the universal cover:

$$\alpha(c = 0) = \lim_{t \to -\infty} \frac{\xi(t)}{t} = 0.$$

In the general case,

(3.49) $$\alpha(c) = \lim_{t \to -\infty} \frac{\xi(t)}{t} = c.$$

Roughly speaking, this is because, if a curve has different asymptotic speed, the cost to the action grows linearly in time, whereas the savings from the random potential can at most grow as $O(\sqrt{t})$. For example, in Lemma 3.3 and Theorem 3.2, we showed that there are large time intervals on which the one-sided minimizers are almost parallel. These results are still true except that the minimizers are parallel with an average slope $1/c$ in the $x - t$ plane. Similarly, the reconnection used in Lemma 3.1 to prove a bound for the velocity of minimizers will also have to be done with curves that have an average slope of $1/c$.

(3.49) is in sharp contrast with the case of periodic (in $x$ and $t$) potential studied in Aubry-Mather theory. There the function $\alpha$ is usually a very complicated function of $c$ with Cantor-like structures.

## 4. Construction and uniqueness of the invariant measure

In this section we prove Theorems 1.2 and 1.3 stated in Section 1. Let $\delta^\omega(du)$ be the atomic measure on $(D_0, \mathcal{D})$ concentrated at $u_+^\omega(\cdot, 0)$.

THEOREM 4.1. *The measure*

$$\mu(d\omega, du) = \delta^\omega(du) P(d\omega) \tag{4.1}$$

*is an invariant measure for the skew-product transformation $F^t$ (see Definition 1.2), and the measure $\kappa(du) = \int_\Omega \mu(d\omega, du)$ is a stationary distribution for the Markov process corresponding to (1.4).*

*Proof.* The second statement of the theorem follows trivially from the first one since the $\delta^\omega(du)$'s are measurable with respect to the $\sigma$-algebra $\mathcal{F}^0_{-\infty}$. It is a general fact that any measure which satisfies the measurability property and which is invariant under $F^t$ generates a stationary distribution for the Markov process (1.4). The first statement is an immediate consequence of the construction of $u_+^\omega(x, t)$. Indeed, $u_+^{\theta^t\omega}(x, 0) = u_+^\omega(x, t)$ by construction. Hence $S_\omega^t \delta^\omega = \delta^{\theta^t\omega}$. This is exactly the condition for the invariance of $\mu$. □

The "one force, one solution" principle not only gives the existence of an invariant measure, it also implies uniqueness.

THEOREM 4.2. *The measure $\mu(d\omega, du)$ is the unique invariant measure on $(\Omega \times D_0, \mathcal{F} \times \mathcal{D})$ with given projection $P(d\omega)$ on $(\Omega, \mathcal{F})$. The measure $\kappa(du) = \int_\Omega \mu(d\omega, du)$ is the unique stationary distribution for the Markov process (1.4).*

*Proof.* Assume $\lambda$ is another invariant measure on $\Omega \times D_0$. Write $\lambda$ as

$$\lambda(d\omega, du) = \int_\Omega \lambda_\omega(du) P(d\omega). \tag{4.2}$$

For $t < 0$, let $H_t^\omega$ be the operator which maps the solution of (1.4) at time $t$ to the solution at time $0$ when the realization of the force is $\omega$. By definition, the invariance of $\lambda$ implies that

$$(H_t^\omega)^* \lambda_{\theta^t \omega}(du) = \lambda_\omega(du) \tag{4.3}$$

for $t < 0$, where $(H_t^\omega)^*$ is the push forward action on the spaces of measures. This means that there exists a subset $B$ of $D_0$ of full measure with respect to $\lambda_\omega(du)$, such that for every $u \in B$ and $n \in \mathbb{N}$, there exists a $v_n$ such that

$$H_{-n}^\omega v_n = u.$$

From Lemma 3.6, if a solution of (1.4) can be extended backward to arbitrary negative times, that solution must coincide with $u_+^\omega$ at $t = 0$, for all $x \in \mathbb{S}^1 \setminus J(\omega, 0)$. In particular, we have

$$u(x) = u_+^\omega(x, 0)$$

for $x \in \mathbb{S}^1 \setminus J(\omega, 0)$. Hence $\lambda_\omega(du) = \delta_{u_+^\omega}(du)$ and $\lambda = \mu$.

To prove the second statement suppose that $\nu(du) \neq \kappa(du)$ is another stationary distribution for the Markov process (1.4). Let

$$A = \{u\colon (u(x_1), \ldots, u(x_k)) \in C \subset \mathbb{R}^k\}$$

be an arbitrary cylindrical set based on the points $x_1, \ldots, x_k \in I^c$. By definition

$$\begin{aligned}
(4.4) \quad \nu(A) &= \int_{D_0} P_n(u_{-n}, A)\nu(du_{-n}) \\
&= \int_{D_0} \left(\int_\Omega \chi_A(u_{-n}, \omega) P(d\omega)\right) \nu(du_{-n}) \\
&= \int_\Omega \left(\int_{D_0} \chi_A(u_{-n}, \omega) \nu(du_{-n})\right) P(d\omega).
\end{aligned}$$

Denote by $\chi_A(\omega)$ the indicator function of the event that $u_+^\omega \in A$. Then for $\eta_n(\omega) = \int_{D_0} \chi_A(u_{-n}, \omega) \nu(du_{-n})$ one has

$$(4.5) \qquad \lim_{n \to \infty} \eta_n(\omega) = \chi_A(\omega)$$

implying $\nu(A) = \int_\Omega \chi_A(\omega) P(d\omega) = \kappa(A)$.

Indeed, in view of the uniqueness of one-sided minimizers, coming to each of the points $(x_1, 0), \ldots, (x_k, 0)$, one observes that

$$(4.6) \qquad \lim_{n \to \infty} \eta_n(\omega) = \begin{cases} 1 & \text{for } (A \setminus \partial A) \ni (u_+^\omega(x_1, 0), \ldots, u_+^\omega(x_k, 0)), \\ 0 & \text{for } (A^c \setminus \partial A) \ni (u_+^\omega(x_1, 0), \ldots, u_+^\omega(x_k, 0)). \end{cases}$$

Thus to obtain (4.3) one need only show that $P\{\omega\colon (u_+^\omega(x_1, 0), \ldots, u_+^\omega(x_k, 0)) \in \partial A\} = 0$. Clearly it is enough to check that:

LEMMA 4.1. *For any $v \in \mathbb{R}^1$ and for any $x \in X$,*

$$(4.7) \qquad P\{\omega\colon u_+^\omega(x, 0) = v\} = 0.$$

*Proof.* For fixed $x$ and $v$ the backward solution (for $t \leq 0$) of the SDE (2.3) with the initial data $\gamma(0) = x$ and $\dot\gamma(0) = v$ is a random process with variance $\operatorname{Var} \dot\gamma(t)$ that grows like $|t|$ as $t \to -\infty$. Thus with probability 1 this solution cannot be a one-sided minimizer. This completes the proof of Theorem 4.2. □

## 5. The two-sided minimizer and the main shock

Sections 5–7 will be devoted to the study of the structure of the solution $u^\omega$. In this section, we will define the two basic objects needed for this study, the two-sided minimizer and the main shock.

Let $(x, t)$ be a point of shock, i.e., $x \in J(\omega, t)$. Denote by $\Delta_{x,t}(t_1)$, $t_1 < t$, an open interval at time $t_1$ generated by the shock at $(x, t)$:

(5.1) $$\Delta_{x,t}(t_1) = (\gamma^-_{x,t}(t_1),\ \gamma^+_{x,t}(t_1))$$

where we use $\gamma^-_{x,t}$ and $\gamma^+_{x,t}$ to denote respectively the left-most and right-most one-sided minimizers starting at $(x,t)$. Roughly speaking, $\Delta_{x,t}(t_1)$ is the set of points at time $t_1$ that merge into the shock at $x$ before time $t$, i.e. the one-sided minimizer that passes through $(y, t_1), y \in \Delta_{x,t}(t_1)$, intersects the past history of the shock at $(x,t)$ before time $t$.

LEMMA 5.1. *For almost all $\omega$ the following statements hold for any fixed $t_1$ and $t$ such that $t_1 < t$.*

(a) *For arbitrary $x_1 \in \mathbb{S}^1$, either there exists a unique one-sided minimizer at time $t$ which passes through $(x_1, t_1)$ or there exists a unique shock at $(x,t)$ for some $x \in \mathbb{S}^1$, such that $(x_1, t_1) \in \Delta_{x,t}(t_1)$. In the second case we will say that $(x_1, t_1)$ is covered by the shock at $(x,t)$. In particular, if $(x_1, t_1)$ is a point of shock, i.e., $x_1 \in J(\omega, t_1)$, then there exists a unique shock $(x,t)$ which covers $(x_1, t_1)$; i.e., $(x_1, t_1) \in \Delta_{x,t}(t_1)$.*

(b) *Let $(x_1, t_1) \in J(\omega, t_1)$ be a point of shock. Denote by $x(t)$, $t \geq t_1$, the position of the shock at time $t$ which covers $(x_1, t_1)$. Then $x(t)$, $t \geq t_1$, is a Lipschitz continuous function.*

*Proof.* Denote by $A(t_1, t)$ the set of points $x_1 \in \mathbb{S}^1$ that correspond to the first situation; i.e., there exists a one-sided minimizer at time $t$ which passes through $(x_1, t_1)$. Obviously this minimizer is unique, and $A(t_1, t)$ is closed. Hence $B(t_1, t) = \mathbb{S}^1 \setminus A(t_1, t)$ is open and consists of nonintersecting open intervals. Let $(x'_1, x''_1)$ be one of these intervals. Both $(x'_1, t_1)$ and $(x''_1, t_1)$ are reached by one-sided minimizers which start at $(x', t)$ and $(x'', t)$. It is easy to see that $x' = x''$. Otherwise a minimizer which starts from some point between $x'$ and $x''$ will reach a point inside $(x'_1, x''_1)$. It follows that $(x', t)$ is a point of shock, and $(x'_1, x''_1) \subset \Delta_{x',t}(t_1)$. Obviously a point of shock cannot be reached by a one-sided minimizer that extends to time $t$. Thus if $t_1 < t$, any point of shock at $(x_1, t_1)$ must be covered by a shock at $(x, t)$ for some $x \in \mathbb{S}^1$. Clearly such a covering shock $(x, t)$ is unique. This completes the proof of (a).

(b) basically follows from the fact that the velocities of minimizers are bounded. It is enough to show that $x(\cdot)$ is Lipschitz continuous at $t = t_1$. It follows from Lemma 3.1 that there exists a constant $C_1(t_1, \omega)$ such that for all one-sided minimizers $\gamma$ and $t \in [t_1, t_1 + 1]$, $|\dot\gamma(t)| \leq C_1(t_1, \omega)$. Therefore for any shock at $(x(t), t)$, $t \in [t_1, t_1 + 1]$, $|x(t) - \gamma^-_{x(t),t}(t_1)| \leq C_1(t_1, \omega)(t - t_1)$, $|x(t) - \gamma^+_{x(t),t}(t_1)| \leq C_1(t_1, \omega)(t - t_1)$. Since $x_1 \in (\gamma^-_{x(t),t}(t_1), \gamma^+_{x(t),t}(t_1))$, we also have $|x_1 - x(t)| \leq C_1(t_1, \omega)(t - t_1)$. This estimate implies that $x(t)$ is Lipschitz continuous. $\square$

*Remark.* It may happen for some points $(x_1, t_1)$ that they are covered by shocks and at the same time there exist one-sided minimizers passing through them. This is possible only if there are more than two minimizers starting at covering shock. This situation occurs when two shocks merge.

In this and following sections we study the detailed structure and regularity properties of solutions supported by the invariant measure. For this purpose, we need certain nondegeneracy conditions on the forcing.

*Nondegeneracy condition.* If $x^*$ is the local maximum of some $F_k$, we will denote by $I(x^*)$ a closed interval on $\mathbb{S}^1$ containing $x^*$ which is contained in the basin of attraction of $x^*$ for the potential $F_k$. In other words, $f_k < 0$ on $I(x^*)$ to the right of $x^*$, and $f_k > 0$ on $I(x^*)$ to the left of $x^*$. Assume that

(A1) There exists a finite set of points $X^* = \{x^*\}$, $x^* \in \mathbb{S}^1$, each of which is a local maximum of some $F_k$ with the following property: for any $(x_1, x_2) \in \mathbb{S}^1 \times \mathbb{S}^1$ there exists an $x^* \in X^*$, such that $x_1, x_2 \in I(x^*)$.

Below we will always assume that (A1) holds. Obviously (A1) fails if there is only one term in the sum of $F(x,t)$. Nevertheless (A1) is fulfilled in generic situations. In particular, it is easy to see that three such intervals $I(x^*)$ suffice. However, by refining the argument in Appendix D, we can show that the basic collision lemma below also holds when the potential contains two shifted cosine functions, for example. We will come back to this point at the end of this section.

Consider two points $x_1^0$ and $x_2^0$ at time $t = 0$. We say that $x_1^0$ and $x_2^0$ merge before $t = \tau > 0$ if there exists a shock at $(y, \tau)$, $y \in \mathbb{S}^1$, which covers both $x_1^0$ and $x_2^0$; i.e., $x_1^0, x_2^0 \in \Delta_{y,\tau}(0)$. The following lemma is of fundamental importance to what follows.

LEMMA 5.2 (basic collision lemma). *For any $\tau > 0$, there exists a positive number $p_0(\tau)$ with the following property. Let $u(\,\cdot\,, 0) \in D_0$ and $x_1^0, x_2^0$ be two positions at $t = 0$ which are measurable with respect to $\mathcal{F}_{-\infty}^0$. Then the conditional probability under $\mathcal{F}_{-\infty}^0$ that $x_1^0$ merges with $x_2^0$ before $t = \tau$ is no less than $p_0(\tau)$.*

The proof of this lemma is given in Appendix D.

LEMMA 5.3. *The set of $\omega$'s for which $u^\omega(\,\cdot\,, t_0)$ is continuous for some $t_0 \in \mathbb{R}^1$ has probability zero.*

*Proof.* It follows from Lemma 5.1 that if $u^\omega(\,\cdot\,, t_0)$ is continuous, then $u^\omega(\,\cdot\,, t)$ is continuous for all $t \leq t_0$. Denote by $C(t)$ the set of $\omega$ such that $u^\omega(\,\cdot\,, s)$ is continuous for $s \leq t$. Then $\theta^s C(t) = C(t+s) \subseteq C(t)$ for all $s \geq 0$. Using ergodicity we conclude that either $P(C(t)) = 1$ for all $t$, or $P(C(t)) = 0$ for all $t$. Assume that $P(C(t)) \equiv 1$. It follows from Lemma 5.2 that there is a positive conditional probability $p(\{\omega(t'), t' < t\})$ under $\mathcal{F}^t_{-\infty}$ that $u(\,\cdot\,, t+1)$ has at least one shock. Hence we have

$$(5.2) \quad P(C(t+1)) = \int_{C(t)} \left(1 - p\{\omega(t'), t' \leq t\}\right) dP < P(C(t)).$$

Therefore $P(C(t+1)) < 1$. Hence $P(C(t)) \equiv 0$ for all $t$. □

It follows from Lemma 5.3 that for almost all $\omega$ and arbitrary $t_0$ there exists at least one shock at $t_0$. Since the number of shocks is at most countable and the sum of their sizes is bounded, i.e.,

$$\sum_{x \in J(\omega, t_0)} (u^\omega_-(x, t_0) - u^\omega_+(x, t_0)) \leq \operatorname{Var} u^\omega_+(x, t_0) < +\infty,$$

we can numerate all the shocks in a measurable way. The first shock is the largest one, the second is next in size and so on. If two or more shocks have the same size then we numerate them according to their order of occurrence on the semi-interval $[0, 1)$. Denote by $\xi_i(t_0) = (x_i(t_0), t_0)$ the position of the $i$-th shock. Obviously $\xi_i(t_0)$ is a measurable function with respect to $\mathcal{F}^{t_0}_{-\infty}$.

We will use $\ell(I)$ to denote the length of the interval $I$.

LEMMA 5.4. *There exist positive constants $C_1, C_2, K_1, K_2 > 0$ such that for all $j$, and $t > t_0$,*

$$(5.3) \quad P\{\omega \colon \ell_j(t) = \ell\left(\Delta_{x_j(t), t}(t_0)\right) \leq 1 - K_1 \exp(-C_1(t - t_0))\} \\ \leq K_2 \exp(-C_2(t - t_0)).$$

*Proof.* Fix any $j \in \mathbb{N}$. The position of the $j$-th shock at time $t$ will be denoted by $x(t)$. The following estimates are independent of $j$.

Consider a sequence of times $t_i = t_0 + i$, $i = 0, 1, 2, \ldots$. For each $i$ let $I_i = \mathbb{S}^1 \setminus \Delta_{x(t_i), t_i}(t_0)$ and $z_i$ be the mid-point of $I_i$. Denote by $y_i$ a point on $\mathbb{S}^1$ at time $t_i$ which corresponds to $z_i$ at $t = t_0$; i.e., either $(y_i, t_i)$ is a point of shock which covers $(z_i, t_0)$, or there is a unique one-sided minimizer at $(y_i, t_i)$ which passes through $(z_i, t_0)$. Clearly, $y_i$ is measurable with respect to $\mathcal{F}^{t_i}_{-\infty}$. Denote by $\eta_i$ a random variable which takes the value 1 if $y_i$ is covered by $\Delta_{x(t_{i+1}), t_{i+1}}(t_i)$ and 0 otherwise. Obviously $\eta_i = 1$ if and only if $y_i$ merges with $x(t_i)$ before $t_{i+1}$.

Notice that if $\eta_i = 1$ then the length of the complement to the interval $\Delta_{x(t_{i+1}), t_{i+1}}(t_0)$ is no more than half of the length of the complement to

$\Delta_{x(t_i),t_i}(t_0)$. For any fixed positive integer $K$, with $n = t - t_0$,

(5.4)
$$P\{\omega: \ell(I_n) \geq 2^{-K}\} \leq P\left(\sum_{i=0}^{n-1} \eta_i \leq K\right)$$
$$\leq \sum_{m=0}^{K} \sum_{0 < i_1 < i_2 < \ldots < i_{n-m} \leq n-1} P\{\eta_{i_1} = \eta_{i_2} = \ldots = \eta_{i_{n-m}} = 0\}$$
$$\leq \sum_{m=0}^{K} C_m^n (1 - p_0(1))^{n-m}$$

where we used Lemma 5.2 and the Markov property to conclude that

(5.5) $$P(\eta_{i_1} = \eta_{i_2} = \ldots = \eta_{i_s} = 0) \leq (1 - p_0(1))^s,$$

for all $0 \leq i_1 < i_2 < \ldots < i_s$.

It follows from (5.4) that for $K \leq \frac{n}{2}$

(5.6) $$P\{\omega: \ell(I_n) \geq 2^{-K}\} \leq (K+1) C_n^K (1 - p_0(1))^{n-K}.$$

Let $K = [\alpha n]$ and choose $\alpha$ so small that

(5.7) $$q_1 = \left(\frac{1}{\alpha}\right)^\alpha \left(\frac{1}{1-\alpha}\right)^{1-\alpha} (1 - p_0(1))^{1-\alpha} < 1.$$

Then,

(5.8) $$P\{\omega: \ell(I_n) \geq e^{-[\alpha n]\ln 2}\} \leq M\sqrt{n}\, q_1^n$$

where $M$ is an absolute constant. It follows that

(5.9) $$P\{\omega: \ell\left(\Delta_{x(t),t}(t_0)\right) \leq 1 - e^{-[\alpha n]\ln 2}\} \leq M\sqrt{t - t_0}\, q_1^{t-t_0}.$$

Take any $q$ such that $q_1 < q < 1$ and let $C_1 = \alpha \ln 2$, $K_1 = 4$, $C_2 = -\ln q$. Then (5.3) follows from (5.9) for large enough $K_2$. □

Using the Borel-Cantelli lemma, one gets from Lemma 5.4 the following:

LEMMA 5.5. *For almost all $\omega$*

(5.10) $$\ell\left(\Delta_{x_j(t),t}(t_0)\right) \to 1 \quad as \quad t \to \infty.$$

*Moreover for any shock at $\xi_j(t_0)$ there exists a random constant $T_j(\omega, t_0)$ such that for all $t > T_j(\omega, t_0)$*

(5.11) $$\ell\left(\Delta_{x_j(t),t}(t_0)\right) \geq 1 - 2K_1 \exp(-C_1(t - t_0)).$$

*Remark.* Since the intervals $\Delta_{x_j(t),t}(t_0)$ do not intersect each other, Lemma 5.5 implies that the shocks $\xi_{j_1}(t_0), \xi_{j_2}(t_0)$ merge with each other after time $T = \max(T_{j_1}, T_{j_2})$.

Let us define now an object which will play a very important role in the remaining part of this paper.

*Definition* 5.1. A $C^1$ curve $\gamma\colon (-\infty, +\infty) \to \mathbb{S}^1$ is called a *two-sided minimizer* if for any $C^1$ compact perturbation $\gamma + \delta\gamma\colon (-\infty, +\infty) \to \mathbb{S}^1$ of $\gamma$

$$\mathcal{A}_{-s,s}(\gamma + \delta\gamma) \geq \mathcal{A}_{-s,s}(\gamma)$$

for all sufficiently large $s$.

In other words, a curve $\gamma$ is called a two-sided minimizer if and only if for arbitrary $t_0 \in \mathbb{R}^1$, its restriction on $(-\infty, t_0]$ is a one-sided minimizer.

THEOREM 5.1. *With probability* 1 *there exists a unique two-sided minimizer.*

*Proof.* Existence of the two-sided minimizer follows from a compactness argument. Consider a sequence of curves $\gamma^{(n)}\colon [-n, n] \to \mathbb{S}^1$ which minimize $\mathcal{A}_{-n,n}(\gamma)$ in the class of $C^1$ curves. It follows from Lemma 3.3 that $|\dot{\gamma}^{(n)}(0)| \leq C(\omega, 0)$. Hence the sequence of points $(x_0^{(n)}, v_0^{(n)}) = (\gamma^{(n)}(0), \dot{\gamma}^{(n)}(0))$ belongs to a compact set $\mathbb{S}^1 \times [-C(\omega, 0), C(\omega, 0)]$. Then there exists at least one limiting point $(x_0, v_0)$. A standard argument as in the proof of Theorem 3.1 shows that the solution of the Euler-Lagrange equation (2.3) with initial conditions $x(0) = x_0$, $v(0) = v_0$ defines a two-sided minimizer.

For uniqueness notice that points on a two-sided minimizer $\gamma$ do not belong to the intervals $\Delta_{x_j(t),t}(t_0)$ for any $j$. Since $\ell\left(\Delta_{x_j(t),t}(t_0)\right) \to 1$ as $t \to \infty$, the two-sided minimizer is unique. $\square$

Denoting the two-sided minimizer by $y_\omega$, we now construct another important object, the main shock. For arbitrary $t_0 \in \mathbb{R}^1$ consider a sequence of non-intersecting intervals $\Delta_{x_j,t_0}(t)$, $t \leq t_0$ corresponding to shocks $\xi_j(t_0) = (x_j, t_0)$ at time $t_0$, $x_j \in J(\omega, t_0)$. Notice that here we consider intervals $\Delta_{x_j,t_0}(t)$ for $t \leq t_0$. It turns out that for almost all $\omega$ there exists a unique shock at $(z(t_0), t_0)$ for which $\ell(t) = \ell(\Delta_{z(t_0),t_0}(t)) \to 1$ as $t \to -\infty$.

THEOREM 5.2. *For almost all $\omega$ the following statements hold.*

(a) *For any $t_0 \in \mathbb{R}^1$ there exists a unique shock at $(z(t_0), t_0)$ such that*

(5.12) $$\ell(t) = \ell\left(\Delta_{z(t_0),t_0}(t)\right) \to 1 \quad \text{as} \quad t \to -\infty.$$

*Moreover, for any $\delta > 0$ there exists a random constant $T_{\delta,t_0}(\omega)$ such that for all $t < T_{\delta,t_0}(\omega)$*

(5.13) $$\ell(t) \geq 1 - \exp(-(C_1 - \delta)(t_0 - t)).$$

*The position of the shock $(z(t_0), t_0)$ is measurable with respect to the $\sigma$-algebra $\mathcal{F}_{-\infty}^{t_0}$.*

(b) *For all other shocks $\xi_j(t_0) = (x_j(t_0), t_0)$*

$$\ell_j(t) = \ell\left(\Delta_{x_j(t_0),t_0}(t)\right) \to 0 \quad \text{as} \quad t \to -\infty.$$

(c) $\{(z(t), t), \ t \in \mathbb{R}^1\}$ *is a Lipschitz continuous curve.*

*Proof.* Consider a sequence of times $\bar{t}_i = t_0 - i$. It follows from Lemma 5.4 that with probability greater than $1 - K_2 \exp(-C_2 i)$ there exists a shock at some point $(x_0(i), t_0)$ such that

(5.14) $$\ell\left(\Delta_{x_0(i),t_0}(\bar{t}_i)\right) \geq 1 - K_1 \exp(-C_1 i).$$

By the Borel-Cantelli lemma, there exists $N_1(\omega)$ such that for all $i > N_1(\omega)$, (5.14) holds for some shock at $(x_0(i), t_0)$. We will show that for $i$ large enough $x_0(i)$ does not depend on $i$. Suppose $x_0(i+1) \neq x_0(i)$ for some $i > N_1$. Then,

$$\ell\left(\Delta_{x_0(i+1),t_0}(\bar{t}_i)\right) \leq K_1 \exp(-C_1 i)$$

and

$$\ell\left(\Delta_{x_0(i+1),t_0}(\bar{t}_{i+1})\right) \geq 1 - K_1 \exp(-C_1(i+1)).$$

Denote by $a_i, b_i, a_{i+1}, b_{i+1}$ the end points of $\Delta_{x_0(i+1),t_0}(\bar{t}_i)$ and $\Delta_{x_0(i+1),t_0}(\bar{t}_{i+1})$, respectively, and by $v(a_i), v(b_i), v(a_{i+1}), v(b_{i+1})$ the velocities of the corresponding one-sided minimizers. It follows from Lemma B.8 that

$$|v(a_i) - v(b_i)| \leq L_i |a_i - b_i|,$$

where $L_i = L_0(\theta^{t_0-i}\omega)$. Thus,

(5.15)
$$\mathcal{D}_i = \operatorname{dist}((a_i, v(a_i)), (b_i, v(b_i))) \leq \sqrt{1 + L_i^2}\, K_1 \exp(-C_1 i),$$

$$\mathcal{D}_{i+1} = \operatorname{dist}((a_{i+1}, v(a_{i+1})), (b_{i+1}, v(b_{i+1}))) \geq 1 - K_1 \exp(-C_1(i+1)).$$

On the other hand, we have $\mathcal{D}_{i+1} \leq \exp(d_i)\mathcal{D}_i$, where $d_i = d_1(\theta^{t_0-i}\omega)$ and $d_1(\omega)$ is as defined in Lemma B.5. It follows that

(5.16) $$\exp(d_i) \geq \frac{\mathcal{D}_{i+1}}{\mathcal{D}_i} \geq \frac{1}{2\sqrt{1 + L_i^2}\, K_1} \exp(C_1 i).$$

Since $L_i$ and $d_i$ have finite expectations, it follows that for any $\varepsilon > 0$ there exists $N_\varepsilon(\omega)$ such that

(5.17) $$|d_i| \leq \varepsilon i, \ |L_i| \leq \varepsilon i, \quad \text{for all} \quad i > N_\varepsilon(\omega)$$

(see Lemma 6.2). Take $\varepsilon < C_1$. Then for $i > \bar{N}_\varepsilon(\omega) > N_\varepsilon(\omega)$, (5.16) and (5.17) contradict each other. Hence, for $i > \bar{N}_\varepsilon(\omega)$, $x_0(i+1) = x_0(i)$. Define $z(t_0) = x_0(i)$, $i > \bar{N}_\varepsilon(\omega)$. Obviously, $\ell(\Delta_{z(t_0),t_0}(\bar{t}_i)) \geq 1 - K_1 \exp(-C_1 i)$, $i > \bar{N}_\varepsilon(\omega)$. It follows from the estimates (5.15)–(5.17) that for any $\delta > 0$ there exists a random constant $T_{\delta,t_0}(\omega)$ such that for all $t < T_{\delta,t_0}(\omega)$, (5.13) holds. A shock that satisfies (5.12) is obviously unique. Clearly $z(t_0)$ is measurable with respect to $\mathcal{F}^{t_0}_{-\infty}$. Notice that for any $t > t_0$, the shock which covers $(z(t_0), t_0)$ also satisfies (5.12), (5.13). Hence for almost all $\omega$ such a shock exists for all $t_0 \in \mathbb{R}^1$. (a) is now proven. (b) follows from (a), since

$$0 \leq \sum_{x_j \in J(\omega, t_0)} \ell\left(\Delta_{x_j(t_0), t_0}(t)\right) \leq 1 \quad \text{and} \quad \ell\left(\Delta_{z(t_0), t_0}(t)\right) \to 1 \quad \text{as} \quad t \to -\infty.$$

Since the shock at $(z(\tilde{t}_0), \tilde{t}_0)$ covers $(z(t_0), t_0)$ for all $\tilde{t}_0 \geq t_0$, we get (c). □

*Definition* 5.2. The shock $(z(t), t)$ constructed in Theorem 5.2 is called the *main shock* at time $t$.

As remarked in the introduction, the two-sided minimizer and the main shock play dual roles. The former acts as a repeller, the latter acts as an attractor. Indeed, it follows from Theorem 5.2 that for any two one-sided minimizers $\gamma_1$, $\gamma_2$, $\text{dist}(\gamma_1(t), \gamma_2(t)) \to 0$ as $t \to -\infty$. One can say that all one-sided minimizers approach the two-sided minimizer as $t \to -\infty$.

LEMMA 5.6. (a) *For any two minimizers* $\gamma_1$, $\gamma_2$

(5.18) $$\text{dist}(\gamma_1(t), \gamma_2(t)) \to 0 \quad \text{as} \quad t \to -\infty.$$

*Moreover, for any* $\delta > 0$ *there exists a random constant* $T^{1,2}_\delta(\omega)$ *such that for all* $t < T^{1,2}_\delta(\omega)$

$$\text{dist}(\gamma_1(t), \gamma_2(t)) \leq \exp(-(C_1 - \delta)t).$$

*If* $\gamma_1$ *starts at time* $t_1$, *and* $\gamma_2$ *starts at* $t_2$, *then* $T^{1,2}_\delta(\omega) = T_{\delta,t^{1,2}}(\omega)$, *where constant* $T_{\delta,t}(\omega)$ *is defined as in Theorem 5.2 and* $t^{1,2} = \min(t_1, t_2)$.

*For minimizers starting at the same time convergence in* (5.18) *is uniform.*

(b) *Any shock at a given time* $t$ *will be eventually absorbed by the main shock.*

This is obvious.
Another way to characterize the curve of the main shock is to say that it is the only shock curve defined for all $t \in \mathbb{R}^1$.

LEMMA 5.7. *For almost all* $\omega$ *there exists a unique shock curve*

$$x^\omega : (-\infty, +\infty) \to \mathbb{S}^1$$

*such that* $u^\omega_+(x^\omega(t), t) < u^\omega_-(x^\omega(t), t)$ *and* $x^\omega(t)$ *is measurable with respect to* $\mathcal{F}^t_{-\infty}$. *This curve is the curve of the main shock.*

*Proof.* The existence follows from the existence of the main shock. Suppose now that there exists another measurable shock curve $x^\omega(\cdot)$ defined on $\mathbb{R}^1$. Fix arbitrary $t_0 \in \mathbb{R}^1$. It follows from Lemma 5.2 that with probability 1 the curves $z(t)$ and $x^\omega(t)$ merge before $t = t_0$. Since $t_0$ is arbitrary, $x^\omega(t)$ coincides with $z(t)$ with probability 1. □

*Remark.* Later we will prove a stronger result: the curve of the main shock is the only shock curve which is defined for all sufficiently negative times.

We end this section with some discussion on the assumption (A1). Obviously, a necessary condition for the main result of this section to hold, namely the existence of a unique main shock and two-sided minimizer, is that the minimum period of all the $F_k$'s is equal to 1. However, this condition is not sufficient, as we show now.

THEOREM 5.3. *If $F(x,t) = \cos 2\pi x\, dB(t)$, then with probability 1, there are at least two main shocks, i.e., shock curves defined for all negative times. There are also at least two two-sided minimizers on $\mathbb{S}^1 \times \mathbb{R}^1$.*

We will give an outline of the proof showing that with probability one, $x = 0$ and $x = 1/2$ are points of shock for any time $t$. The main point is:

LEMMA 5.8. *With probability 1, $\xi\colon (-\infty, t] \to \mathbb{R}^1$, $\xi(s) \equiv 0$, is not a one-sided minimizer.*

This follows from the observation that with probability 1, there are large intervals on which $B(t) - B(s) > 0$. On such intervals, one-sided minimizers are close to $x = \frac{1}{2}$. Hence $\xi \equiv 0$ is not minimizing.

As a consequence of symmetry, if $\xi\colon (-\infty, t] \to \mathbb{R}^1$ is a one-sided minimizer such that $\xi(t) = 0$, then $-\xi$ is also a one-sided minimizer. Therefore, with probability 1, $x = 0$ is a point of shock for all $t \in \mathbb{R}^1$. The same argument applies to $x = \frac{1}{2}$. So there are at least two main shocks. The rest of the statement in Theorem 5.3 follows from the same argument as in the proof of Theorems 5.1 and 5.2.

It is easy to check that (A1) holds for
(5.19)
$$F(x,t) = \cos 2\pi(x + x_1)dB_1(t) + \cos 2\pi(x + x_2)dB_2(t) + \cos 2\pi(x + x_3)dB_3(t)$$

where $x_1, x_2, x_3$ are fixed constants such that their differences are not integral multiples of $\frac{1}{2}$. By refining the argument in Appendix D, one can actually show that Lemma 5.2 also holds if there are only two terms in (5.19). On the other hand, without shifting phases, (A1) does not hold if all of the $F_k$'s are of the form $\cos 2\pi kx$. It fails when $x_1^0 = 0$, $x_2^0 = \frac{1}{2}$. However, the following argument shows that Lemma 5.2 still holds if:

The set $\{F_k\}$ contains either

$\{\sin 2\pi x, \cos 2\pi l x, \text{ for some } l \neq 0\}$, or $\{\cos 2\pi x, \sin 2\pi l x, \text{ for some } l \neq 0\}$.

We will illustrate how this claim can be proved when $\{F_k\}$ contains $\{\cos 2\pi x, \sin 4\pi x\}$. The only situation we have to reconsider is when $x_1^0$ is close to a critical point of $F_1(x) = \cos 2\pi x$, and $x_2^0$ is close to a critical point of $F_2(x) = \sin 4\pi x$. Without loss of generality, let us assume that $x_1^0$ is close to 0, and $x_2^0$ is close to $\frac{1}{8}$. Heuristically we can first use $F_1$ to move $x_2^0$ to a small neighborhood of $\frac{1}{2}$. If in this process $x_1^0$ has moved out of the neighborhood of 0, then we can use Lemma 5.2 with $F_1$. If not, we use $F_2$ to move $x_1^0$ to a small neighborhood of $-\frac{1}{8}$. The forces $dB_2$ can be chosen such that $x_2^0$ will stay inside $\left(\frac{1}{8}, \frac{5}{8}\right)$. Now both $x_1^0$ and $x_2^0$ are inside the region where $F_1'$ is bounded away from 0, so we can apply the proof of Lemma 5.2 to $x_1^0$ and $x_2^0$ with the potential $F_1$. We will omit the detailed proof of these statements since they follow closely the proof of Lemma 5.2.

## 6. Hyperbolicity and unstable manifolds of the two-sided minimizer

In this section we prove that the two-sided minimizer $y_\omega(\cdot)$ constructed in Section 5 is a hyperbolic trajectory of the dynamical system (2.3), and we establish the existence of its stable and unstable manifolds. The main technical difficulty is associated with the fact that the hyperbolicity is nonuniform, as in many other dynamical systems with noises. This is overcome by using Pesin's theory (see [Pes], [En]). Note in passing that

$$y_\omega(t+s) = y_{\theta^s \omega}(t).$$

Denote by $G_t^\omega$ the stochastic flow generated by the solutions of (2.3). Let $J_s^t(\omega)$ be the Jacobi map, i.e., the tangent map that maps the tangent plane $T(y_\omega(s), u^\omega(y_\omega(s), s))$ onto the tangent plane $T(y_\omega(t), u^\omega(y_\omega(t), t))$. This is well-defined since $y_\omega(t)$ is a point of continuity of $u^\omega(\,\cdot\,, t)$ for all $t$. Moreover the Jacobi map has determinant 1 since the dynamical system (2.3) preserves the Lebesgue measure. Obviously, we have

$$J_0^{t_2}(\omega) = J_0^{t_2-t_1}(\theta^{t_1}\omega) J_0^{t_1}(\omega)$$

for all $t_1, t_2$. In the terminology of ergodic theory, $\{J_0^t(\omega)\}$ is a cocycle (see [O]).

LEMMA 6.1. *Define* $\log^+ x = \max(\log x, 0)$, *for* $x > 0$. *Then*

$$\sup_{-1 \leq t \leq 1} \log^+ \|J_0^t(\omega)\| \in L^1(dP).$$

This result, together with some other technical estimates, is proved in Appendix B (Lemma B.5).

As a consequence of the multiplicative ergodic theorem [O], [En], we conclude that with probability 1:

(A) either
$$\lim_{t \to \pm\infty} \frac{1}{t} \ln \|J_{t_1}^{t+t_1}(\omega)e\| = 0$$
for all $e \in T_{t_1} = T(y_\omega(t_1), u^\omega(y_\omega(t_1), t_1))$;

(B) or there exist a constant $\lambda > 0$ and a measurable normalized basis $\{e_t^u(\omega), e_t^s(\omega)\}$ of $T_t = T(y_\omega(t), u^\omega(y_\omega(t), t))$, such that
$$J_{t_1}^t(\omega) e_{t_1}^u(\omega) = a^u(t, t_1; \omega) e_t^u(\omega), \quad J_{t_1}^t(\omega) e_{t_1}^s(\omega) = a^s(t, t_1; \omega) e_t^s(\omega),$$
where the functions $a^u(t, t_1; \omega)$ and $a^s(t, t_1; \omega)$ are also cocycles satisfying
$$a^u(t+s, 0; \omega) = a^u(s, 0; \theta^t \omega) a^u(t, 0; \omega),$$
$$a^s(t+s, 0; \omega) = a^s(s, 0; \theta^t \omega) a^s(t, 0; \omega).$$
Furthermore,
$$\lim_{t \to \infty} \frac{\ln a^u(t, t_1; \omega)}{t - t_1} = \lambda, \quad \lim_{t \to \infty} \frac{\ln a^s(t, t_1; \omega)}{t - t_1} = -\lambda.$$

If (B) holds, the cocycle $\{J_s^t(\omega)\}$ is said to be hyperbolic and the basis $\{e_t^u(\omega), e_t^s(\omega)\}$ is called the Oseledetz basis.

THEOREM 6.1. *With probability 1, the cocycle $\{J_s^t(\omega)\}$ is hyperbolic.*

We will prove Theorem 6.1 later. It is useful to recall the following simple result:

LEMMA 6.2. *Let $\{\eta_i\}$ be a sequence of identically distributed random variables such that $E|\eta_i| < +\infty$. Then for any $\varepsilon > 0$, there exists a random variable $N_\varepsilon > 0$, such that for all $i$, $|i| \geq N_\varepsilon$,*
$$|\eta_i| \leq \varepsilon |i|.$$

This is a simple consequence of the Chebyshev inequality and the Borel-Cantelli lemma. Lemma 6.2 is equivalent to the statement that
$$\lim_{i \to \infty} \frac{\eta_i}{i} = 0$$
with probability 1. However, we will use it in the form of Lemma 6.2.

Let $x(\cdot)$ be an arbitrary one-sided minimizer defined on $(-\infty, 0]$. Fix a positive integer $k$ and consider a sequence of times $t_i = -ki$.

Denote $(y_i, u_i) = (y_\omega(t_i), u^\omega(y_\omega(t_i), t_i))$, $(x_i, v_i) = (x(t_i), u^\omega(x(t_i), t_i))$, $J_i = J_{t_i}^{t_{i-1}}(\omega)$, $\ell_i = \text{dist}(x_i, y_i)$, $\rho_i = \text{dist}((x_i, v_i), (y_i, u_i))$.

LEMMA 6.3. *For any $\varepsilon > 0$ there exists a random constant $\ell_{\varepsilon,k}(\omega)$ such that with probability 1*

$$\rho_i \leq (1+\varepsilon)\|J_{i+1}\|\rho_{i+1}, \quad i \geq 0, \tag{6.1}$$

*provided that $\ell_0 \leq \ell_{\varepsilon,k}(\omega)$.*

*Proof.* Let $L_0(i) = L_0(\theta^{t_i}\omega)$, $d(i) = d_k(\theta^{t_i}\omega)$, $\bar{d}(i) = \bar{d}_k(\theta^{t_i}\omega)$, where $L_0$ is as defined in Lemma B.8, and $d_k$, $\bar{d}_k$ are as defined in (B.29–30) and Lemma B.5. It follows from Theorem 5.2 that for any $\delta > 0$ there exists a random constant $N_\delta(\omega)$ such that for all $i > N_\delta(\omega)$

$$\ell_i \leq \exp((C_1 - \delta)t_i). \tag{6.2}$$

Since $|v_i - u_i| \leq L_0(i)\ell_i$, we have

$$\rho_i \leq \sqrt{L_0^2(i) + 1}\,\ell_i \leq \sqrt{L_0^2(i) + 1}\,\exp((C_1 - \delta)t_i). \tag{6.3}$$

Let $\Delta_i = \{(x, v) = \alpha(x_i, v_i) + (1 - \alpha)(y_i, u_i), 0 \leq \alpha \leq 1\}$ be the interval connecting $(x_i, v_i)$ and $(y_i, u_i)$. Clearly $\Delta_i \in B_k(\theta^{t_i}\omega)$. It follows from the definition of $d(i)$ that for any $(x, v) \in \Delta_i$

$$\|G_t^{\theta^{t_i}\omega}(x, v) - G_t^{\theta^{t_i}\omega}(y_i, u_i)\| \leq \exp(d(i))\rho_i, \quad -k \leq t \leq 0. \tag{6.4}$$

Since $d(i)$, $\bar{d}(i)$, $L_0(i)$ have finite expectations, for any $\nu > 0$, there exists $N_\nu(\omega)$ such that

$$|d(i)|, \ |\bar{d}(i)|, \ |L_0(i)| \leq \nu i \quad \text{for} \quad i > N_\nu(\omega). \tag{6.5}$$

Hence, for $i > \max(N_\delta(\omega), N_\nu(\omega))$
$$\|G_t^{\theta^{t_i}\omega}(x, v) - G_t^{\theta^{t_i}\omega}(y_i, u_i)\| \leq \sqrt{\nu^2 i^2 + 1}\,\exp(\nu i + (C_1 - \delta)t_i), \quad -k \leq t \leq 0. \tag{6.6}$$

Take $\nu < (C_1 - \delta)$. Then (6.6) implies that there exists

$$N_{\delta,\nu}(\omega) > \max(N_\delta(\omega), N_\nu(\omega))$$

such that $\Delta_i \subset O_k(\theta^{t_i}\omega)$ for all $i > N_{\delta,\nu}(\omega)$. Clearly, if $\rho_0$ is small enough, then $\Delta_i \subset O_k(\theta^{t_i}\omega)$ for $i \leq N_{\delta,\nu}(\omega)$. Since the two-sided minimizer corresponds to a point of continuity of $u_+^\omega$, we have $\rho_0 \to 0$ as $\ell_0 \to 0$. Thus, there exists $\overline{\ell_0}(\omega) > 0$ such that $\Delta_i \in O_k(\theta^{t_i}\omega)$ for all $i$, provided that $\ell_0 \leq \overline{\ell_0}(\omega)$. Denote now $D(i) = D_{k,2}(\omega)$, $\bar{D}(i) = \bar{D}_{k,2}(\omega)$, where $D_{T,r}(\omega)$, $\bar{D}_{T,r}(\omega)$ are defined as in (B.22). Since $\Delta_i \in O_k(\theta^{t_i}\omega)$, we have for all $i \geq 0$

$$\rho_i \leq \|J_{i+1}\|\rho_{i+1} + \frac{1}{2}\exp(\bar{D}(i))\rho_{i+1}^2 = \|J_{i+1}\|\rho_{i+1}\left(1 + \frac{\exp(\bar{D}(i))}{2\|J_{i+1}\|}\rho_{i+1}\right).$$

Since $\frac{1}{\|J_{i+1}\|} \leq \|J_{i+1}^{-1}\| \leq \exp(D(i))$,

$$\text{(6.7)} \qquad \rho_i \leq \|J_{i+1}\| \rho_{i+1} \left(1 + \frac{1}{2}\exp(D(i) + \bar{D}(i))\rho_{i+1}\right).$$

Again, since $D(i)$ and $\bar{D}(i)$ have finite expectations, for any $\nu > 0$, there exists $\overline{N_\nu}(\omega)$ such that

$$\text{(6.8)} \qquad |D(i)|, \ |\bar{D}(i)| \leq \nu i \quad \text{for} \quad i > \overline{N_\nu}(\omega).$$

Thus, by (6.3), (6.5), (6.8), for $i > \max(N_\delta(\omega), N_\nu(\omega), \overline{N_\nu}(\omega))$:

$$\text{(6.9)} \qquad \frac{1}{2}\exp(D(i) + \bar{D}(i))\rho_{i+1} \leq \frac{1}{2}\sqrt{\nu^2 i^2 + 1} \ \exp(2\nu i + (C_1 - \delta)t_i).$$

Take $\nu < \frac{C_1 - \delta}{2}$. It follows from (6.9) that there exists $N(\omega)$ such that for $i > N(\omega)$

$$\text{(6.10)} \qquad \frac{1}{2}\exp(D(i) + \bar{D}(i))\rho_{i+1} \leq \varepsilon.$$

This implies (6.1) for $i > N(\omega)$. Now, in order to get (6.1) for all $i$, take $\rho_0(\omega)$ so small that for $\rho_0 \leq \rho_0(\omega)$

$$\text{(6.11)} \qquad \frac{1}{2}\exp(D(i) + \bar{D}(i))\rho_{i+1} \leq \varepsilon, \quad 0 \leq i \leq N(\omega).$$

As above, we can choose $\ell_{\varepsilon,k}(\omega) < \overline{\ell_0}(\omega)$ so small that $\ell_0 < \ell_{\varepsilon,k}(\omega)$ implies $\rho_0 \leq \rho_0(\omega)$. (6.1) obviously follows from (6.10), (6.11). □

*Proof of Theorem* 6.1. Assume that (A) holds. It follows from the subadditive ergodic theorem that

$$\text{(6.12)} \qquad \lim_{n \to \infty} \frac{\int \ln \|J_0^n(\omega)\| P(d\omega)}{n} = 0.$$

Then, for any $\varepsilon > 0$ there exists $k \in \mathbb{N}$ such that $A_k = \frac{1}{k}\int \ln \|J_0^k(\omega)\| P(d\omega) < \varepsilon$. By the ergodic theorem, with probability 1,

$$\text{(6.13)} \qquad \frac{1}{kn}\sum_{i=1}^n \ln \|J_0^k(\theta^{t_i}\omega)\| \xrightarrow[n \to \infty]{} \frac{1}{k}\int \ln \|J_0^k(\omega)\| P(d\omega) = A_k.$$

Hence there exists a random constant $n_\varepsilon(\omega)$ such that, with probability 1,

$$\text{(6.14)} \qquad \frac{1}{kn}\sum_{i=1}^n \ln \|J_0^k(\theta^{t_i}\omega)\| \leq A_k + \varepsilon \leq 2\varepsilon$$

for all $n > n_\varepsilon(\omega)$.

Consider now a one-sided minimizer $x(\cdot)$ at time $t_0 = 0$ such that $\ell_0 = |x(0) - y_\omega(0)| \leq \ell_{\varepsilon,k}(\omega)$, where $\ell_{\varepsilon,k}(\omega)$ is defined as in Lemma 6.3. Then, by Lemma 6.3 for all $n > 0$:

$$\ell_0 \leq \rho_0 \leq (1+\varepsilon)^n \prod_{i=1}^n \|J_{t_i}^{t_{i-1}}(\omega)\| \rho_n. \tag{6.15}$$

Thus for $n > n_\varepsilon(\omega)$

$$\rho_n \geq \frac{\rho_0}{(1+\varepsilon)^n \prod_{i=1}^n \|J_{t_i}^{t_{i-1}}(\omega)\|} = \frac{\rho_0}{(1+\varepsilon)^n \exp\left(\sum_{i=1}^n \ln \|J_0^k(\theta^{t_i}\omega)\|\right)} \tag{6.16}$$
$$\geq \rho_0 e^{-\varepsilon n} \exp(-2\varepsilon k n).$$

On the other hand, it follows from Theorem 5.2 that for large enough $n$

$$\rho_n \leq \sqrt{L_0^2(n) + 1} \, \ell_n \leq \sqrt{L_0^2(n) + 1} \, \exp(-(C_1 - \delta)kn) \tag{6.17}$$
$$\leq \sqrt{\nu^2 n^2 + 1} \, \exp(-(C_1 - \delta)kn).$$

Here, as in the proof of Lemma 6.3, we used again that $|L_0(n)| \leq \nu n$ for $n$ large enough. Take $\varepsilon$ so small that $3\varepsilon < C_1 - \delta$. Then, (6.16) and (6.17) are contradictory to each other. □

*Remark.* It follows from the proof of Theorem 6.1 that $\lambda \geq C_1 - \delta$. Since $\delta$ is arbitrarily small, $\lambda \geq C_1$.

Next we construct stable and unstable manifolds of the two-sided minimizer. We will denote by $\Gamma_\omega$ the trajectory in the phase space of the two-sided minimizer $\Gamma_\omega = \{(y_\omega(t), u^\omega(y_\omega(t), t)), t \in \mathbb{R}^1\}$, and let $(x(t; x_0, u_0), u(t; x_0, u_0))$ be the solution of the SDE (2.3) with initial data $x(0) = x_0$, $u(0) = u_0$. We will concentrate on $t = 0$ but the same holds for any other $t$.

*Definition* 6.1. A *local stable manifold* of $\Gamma_\omega$ at $t = 0$ is the set

$$W_{\delta,\varepsilon}^s = \{(x_0, u_0), \ \mathrm{dist}((x(t; x_0, u_0), u(t; x_0, u_0)),$$
$$(y_\omega(t), u^\omega(y_\omega(t), t))) \leq \delta e^{-(\lambda-\varepsilon)t}\}$$

for some $\varepsilon > 0$, $\delta > 0$ and all $t > 0$. A *local unstable manifold* of $\Gamma_\omega$ at $t = 0$ is the set

$$W_{\delta,\varepsilon}^u = \{(x_0, u_0), \ \mathrm{dist}((x(t; x_0, u_0), u(t; x_0, u_0)),$$
$$(y_\omega(t), u^\omega(y_\omega(t), t))) \leq \delta e^{-(\lambda-\varepsilon)|t|}\}$$

for some $\varepsilon > 0$, $\delta > 0$ and all $t < 0$.

Pesin [Pes] gave general conditions under which such local stable and unstable manifolds exist for smooth maps of compact manifolds. It is easy to check that his results can be extended directly to the current situation of stochastic flows. Below we will formulate Pesin's theorem and later verify that its conditions are satisfied for our problem.

Denote by $S_i$ the Poincaré map at $t = i$ associated with the SDE (2.3). In other words, $S_i$ maps $(x_i, u_i)$ at $t = i$ to the solution of (2.3), $(x_{i+1}, u_{i+1})$ at $t = i + 1$. Similarly we denote by $S_i^n, S_i^{-n}$ the maps that map the solution of (2.3) at $t = i$ to the solution at $t = i + n$, $t = i - n$ respectively.

Define $\lambda_i^u, \lambda_i^s$ by the relations

$$J_i^{i+1}(\omega)e_i^u = e^{\lambda_i^u} e_{i+1}^u, \quad J_i^{i+1}(\omega)e_i^s = e^{-\lambda_i^s} e_{i+1}^s$$

where $\{e_i^s, e_i^u\}$ constitutes the Oseledetz basis.

PESIN'S THEOREM. *Assume that there exist constants $\lambda, \mu > 0$, and $\varepsilon_0 \in (0, 1)$, and for $\varepsilon \in (0, \varepsilon_0)$, one can find a positive random variable $C(\varepsilon, \omega)$ such that for $i \in \mathbb{Z}$*

(I) $$\|DS_i^n e_i^s\| \leq C(\varepsilon, \omega)e^{-(\lambda-\varepsilon)n} e^{\varepsilon|i|}$$
$$\|DS_i^{-n} e_i^u\| \leq C(\varepsilon, \omega)e^{-(\mu-\varepsilon)n} e^{\varepsilon|i|},$$

(II) $$|\sin \langle e_i^s, e_i^u \rangle| \geq \frac{1}{C(\varepsilon, \omega)} e^{-\varepsilon|i|}.$$

(III) *Let $r_i = \frac{1}{C(\varepsilon,\omega)} e^{-\varepsilon|i|}$, and*

$$\mathcal{B}_i(\omega) = \{(x, u), \|(x, u) - (y_\omega(i), u^\omega(y_\omega(i), i))\| \leq r_i\}.$$

*Then for some $r \geq 2$,*

$$\sup_{(x,u) \in \mathcal{B}_i} \max_{1 \leq j \leq r} \left( \|D^j S_i(x, u)\|, \|D^j S_i^{-1}(x, u)\| \right) \leq C(\varepsilon, \omega) e^{\varepsilon|i|}.$$

*Under these assumptions, there exist positive $\varepsilon_1(\lambda, \mu, \varepsilon_0)$ and $\delta(\varepsilon)$, defined for $0 < \varepsilon < \varepsilon_1$, and $C^{r-1}$ curves $W_{\delta,\varepsilon}^s$, $W_{\delta,\varepsilon}^u$ in the phase space of the dynamical system (2.3), such that*

(i) *$W_{\delta,\varepsilon}^s$ and $W_{\delta,\varepsilon}^u$ are respectively the stable and unstable manifolds of $\Gamma_\omega$ at $t = 0$. Moreover, they are $C^{r-1}$ graphs on the interval $[-\delta_1(\varepsilon), \delta_1(\varepsilon)]$ for some $\delta_1(\varepsilon) > 0$.*

(ii) *$W_{\delta,\varepsilon}^s \cap W_{\delta,\varepsilon}^u = (y_\omega(0), u^\omega(y_\omega(0), 0))$.*

(iii) *The tangent vectors to $W_{\delta,\varepsilon}^s$ and $W_{\delta,\varepsilon}^u$ at $(y_\omega(0), u^\omega(y_\omega(0), 0))$ are respectively $e_0^s$ and $e_0^u$.*

(iv) *If* $(x, u) \in \mathcal{B}_0$, *and* $n \geq 0$,

$$\mathrm{dist}\left(S_0^{-n}(x,u), (y_\omega(-n), u^\omega(y_\omega(-n), -n))\right) \leq \bar{\delta} e^{-\chi n}$$

*for some constants* $\chi > 0$ *and* $\bar{\delta} > 0$, *then* $(x, u) \in W^u_{\delta, \varepsilon}$.

Our task is reduced to checking the assumptions (I), (II) and (III) in Pesin's Theorem.

To begin with, let us observe that (II) follows from the next argument (see [R]). Since

$$\lim_{t \to \infty} \frac{\ln a^s(0, t)}{t} = \lambda, \quad \lim_{t \to \infty} \frac{\ln a^u(0, t)}{t} = -\lambda$$

and the area of the parallelogram generated by $e_t^s$ and $e_t^u$ is independent of $t$, we have

$$\lim_{t \to \infty} \frac{\ln |\sin \langle e_t^s, e_t^u \rangle|}{t} = 0.$$

To see that (III) holds, define

$$d_i(\omega) = \sup_{-1 \leq t \leq 0} \sup_{(x,u) \in \mathcal{B}_0(\theta^i \omega)} \|DG_t^{\theta^i \omega}(x, u)\|,$$

$$\bar{d}_i(\omega) = \sup_{0 \leq t \leq 1} \sup_{(x,u) \in \mathcal{B}_0(\theta^i \omega)} \|DG_t^{\theta^i \omega}(x, u)\|,$$

where $G_t^\omega$ is the stochastic flow defined earlier.

LEMMA 6.4. *For any $\varepsilon > 0$, there exist random constants $C_1(\varepsilon, \omega)$, $C_2(\varepsilon, \omega)$ such that, with probability 1,*

$$d_i(\omega) \leq C_1(\varepsilon, \omega) e^{\varepsilon |i|}, \quad \bar{d}_i(\omega) \leq C_2(\varepsilon, \omega) e^{\varepsilon |i|}$$

*for $i \in \mathbb{Z}$.*

*Proof.* Assume $C(\varepsilon, \omega) > 1$. Then it follows from Lemma B.5 that

$$\int \log^+ d_0(\omega) dP < +\infty, \quad \int \log^+ \bar{d}_0(\omega) dP < +\infty.$$

Now Lemma 6.4 follows directly from Lemma 6.2. □

Let $C(\varepsilon, \omega) > \max(C_1(\varepsilon, \omega), C_2(\varepsilon, \omega))$, and

$$d_{i,r}(\omega) = \sup_{(x,v) \in \mathcal{B}_i} \max_{1 \leq j \leq r} \|D^j S_i(x, v)\|,$$

$$\bar{d}_{i,r}(\omega) = \sup_{(x,v) \in \mathcal{B}_i} \max_{1 \leq j \leq r} \|D^j (S_{i-1})^{-1}(x, v)\|.$$

Statement (III) follows from the next result.

LEMMA 6.5. *For any $\varepsilon > 0$, there exist random constants $C_3(\varepsilon, \omega)$ and $C_4(\varepsilon, \omega)$ such that*

$$d_{i,r}(\omega) \leq C_3(\varepsilon, \omega) e^{\varepsilon|i|}, \quad \bar{d}_{i,r}(\omega) \leq C_4(\varepsilon, \omega) e^{\varepsilon|i|}$$

*for $i \in \mathbb{Z}$.*

*Proof.* Let $(x_0, v_0) \in \mathcal{B}_i(\omega)$. Consider the solution of (2.3), $(x(t), v(t))$, such that $x(i) = x_0$, $v(i) = v_0$, for $t \in [i, i+1]$,

$$|v(t) - u^\omega(y_\omega(i+t), i+t)| \leq d_i(\omega) \|(x_0, v_0) - (y_\omega(i), u^\omega(y_\omega(i), i))\| \leq 1$$

from Lemma 6.4. Therefore $\mathcal{B}_i(\omega) \subset O_1(\theta^i \omega)$ where $O_1(\omega)$ is defined as in Appendix B. Lemma 6.5 now follows directly from Lemma B.4 and Lemma 6.2. The second estimate can be proved in the same way. □

Finally, we prove statement (I).

LEMMA 6.6. *For any $\varepsilon > 0$, one can find random constants $C_5(\varepsilon, \omega)$ and $C_6(\varepsilon, \omega)$ such that*

$$\|DS_i^n e_i^s\| \leq C_5(\varepsilon, \omega) e^{-(\lambda-\varepsilon)n} e^{\varepsilon|i|}, \qquad n \geq 1,$$
$$\|DS_i^{-n} e_i^u\| \leq C_6(\varepsilon, \omega) e^{-(\lambda-\varepsilon)|n|} e^{\varepsilon|i|}, \qquad n \leq -1.$$

*Proof.* We will prove the first statement. The second one can be proved in the same way.

From the ergodic theorem,

$$\lim_{n \to +\infty} \frac{1}{n} \sum_{j=0}^{n-1} \lambda_j^s = -\lambda < 0.$$

Thus for any $\varepsilon \in (0, \lambda)$, there exists a constant $C_7(\varepsilon, \omega) \geq 0$ such that

$$e^{\sum_{j=0}^{n-1} \lambda_j^s} \leq C_7(\varepsilon, \omega) e^{-n(\lambda-\varepsilon)}.$$

For any $\delta \in (0,1)$, define $K(\delta)$ by:

$$K(\delta) = \inf\{K : \; P(C(\omega) \leq K) \geq \delta\}.$$

Denote

$$\delta_1 = P(C(\omega) \leq K(\delta)) \geq \delta, \; m_1(m, \omega) = \max\{i : \; 1 \leq i \leq m, \; C(\theta^i \omega) \leq K(\delta)\}.$$

Notice that $m_1(m, \omega)$ is defined for large enough $m$. By the ergodic theorem,

$$\lim_{m \to \infty} \frac{\#\{i : \; 1 \leq i \leq m, \; C(\theta^i \omega) \leq K(\delta)\}}{m} = \delta_1$$

where # denotes cardinality. Thus for any $\delta_2 \in (0, \delta_1)$, there exists a random constant $M(\delta_2, \omega)$ such that for all $m > M(\delta_2, \omega)$,

$$\frac{\#\{i: \ 1 \leq i \leq m, \ C(\theta^i \omega) \leq K(\delta)\}}{m} \geq \delta_1 - \delta_2\,.$$

Hence

$$m_1 \geq (\delta_1 - \delta_2)m\,.$$

Consequently for $m > M(\delta_2, \omega)$,

$$\|DS_m^n e_m^s\| = e^{\sum_{j=m}^{m+n-1} \lambda_j^s} = e^{\sum_{j=m_1}^{m+n-1} \lambda_j^s} e^{-\sum_{j=m_1}^{m-1} \lambda_j^s}$$

$$\leq C(\theta^{m_1}\omega) e^{-(m+n-m_1)(\lambda-\varepsilon)} e^{-\sum_{j=m_1}^{m-1} \lambda_j^s}$$

$$\leq K(\delta) e^{-n(\lambda-\varepsilon)} e^{-(m-m_1)(\lambda-\varepsilon)} e^{\sum_{j=m_1}^{m-1}(-\lambda_j^s)}\,.$$

We also have, with $\delta_3 = \delta_1 - \delta_2$

$$\sum_{j=m_1}^{m-1}(-\lambda_j^s) \leq \max_{\delta_3 m \leq k \leq m-1} \sum_{k}^{m-1}(-\lambda_j^s)$$

$$\leq \max_{\delta_3 m \leq k \leq m-1} \sum_{j=k+1}^{m} \log^+ \|D(S_j^{-1})(y_\omega(j), u^\omega(y_\omega(j), j))\|$$

$$\leq \max_{\delta_3 m \leq k \leq m-1} \sum_{j=k+1}^{m} d_1(\theta^j \omega)$$

where $d_1(\omega)$ is defined as in (B.29). Using Lemma B.5 and standard probabilistic estimates, one can show that for appropriate $\delta_3 < 1$ there exists a constant $M_1(\omega)$ such that for all $m > M_1(\omega)$,

$$\sum_{j=m_1}^{m-1}(-\lambda_j^s) \leq \max_{\delta_3 m \leq k \leq m-1} \sum_{j=k+1}^{m} d_1(\theta^j \omega) \leq \varepsilon m\,.$$

In fact, it is enough to have $\delta_3$ so close to 1 that $(1-\delta_3)\mathbb{E}d_1(\omega) < \frac{\varepsilon}{2}$, where $\mathbb{E}d_1(\omega) = \int d_1(\omega) P(d\omega)$. Hence we choose $\delta = 1 - \frac{\varepsilon}{8\mathbb{E}d_1(\omega)}, \delta_2 = \frac{\varepsilon}{8\mathbb{E}d_1(\omega)}$. Then for $m > \max(M_1(\omega), M(\delta_2, \omega))$,

$$\|DS_m^n e_m^s\| \leq K(\delta) e^{-n(\lambda-\varepsilon)} e^{\varepsilon m}\,. \qquad \square$$

This completes the verification of the assumptions in Pesin's theorem, and establishes the existence of local stable and unstable manifolds $W_{\delta,\varepsilon}^s$, $W_{\delta,\varepsilon}^u$.

One can also define, in a standard way, global stable and unstable manifolds $W^s$, $W^u$:

$$W^s = \bigcup_{i=1}^{\infty} S_i^{-i} W_{\delta,\varepsilon}^s(y_{\theta^i\omega}(0), u^\omega(y_{\theta^i\omega}(0), 0)),$$

$$W^u = \bigcup_{i=1}^{\infty} S_{-i}^i W_{\delta,\varepsilon}^u(y_{\theta^{-i}\omega}(0), u^\omega(y_{\theta^{-i}\omega}(0), 0)).$$

Obviously, $W^s$, $W^u$ are also $C^{r-1}$ curves which coincide with $W_{\delta,\varepsilon}^s$, $W_{\delta,\varepsilon}^u$ in some neighborhood of $(y_\omega(0), u^\omega(y_\omega(0), 0))$.

The following theorem is a consequence of the properties of unstable manifolds.

THEOREM 6.2. a. *The graph $\{(x, u_+^\omega(x, 0)), x \in \mathbb{S}^1\}$ is a subset of the global unstable manifold $W^u$.*

b. *There exists a (random) neighborhood of $(y_\omega(0), u^\omega(y_\omega(0), 0))$ such that $W_{\delta,\varepsilon}^u$ consists of one-sided minimizers in this neighborhood, i.e., the solutions of (2.3) with initial data on $W_{\delta,\varepsilon}^u$ in this neighborhood give rise to one-sided minimizers.*

*Proof.* a. As was shown earlier, any one-sided minimizer $(x, u)$ converges exponentially fast to $(y_\omega(t), u^\omega(y_\omega(t), t))$ as $t \to -\infty$. It follows that $S_0^{-i}(x, u) \in \mathcal{B}_0(\theta^{-i}\omega)$ for some $i > 0$. Hence $S_0^{-i}(x, u)$ lies on the local unstable manifold $W_{\delta,\varepsilon}^u(y_{\theta^{-i}\omega}(0), u^\omega(y_{\theta^{-i}\omega}(0), 0))$, as a consequence of Pesin's theorem (iv), and $(x, u)$ lies on the global unstable manifold.

b. The local unstable manifold $W_{\delta,\varepsilon}^u$ is a $C^{r-1}$ curve with the tangent vector $e_0^u$ at $(y_\omega(0), u^\omega(y_\omega(0), 0))$. Let $M_\alpha = \{(x, u) : x \in (y_\omega(0) - \alpha, y_\omega(0) + \alpha), (x, u)$ corresponds to a one-sided miminizer$\}$. Since $y_\omega(0)$ is a point of (Lipschitz) continuity of $u_+^\omega(x)$ it follows that there exists $\alpha_0(\omega)$ such that for all $\alpha < \alpha_0(\omega)$, $M_\alpha \subset W_{\delta,\varepsilon}^u$. Hence $e_0^u$ is not a vertical vector; i.e., $e_0 \neq (0, 1)$. Therefore there exists a neighborhood $O_\omega$ of $(y_\omega(0), u^\omega(y_\omega(0), 0))$ such that in this neighborhood $W_{\delta,\varepsilon}^u$ is a graph of a $C^{r-1}$ function; i.e.,

$$W_{\delta,\varepsilon}^u \cap O_\omega = \{(x, u) : x \in (y_\omega(0) - \alpha_1(\omega), y_\omega(0) + \alpha_2(\omega)), \ u = \bar{u}(x)\},$$

where $\alpha_1(\omega)$, $\alpha_2(\omega) > 0$, and $\bar{u}(x)$ is a $C^{r-1}$ function. Now choose $\alpha$ so small that $M_\alpha \subset W_{\delta,\varepsilon}^u \cap O_\omega$. It follows that for $x \in (y_\omega(0) - \alpha, y_\omega(0) + \alpha)$, $\bar{u}(x) = u_+^\omega(x, 0)$, which proves b. □

COROLLARY 6.3. *There exists $\alpha(\omega) > 0$ such that there are no shocks inside the interval $(y_\omega(0) - \alpha, y_\omega(0) + \alpha)$. Moreover,*

$$u_+^\omega \in C^{r-1}(y_\omega(0) - \alpha, y_\omega(0) + \alpha).$$

## 7. Regularity of solutions

In this section, we give a complete description of the solution $u_+^\omega$ in terms of the unstable manifold $W^u$, and prove that the number of shocks is finite for almost all $\omega$. We also prove a stronger version of Lemma 5.8, namely, that all shocks except the main shock have finite prehistory. We will start with the latter statement.

THEOREM 7.1. *Fix arbitrary $t_0 \in \mathbb{R}^1$. With probability 1 there exists a random constant $T_{t_0}(\omega)$ such that all shocks at time $t_0$, except the main shock, are generated after the time $t_0 - T_{t_0}(\omega)$. In other words, all shocks at time $t_0 - T_{t_0(\omega)}$ merge with the main shock before $t = t_0$.*

*Proof.* Fix any $t \in \mathbb{R}^1$. It follows from Theorem 6.2 that there exists $\varepsilon_1(t, \omega) > 0$ such that the velocities of all one-sided minimizers in $(y_\omega(t) - \varepsilon_1, y_\omega(t) + \varepsilon_1)$ lie on a $C^{r-1}$ curve. Hence there are no shocks in the $\varepsilon_1$-neighborhood of $y_\omega(t)$. Notice that the random constant $\varepsilon_1(t, \omega)$ has stationary distribution. We can choose an $\alpha > 0$ so small that $P(\omega: \varepsilon_1(t, \omega) > \alpha) > 0$. Then there exists an infinite sequence $t_i \to -\infty$, such that $\varepsilon_1(t_i, \omega) > \alpha, i \in \mathbb{N}$. Since minimizers at $t = t_0$ converge uniformly to the two-sided minimizer as $t \to -\infty$, there exists $I_{t_0}(\omega)$ such that for $i \geq I_{t_0}(\omega)$, all minimizers starting at $t = t_0$ pass through the $\varepsilon_1$-neighborhood of $y_\omega(t_i)$. Now let $T_{t_0}(\omega) = t_0 - t_{I_{t_0}(\omega)}$. We conclude that the complement (on $\mathbb{S}^1$) of the $\varepsilon_1$-neighborhood of $y_\omega(t_{I_{t_0}(\omega)})$ will merge into the main shock before time $t_0$. Since the $\varepsilon_1$-neighborhood of $y_\omega(t_{I_{t_0}(\omega)})$ contains no shocks, this completes the proof of Theorem 7.1. □

Let $s$ be the signed arc-length parameter for the unstable manifold $W^u$ of the two-sided minimizer at $t = 0$:

$$(7.1) \qquad W^u = \{(x(s), u(s)),\ x(s) \in \mathbb{S}^1,\ u(s) \in \mathbb{R}^1\},$$

with $s = 0$ at $(y_\omega(0), u^\omega(y_\omega(0), 0))$. From the proof of Theorem 6.2, $\frac{dx}{ds}(0) \neq 0$. We will fix orientation of the parameter $s$ by the assumption $\frac{dx}{ds}(0) > 0$. Let $\widetilde{\Gamma}_0$ be the lifting of $W^u$ to the universal cover

$$(7.2) \qquad \widetilde{\Gamma}_0 = \{(\tilde{x}(s), u(s)),\ \tilde{x}(s) \in \mathbb{R}^1,\ u(s) \in \mathbb{R}^1\}.$$

Also, denote by $(x_s(t, \omega), v_s(t, \omega))$ the solution of (2.3) with initial data $x_s(0, \omega) = x(s), v_s(0, \omega) = u(s)$. Since for all $s$, the solutions $(x_s(t, \omega), v_s(t, \omega))$ converge exponentially fast to $(y_\omega(t), \dot{y}_\omega(t))$ as $t \to -\infty$, we can define the function

$$(7.3) \quad A(s) = \int_{-\infty}^{0} \left\{ \frac{1}{2}(v_s^2(t, \omega) - \dot{y}_\omega^2(t)) + (F(x_s(t, \omega), t) - F(y_\omega(t), t)) \right\} dt.$$

Since $W^u$ is almost surely a $C^{r-1}$ manifold, $A(s)$ is almost surely a $C^{r-1}$ function. Let

$$\bar{A}(x) = \min_{s:\, x(s)=x} A(s). \tag{7.4}$$

In the following, we will enumerate the shocks as in Section 5, except we number the main shock (which is not necessarily the strongest shock at a given time) as the zeroth shock.

Our next theorem describes the following picture. When we view the unstable manifold $W^u$ as a curve on the cylinder $\{x \in \mathbb{S}^1,\ u \in \mathbb{R}^1\}$, the two-sided minimizer divides $W^u$ into left and right pieces. It turns out that all shocks correspond to double folds of $W^u$, i.e., graphs of a multi-valued function. A single-valued function is obtained by introducing jump discontinuities which are vertical cuts on the double fold. These are the shocks in the solution. The end points of the cut define two points on $W^u$ with the same $x$-coordinate (namely the position of the shock) and the same value of the action $A$ in (7.3). If $x$ denotes the position of the shock, then the end points of the cut are $(x, u_+^\omega(x, 0))$ and $(x, u_-^\omega(x, 0))$. Except for the main shock, the one-sided minimizers starting from $(x, u_+^\omega(x, 0))$ and $(x, u_-^\omega(x, 0))$ approach the two-sided minimizer as $t \to -\infty$ from the same side. However, for the main shock, they approach the two-sided minimizer from different sides. We formulate this as:

THEOREM 7.2. *Fix arbitrary $t_0 \in \mathbb{R}^1$.*

I. *Let $(x(s), u(s)) \in W^u$. Now $(x(s), u(s))$ gives rise to a one-sided minimizer if and only if $A(s) = \bar{A}(x(s))$. With probability 1, $\bar{A}(x)$ is defined for all $x \in \mathbb{S}^1$; i.e., the minimum in (7.4) is attained. Moreover $\bar{A}$ is a continuous function on $\mathbb{S}^1$.*

II. *Let $x_i$ be the position of the $i^{\text{th}}$ shock, $i \geq 1$ (not the main shock!). Then there exists an interval $\ell_i = [\underline{s}_i, \overline{s}_i]$ such that*

$$\underline{s}_i = \min\{s:\, x(s) = x_i,\ A(s) = \bar{A}(x_i)\},$$
$$\overline{s}_i = \max\{s:\, x(s) = x_i,\ A(s) = \bar{A}(x_i)\}.$$

*Also, $\ell_i$ lies either to the left or to the right of the two-sided minimizer, i.e. $s = 0 \notin \ell_i$, $i \geq 1$ and $\ell_i \cap \ell_j = \emptyset$, $i \neq j$. If $\tilde{x}$ is the $x$-coordinate of points on the unstable manifold lifted to the universal cover, then $\tilde{x}(\underline{s}_i) = \tilde{x}(\overline{s}_i)$.*

III. *The main shock corresponds to the only point $z(t_0, \omega) \in \mathbb{S}^1$ such that there exist $s^{(1)} < 0$, $s^{(2)} > 0$ for which $A(s^{(1)}) = A(s^{(2)}) = \bar{A}(z(t_0, \omega))$, $x(s^{(1)}) = x(s^{(2)}) = z(t_0, \omega)$. Denote*

$$\underline{S} = \max\{s < 0:\, x(s) = z(t_0, \omega),\ A(s) = \bar{A}(z(t_0, \omega))\}$$
$$\overline{S} = \min\{s > 0:\, x(s) = z(t_0, \omega),\ A(s) = \bar{A}(z(t_0, \omega))\}.$$

*Then, $\tilde{x}(\overline{S}) - \tilde{x}(\underline{S}) = 1$.*

IV. Let $\Delta^\omega = [\underline{S}, \overline{S}] \setminus \bigcup_{i \geq 1} [\underline{s}_i, \overline{s}_i)$. Then for almost all $\omega$ the graph of $u_+^\omega(x, 0)$, $x \in \mathbb{S}^1$, coincides with $\{(x(s), u(s)), s \in \Delta^\omega\}$.

*Proof.* I. Clearly minimizers correspond to minima of $A$ when the $x$ coordinate is fixed. Since with probability 1, minimizers exist for all $x \in \mathbb{S}^1$, $\bar{A}(x)$ attains its minimum. Since the set of minimizers is closed, $\bar{A}(x)$ is continuous.

II. For every shock (except the main shock), denote by $\underline{s}_i$, $\overline{s}_i$ the values of the parameter $s$ corresponding to the left-most and right-most minimizers. Since both the left-most and right-most minimizers approach the two-sided minimizer from the same side, the interval $l_i = [\underline{s}_i, \overline{s}_i]$ does not contain $s = 0$. Since minimizers do not intersect, the intervals $\ell_i$ do not intersect. It follows from Theorem 7.1 that all shocks, except the main shock, have finite past history. Notice that, at the moment of creation of a shock, $\underline{s}_i = \overline{s}_i$. Hence, $\tilde{x}(\underline{s}_i) = \tilde{x}(\overline{s}_i)$. Since $\tilde{x}(\underline{s}_i) - \tilde{x}(\overline{s}_i)$ is a continuous function of time between merges and it takes only integer values, $\tilde{x}(\underline{s}_i) - \tilde{x}(\overline{s}_i) \equiv 0$ for all shocks except the main shock.

III. The main shock is the only shock for which the two extreme one-sided minimizers approach the two-sided minimizer from different sides. Thus $[\underline{S}, \overline{S}]$ has nonempty interior. Clearly the intervals $\ell_i$ constructed above belong to $[\underline{S}, \overline{S}]$. As a consequence of periodicity, we have $\tilde{x}(\overline{S}) = 1 + \tilde{x}(\underline{S})$.

IV. IV follows easily from I–III. $\square$

We next prove that for fixed time $t_0$ the number of shocks is finite. Consider time $t_0 - 1$. Although the position of the two-sided minimizer at time $t_0 - 1$ is not measurable with respect to $\mathcal{F}_{-\infty}^{t_0-1}$, the position of the main shock and the unstable manifold are measurable with respect to $\mathcal{F}_{-\infty}^{t_0-1}$. Consider the unstable manifold $W^u(t_0)$ at time $t_0$ as the image of $W^u(t_0 - 1)$ under the time-1 stochastic flow $G_1 = G_1^{\theta^{t_0-1}\omega}$:

(7.5)
$$W^u(t_0) = \{(x(s), v(s)) = G_1(y(s), w(s)), \ (y(s), w(s)) \in W^u(t_0 - 1), s \in \mathbb{R}^1\}.$$

Let $E$ be the event that there exists $s_0 \in \mathbb{R}^1$ such that $\frac{dx}{ds}(s_0) = \frac{d^2x}{ds^2}(s_0) = 0$.

LEMMA 7.3.
$$P\{E|\mathcal{F}_{-\infty}^{t_0-1}\} = 0$$

*for almost all conditions.*

*Proof.* Consider an arbitrary interval $[s_1, s_2]$, $s_1, s_2 \in \mathbb{R}^1$. Denote by $E_1$ the event that there exists $s_0 \in [s_1, s_2]$ such that $\frac{dx}{ds}(s_0) = \frac{d^2x}{ds^2}(s_0) = 0$. It is enough to show that $P\{E_1|\mathcal{F}_{-\infty}^{t_0-1}\} = 0$ for all $s_1, s_2$. Let $G_1 = (G_1^{(1)}, G_1^{(2)})$. We

have

$$\frac{dx}{ds} = \frac{\partial G_1^{(1)}}{\partial y}\frac{dy}{ds} + \frac{\partial G_1^{(1)}}{\partial w}\frac{dw}{ds}, \quad \frac{dv}{ds} = \frac{\partial G_1^{(2)}}{\partial y}\frac{dy}{ds} + \frac{\partial G_1^{(2)}}{\partial w}\frac{dw}{ds} \tag{7.6}$$

$$\frac{d^2x}{ds^2} = \frac{\partial^2 G_1^{(1)}}{\partial y^2}\left(\frac{dy}{ds}\right)^2 + 2\frac{\partial^2 G_1^{(1)}}{\partial y \partial w}\frac{dy}{ds}\frac{dw}{ds} + \frac{\partial^2 G_1^{(1)}}{\partial w^2}\left(\frac{dw}{ds}\right)^2$$

$$+ \frac{\partial G_1^{(1)}}{\partial y}\frac{d^2y}{ds^2} + \frac{\partial G_1^{(1)}}{\partial w}\frac{d^2w}{ds^2}.$$

Denote by

$$A(\omega) = \max_{s_1 \le s \le s_2}\left\{\left|\frac{\partial^k G_1^{(1,2)}}{\partial y^i \partial w^j}(y(s), w(s))\right|, \left|\frac{d}{ds}\left(\frac{\partial^k G_1^{(1,2)}}{\partial y^i \partial w^j}(y(s), w(s))\right)\right|, \right. \tag{7.7}$$

$$\left. 1 \le k \le 2, \; i+j=k, \; i \ge 0, \; j \ge 0\right\},$$

$$B(\omega) = \max_{s_1 \le s \le s_2}\left\{1, \left|\frac{d^2y}{ds^2}\right|, \left|\frac{d^2w}{ds^2}\right|\right\}.$$

Notice that $B(\omega)$ is measurable with respect to $\mathcal{F}_{-\infty}^{t_0-1}$. Take a small $\varepsilon > 0$ and divide the interval $[s_1, s_2]$ into subintervals of length $\varepsilon$. Denote by $s(j)$, $1 \le j \le \left[\frac{s_2-s_1}{\varepsilon}\right] + 2$, the end-points of the elements of the partition. Assume that there exists $s_0 \in [s(j), s(j+1)]$, such that $\frac{dx}{ds}(s_0) = \frac{d^2x}{ds^2}(s_0) = 0$. Then

$$\left|\frac{dx}{ds}(s(j))\right| = \left|\frac{dx}{ds}(s_0) + \frac{d^2x}{ds^2}(\xi)(s(j) - s_0)\right| = \left|\frac{d^2x}{ds^2}(\xi)(s(j) - s_0)\right|, \tag{7.8}$$

where $\xi \in (s(j), s_0)$. Denote by

$$V(\varepsilon, \omega) = \max_{s_1 \le s', s'' \le s_2, |s'-s''| \le \varepsilon}\left\{\left|\frac{dy}{ds}(s') - \frac{dy}{ds}(s'')\right|, \left|\frac{dw}{ds}(s') - \frac{dw}{ds}(s'')\right|, \right. \tag{7.9}$$

$$\left. \left|\frac{d^2y}{ds^2}(s') - \frac{d^2y}{ds^2}(s'')\right|, \left|\frac{d^2w}{ds^2}(s') - \frac{d^2w}{ds^2}(s'')\right|\right\}.$$

Obviously $V(\varepsilon, \omega)$ is measurable with respect to $\mathcal{F}_{-\infty}^{t_0-1}$ and $V(\varepsilon, \omega) \to 0$ as $\varepsilon \to 0$. Using $\left(\frac{dy}{ds}\right)^2 + \left(\frac{dw}{ds}\right)^2 = 1$ and (7.6) it is easy to show that for all $s_1 \le s', s'' \le s_2$, $|s' - s''| \le \varepsilon$,

$$\left|\frac{d^2x}{ds^2}(s') - \frac{d^2x}{ds^2}(s'')\right| \le A(\omega)(10 V(\varepsilon, \omega) + 2\varepsilon B(\omega) + 4\varepsilon). \tag{7.10}$$

Notice that $|\xi - s_0| \leq \varepsilon$. Hence, using (7.8), (7.9), we have

$$\text{(7.11)} \qquad \left|\frac{dx}{ds}(s(j))\right| \leq A(\omega)(10V(\varepsilon, \omega) + 2\varepsilon B(\omega) + 4\varepsilon)\varepsilon.$$

Fix arbitrary $\delta > 0$. We will show that the conditional probability that $s_0$ exists is less than $\delta$ for almost all conditions. Clearly there exists a random constant $K(\omega) > 0$ which is measurable with respect to $\mathcal{F}_{-\infty}^{t_0-1}$ such that

$$\text{(7.12)} \qquad P(A(\omega) > K(\omega)|\mathcal{F}_{-\infty}^{t_0-1}) < \frac{\delta}{2}$$

for almost all conditions. If $A(\omega) \leq K(\omega)$ then

$$\text{(7.13)} \qquad \left|\frac{dx}{ds}(s(j))\right| \leq T(\varepsilon, \omega) = K(\omega)(10V(\varepsilon, \omega) + 2\varepsilon B(\omega) + 4\varepsilon)\varepsilon,$$

where $T(\varepsilon, \omega)$ is measurable with respect to $\mathcal{F}_{-\infty}^{t_0-1}$.

Fix $\tilde{s} \in [s_1, s_2]$ and consider the random process $x_s(t) = \frac{\partial x}{\partial s}(t, \tilde{s})$, $v_s(t) = \frac{\partial v}{\partial s}(t, \tilde{s})$, where $(x(t,s), v(t,s)) = G_t^{\theta^{t_0-1}\omega}(y(s), w(s))$. Clearly $(x_s(t), v_s(t))$ satisfies the stochastic differential equations

$$\text{(7.14)} \qquad \dot{x}_s(t) = v_s(t), \ x_s(0) = \frac{dy}{ds}(\tilde{s}),$$

$$\dot{v}_s(t) = \sum_k f'_k(x(t, \tilde{s})) x_s(t) \dot{B}_k(t), \ v_s(0) = \frac{dw}{ds}(\tilde{s}).$$

It follows from Lemma B.9 that the joint probability distribution for $\left(\frac{\partial x}{\partial s}(1, \tilde{s}) = x_s(1), \frac{\partial v}{\partial s}(1, \tilde{s}) = v_s(1)\right)$ has density $p(x_s, v_s)$ which is uniformly bounded inside any compact set for all $\tilde{s} \in [s_1, s_2]$. If $A(\omega) \leq K(\omega)$ then, as follows from (7.6),

$$\text{(7.15)} \qquad \max_{s_1 \leq s \leq s_2} \max\left(\left|\frac{dx}{ds}\right|, \left|\frac{dv}{ds}\right|\right) \leq 2K(\omega).$$

Denote by

$$O(\omega) = \{(x_s, v_s) \in \mathbb{R}^2 : x_s^2 + v_s^2 \leq 8K^2(\omega)\},$$
$$\Pi_\varepsilon(\omega) = \{(x_s, v_s) \in \mathbb{R}^2 : |x_s| \leq T(\varepsilon, \omega)\}.$$

Let $R(\omega) = \max_{s_1 \leq \tilde{s} \leq s_2} \sup_{(x_s, v_s) \in O(\omega)} p(x_s, v_s)$. Then, for any $\tilde{s} \in [s_1, s_2]$,

$$\text{(7.16)}$$
$$P\left(\left(\frac{\partial x}{\partial s}(1, \tilde{s}), \frac{\partial v}{\partial s}(1, \tilde{s})\right) \in O(\omega) \cap \Pi_\varepsilon(\omega) \Big| \mathcal{F}_{-\infty}^{t_0-1}\right) \leq R(\omega)(4K(\omega))2T(\varepsilon, \omega).$$

Clearly, $R(\omega)$ is measurable with respect to $\mathcal{F}_{-\infty}^{t_0-1}$. Using (7.12), (7.16), we have

$$\text{(7.17)} \qquad P(E_1|\mathcal{F}_{-\infty}^{t_0-1}) \leq \frac{\delta}{2} + 8K(\omega)R(\omega)T(\varepsilon, \omega)\left(\left[\frac{s_2 - s_1}{\varepsilon}\right] + 1\right).$$

Choose $\varepsilon$ so small that

$$8K^2(\omega)R(\omega)\varepsilon\left(\left[\frac{s_2-s_1}{\varepsilon}\right]+1\right)(10V(\varepsilon,\omega)+2\varepsilon B(\omega)+4\varepsilon)\leq\frac{\delta}{2}.$$

Then, $P(E_1|\mathcal{F}_{-\infty}^{t_0-1})\leq\delta$ for almost all conditions. □

Lemma 7.3 easily implies the following theorem.

THEOREM 7.4. *Fix $t_0\in\mathbb{R}^1$. With probability* 1, *the number of shocks at time $t_0$ is finite.*

*Proof.* As above, consider the unstable manifold $W^u(t_0)$ parametrized by the arc-length parameter of the unstable manifold $W^u(t_0-1)$. Denote by $\underline{S}',\overline{S}'$ the values of the parameter corresponding to the main shock. For every shock at time $t_0$ (except the main shock) there exists an interval $\ell=[s',s'']\subset[\underline{S}',\overline{S}']$, such that $x(s')=x(s'')$. Thus there exists a point $\hat{s}\in(s',s'')$ for which $\frac{dx}{ds}(\hat{s})=0$. Notice that the intervals $\ell$ corresponding to different shocks do not intersect. If there are infinitely many shocks, then there exists an infinite sequence of $\hat{s}_i$'s in $[\underline{S}',\overline{S}']$, such that $\frac{dx}{ds}(\hat{s}_i)=0$. Let $s_0$ be an accumulation point for the sequence $\{\hat{s}_i\}$. Obviously, $\frac{dx}{ds}(s_0)=\frac{d^2x}{ds^2}(s_0)=0$. It follows from Lemma 7.3 that the conditional probability for the existence of such an $s_0$ is equal to zero. This immediately implies the theorem. □

## 8. The zero viscosity limit

In this section, we study the limit as $\varepsilon\to 0$ of the invariant measures for the viscous equation

$$(8.1)\qquad \frac{\partial u}{\partial t}+\frac{\partial}{\partial x}\Big(\frac{u^2}{2}\Big)=\frac{\varepsilon}{2}\frac{\partial^2 u}{\partial x^2}+\frac{\partial F}{\partial x}.$$

Under the same assumptions on $F$, it was proved in [S2] that for $\varepsilon>0$, there exists a unique invariant measure $\kappa_\varepsilon$ defined on the $\sigma$-algebra of Borel sets of $D_0$. Furthermore, as in the inviscid case studied in this paper, $\kappa_\varepsilon$ can be constructed as the probability distribution of an invariant functional

$$(8.2)\qquad u_\varepsilon^\omega(\cdot,0)=\Phi_0^\varepsilon(\omega)(\cdot)$$

such that $u_\varepsilon^\omega$ is a solution of (8.1) when the realization of the forces is given by $\omega$. The main result of this section is the following:

THEOREM 8.1. *With probability* 1,

$$u_\varepsilon^\omega(x,0)\to u^\omega(x,0)\,,$$

*for almost all* $x \in \mathbb{S}^1$, *as* $\varepsilon \to 0$. *More precisely, let* $x \in I(\omega) = \{y \in [0, 1]:$ *there exists a unique one-sided minimizer passing through* $(y, 0)\}$. *Then*

$$(8.3) \qquad u_\varepsilon^\omega(x, 0) \to u^\omega(x, 0) \qquad as \qquad \varepsilon \to 0.$$

As a simple corollary, we have:

THEOREM 8.2. *$\kappa_\varepsilon$ converges weakly to $\kappa$ as $\varepsilon \to 0$.*

Our proof of Theorem 8.1 relies on the Hopf-Cole transformation: $u_\varepsilon^\omega = -\varepsilon(\log \varphi)_x$ where $\varphi$ satisfies the stochastic heat equation

$$(8.4) \qquad \frac{\partial \varphi}{\partial t} = \frac{\varepsilon}{2}\frac{\partial^2 \varphi}{\partial x^2} - \frac{1}{\varepsilon}\varphi \circ F.$$

As explained in Appendix C, the product $F \circ \varphi$ should be understood in the Stratonovich sense. The solution of (8.4) has the Feynman-Kac representation

$$(8.5) \qquad \varphi(x, t) = \mathbb{E}\left\{e^{-\frac{1}{\varepsilon}\int_s^t F(x+\beta(\tau),\tau)d\tau}\varphi(x+\beta(s), s)\right\}$$

for $s < t$, where $\mathbb{E}$ denotes expectation with respect to the Wiener measure with variance $\varepsilon$, $\beta(t) = 0$, and

$$(8.6) \qquad \int_s^t F(x+\beta(\tau),\tau)d\tau = \sum_k F_k(x+\beta(t))B_k(t) - \sum_k F_k(x+\beta(s))B_k(s)$$
$$- \sum_k \int_s^t f_k(x+\beta(\tau))B_k(\tau)d\beta(\tau).$$

The integrals in (8.6) are understood in the Ito sense.

For $x, y \in \mathbb{R}^1$, $\tau_1 > \tau_2$, define

$$K_\varepsilon(x, \tau_1, y, \tau_2) = e^{\frac{1}{\varepsilon}\left(\sum_k F_k(x)B_k(\tau_1) - \sum_k F_k(y)B_k(\tau_2)\right)}$$
$$\times \int e^{\frac{1}{\varepsilon}\sum_k \int_{\tau_2}^{\tau_1} f_k(\beta(s))B_k(s)d\beta(s)} dW_{(y,\tau_2)}^{(x,\tau_1)}(\beta)$$

where $dW_{(y,\tau_2)}^{(x,\tau_1)}(\beta)$ is the probability measure defined by the Brownian bridge: $\beta(\tau_1) = x$, $\beta(\tau_2) = y$, with variance $\varepsilon$. Using (8.5), for $s < t$, we can write the solution of (8.1) as

$$(8.7) \qquad u^\varepsilon(x, t) = -\varepsilon \frac{\int_0^1 \frac{\partial}{\partial x} M(x, t, y, s) e^{-\frac{1}{\varepsilon}h^\varepsilon(y,s)} dy}{\int_0^1 M(x, t, y, s) e^{-\frac{1}{\varepsilon}h^\varepsilon(y,s)} dy}$$

where $h^\varepsilon(y, s) = \int_0^y u^\varepsilon(z, s)dz$, and

$$M(x, \tau_1, y, \tau_2) = \sum_{m=-\infty}^{\infty} K_\varepsilon(x, \tau_1, y+m, \tau_2)$$

for $x, y \in [0, 1]$. $M$ is the transfer matrix for Brownian motion on the circle $\mathbb{S}^1$.

Define also
$$A(x, \tau_1, y, \tau_2) = \inf_{\xi(\tau_1)=x, \xi(\tau_2)=y} \mathcal{A}_{\tau_2, \tau_1}(\xi).$$

LEMMA 8.1. *For almost every $\omega \in \Omega$, there exists a $\tau = \tau(\omega) > 0$, and $C_0 = C_0(\omega) > 0$, such that*

$$(8.8) \quad \frac{1}{g_m} \leq K_\varepsilon(x, \tau, z+m, 0) e^{\frac{1}{\varepsilon} A(x, \tau, z+m, 0)} \leq g_m$$

*for all $x, z \in [0,1]$, $m \in \mathbb{N}$, where*

$$(8.9) \quad g_m = \begin{cases} C_1(\omega), & \text{if } (|m|+1)\|F\|_\tau \leq C_0(\omega) \\ C_2(\omega) e^{\frac{8(|m|+1)^2}{\varepsilon} \|F\|_\tau^2 \tau}, & \text{if } (|m|+1)\|F\|_\tau > C_0(\omega). \end{cases}$$

$C_1(\omega)$ *and* $C_2(\omega)$ *are constants depending only on $\omega$ and $\|F\|_\tau$ is defined as in Appendix B.*

*Proof.* We will assume $m \geq 0$. Let $\gamma: [0, \tau] \to R^1$ be a minimizer such that $\gamma(0) = z + m$, $\gamma(\tau) = x$, and
$$A(x, \tau, z+m, 0) = \mathcal{A}_{0,\tau}(\gamma).$$
Now, $\gamma$ satisfies the Euler-Lagrange equation:

$$(8.10) \quad \int_0^\tau \dot\gamma(s) d\eta(s) - \int_0^\tau \sum_k B_k(s) \Big(f_k'(\gamma(s))\dot\gamma(s) ds + f_k(\gamma(s))d\eta(s)\Big) = 0$$

for test functions $\eta$ on $[0, \tau]$. Performing a change of variable $\beta = \gamma + \sqrt{\varepsilon}\eta$ in the functional integral in $K_\varepsilon$, we obtain, using (8.10) and the Cameron-Martin-Girsanov formula:

$$(8.11) \quad K_\varepsilon(x, \tau, z+m, 0) = e^{-\frac{1}{\varepsilon} A(x, \tau, z+m, 0)} \mathbb{E}_\eta e^{\frac{1}{\varepsilon} H}$$

where the exponent $H$ is given by
(8.12)
$$H = H_1 + \sqrt{\varepsilon} H_2,$$
$$H_1 = \int_0^\tau \sum_k \{f_k(\gamma + \sqrt{\varepsilon}\eta) - f_k(\gamma) - f_k'(\gamma)\sqrt{\varepsilon}\eta\} B_k(s) \dot\gamma(s) ds,$$
$$H_2 = \int_0^\tau \sum_k B_k(s) \{f_k(\gamma + \sqrt{\varepsilon}\eta) - f_k(\gamma)\} d\eta(s).$$

In (8.11), $\mathbb{E}_\eta$ denotes expectation with respect to the standard Brownian bridge $\eta(0) = 0$, $\eta(\tau) = 0$. □

We now estimate $H$. A simple Taylor expansion to second order gives:

$$(8.13) \quad |H_1| \leq \varepsilon \|F\|_\tau \max_{0 \leq s \leq \tau} |\dot\gamma(s)| \int_0^\tau \eta^2(s) ds.$$

Using Lemma B.1, we get for $\tau = \tau(\omega)$,

$$(8.14) \quad |H_1| \leq \varepsilon \|F\|_\tau \frac{C(\omega)(|m|+1)}{\tau} \int_0^\tau \eta^2(s)ds \leq C(\omega)\varepsilon\|F\|_\tau(|m|+1)(\eta^*)^2,$$

where $\eta^* = \max_{0 \leq s \leq \tau} \eta(s)$. For $H_2$, we use the mean value theorem to write

$$(8.15) \quad \frac{1}{\sqrt{\varepsilon}} H_2 = \int_0^\tau \sum_k B_k(s) f'_k\left(\gamma + \sqrt{\varepsilon}\theta_k \eta\right) \eta(s) d\eta(s) = H_{21} + H_{22}$$

where $\theta_k \in [0,1]$, $H_{21} = \frac{1}{\sqrt{\varepsilon}} H_2 - H_{22}$ with

$$(8.16) \quad H_{22} = \frac{\alpha}{2} \int_0^\tau \left(\sum_k B_k(s) f'_k\left(\gamma + \sqrt{\varepsilon}\theta_k \eta\right)\right)^2 \eta^2(s) ds \leq \frac{\alpha}{2} \|F\|_\tau^2 \tau (\eta^*)^2.$$

We will choose the value of $\alpha$ later ($\alpha = 3$ will suffice). Note that (8.14) and (8.16) can be combined to give:

$$\frac{1}{\varepsilon}|H_1| + |H_{22}| \leq C(\omega)\|F\|_\tau(|m|+1)(\eta^*)^2.$$

The constant $C(\omega)$ is changed to a different value in the last step. Now we have
(8.17)
$$\mathbb{E}_\eta e^{\frac{1}{\varepsilon}H} = \mathbb{E}_\eta e^{\frac{1}{\varepsilon}H_1 + H_{21} + H_{22}} \leq \left(\mathbb{E}_\eta e^{C(\omega)\|F\|_\tau(|m|+1)(\eta^*)^2}\right)^{1/2} \left(\mathbb{E}_\eta e^{2H_{21}}\right)^{1/2}.$$

Using the fact that for $a > 0$

$$(8.18) \quad P\{\eta^* > a\} \leq \frac{C}{\sqrt{\pi\tau}} \int_a^{+\infty} e^{-\frac{\lambda^2}{2\tau}} d\lambda$$

we have

$$(8.19) \quad \mathbb{E}_\eta e^{C(\omega)\|F\|_\tau(|m|+1)(\eta^*)^2} \leq \text{Constant}$$

if

$$(8.20) \quad C(\omega)\|F\|_\tau(|m|+1) < \frac{1}{2\tau}.$$

For the second factor in (8.17), we use the inequality (see [McK])

$$(8.21) \quad P\{H_{21} > \beta\} < e^{-\alpha\beta}.$$

Then

$$(8.22) \quad \mathbb{E}_\eta e^{2H_{21}} \leq \sum_{N=0}^\infty e^{2(N+1)} P\{H_{21} > N\} \leq \sum_{N=1}^\infty e^{2(N+1)-\alpha N} < +\infty$$

if we choose $\alpha > 2$.

When $m$ violates (8.20), we estimate $H_1$ using Lemma B.1:

(8.23)
$$H_1 = \sqrt{\varepsilon} \int_0^\tau \sum_k B_k(s)\Big(f'_k\big(\gamma + \sqrt{\varepsilon}\theta_k\eta\big) - f'_k(\gamma)\Big)\eta(s)\dot{\gamma}(s)ds,$$

(8.24)
$$|H_1| \leq 2\sqrt{\varepsilon}\,\|F\|_\tau \frac{(|m|+1)}{\tau} \int_0^\tau |\eta(s)|ds \leq 2\sqrt{\varepsilon}\,\|F\|_\tau (|m|+1)\eta^*.$$

Hence

(8.25)
$$\mathbb{E}_\eta \left(e^{\frac{1}{\varepsilon}H_1}\right)^2 \leq C(\omega) \int_0^{+\infty} e^{\frac{4(|m|+1)}{\sqrt{\varepsilon}}\|F\|_\tau \sqrt{\tau}\lambda - \frac{\lambda^2}{2}} d\lambda$$
$$\leq C(\omega) e^{\frac{16(|m|+1)^2}{\varepsilon}\|F\|_\tau^2 \tau}.$$

As before we have then
$$\mathbb{E}_\eta e^{\frac{1}{\varepsilon}H} \leq C(\omega) e^{\frac{8(|m|+1)^2}{\varepsilon}\|F\|_\tau^2 \tau}.$$

These give the upper bounds.

Similarly we can prove the same bounds for $(\mathbb{E}_\eta e^{\frac{1}{\varepsilon}H})^{-1}$. This completes the proof of Lemma 8.1. $\square$

It is easy to see that for fixed $\tau$, $z$ and $m$, $A(x,\tau,z+m,0)$ is differentiable at $x$, if and only if there exists a unique minimizer $\gamma$ such that $A(x,\tau,z+m,0) = \mathcal{A}_{0,\tau}(\gamma)$ and $\gamma(\tau) = x$, $\gamma(0) = z+m$. In this case we have

(8.26)
$$\dot{\gamma}(\tau) = \frac{\partial}{\partial x}A(x,\tau,z+m,0).$$

When the minimizer is not unique, $A(x,\tau,z+m,0)$ has both left and right derivatives. Moreover
$$D_x^+ A(x,\tau,z+m,0) = \dot{\gamma}_+(\tau),$$
$$D_x^- A(x,\tau,z+m,0) = \dot{\gamma}_-(\tau)$$

where $\gamma_+$ and $\gamma_-$ are the right-most and left-most minimizers. In either case, let us define
$$v(x,z+m,\tau) = D_x^- A(x,\tau,z+m,0).$$

LEMMA 8.2. *The following inequality holds*:

(8.27)
$$\left|\frac{1}{K_\varepsilon(x,\tau,z+m,0)}\left(\varepsilon\frac{\partial K_\varepsilon}{\partial x}(x,\tau,z+m,0) + v(x,z+m,\tau)K_\varepsilon(x,\tau,z+m,0)\right)\right|$$
$$\leq \sqrt{\varepsilon}\,\|F\|_\tau g_m$$

where $g_m$ is as defined in Lemma 8.1.

*Proof.* For simplicity, we will write

$$K_m = K_\varepsilon(x, \tau, z+m, 0), \quad v_m = v(x, z+m, \tau).$$

A straightforward computation gives
(8.28)
$$\varepsilon \frac{\partial K_m}{\partial x} = e^{-\frac{1}{\varepsilon} \sum_k F_k(x) B_k(\tau)} \int G(\beta) e^{\frac{1}{\varepsilon} \sum_k \int_0^\tau f_k(\beta(s)) B_k(s) d\beta(s)} dW^{(x,\tau)}_{(z+m,0)}(\beta)$$

where

(8.29)
$$G(\beta) = -\sum_k f_k(x) B_k(\tau) - \frac{x-(z+m)}{\tau}$$
$$+ \sum_k \int_0^\tau B_k(s) \left[ f_k'(\beta) \frac{s}{\tau} d\beta + f_k(\beta) \frac{1}{\tau} ds \right],$$

and

$$\varepsilon \frac{\partial K_m}{\partial x} + v_m K_m = e^{-\frac{1}{\varepsilon} \sum_k F_k(x) B_k(\tau)}$$
$$\times \int (G(\beta) - G(\gamma_-)) e^{\frac{1}{\varepsilon} \sum_k \int_0^\tau f_k(\beta(s)) B_k(s) d\beta(s)} dW^{(x,\tau)}_{(z+m,0)}(\beta).$$

Performing a change of variables $\beta = \gamma_- + \sqrt{\varepsilon}\eta$, we get

(8.30) $$\varepsilon \frac{\partial K_m}{\partial x} + v_m K_m = e^{-\frac{1}{\varepsilon} A(x,\tau,z+m,0)} \mathbb{E}_\eta \left( (G(\gamma_- + \sqrt{\varepsilon}\eta) - G(\gamma_-)) e^{\frac{1}{\varepsilon} H} \right)$$

where $H$ is defined as before. Write $G(\gamma_- + \sqrt{\varepsilon}\eta) - G(\gamma_-)$ as

$$G(\gamma_- + \sqrt{\varepsilon}\eta) - G(\gamma_-)$$
$$= \frac{\sqrt{\varepsilon}}{\tau} \sum_k \left[ \int_0^\tau B_k(\tau) f_k''(\gamma_- + \theta_k \sqrt{\varepsilon}\eta) s\eta(s) \dot{\gamma}_-(s) ds \right.$$
$$\left. + \int_0^\tau B_k(s) f_k'(\gamma_- + \sqrt{\varepsilon}\eta) s d\eta(s) + \int_0^\tau B_k(\tau) f_k'(\gamma_- + \theta_k \sqrt{\varepsilon}\eta) \eta(s) ds \right].$$

We can then follow the steps in the proof of Lemma 8.1 to establish (8.27). □

*Remark.* The estimates in Lemmas 8.1 and 8.2 are proved for the time interval $[0, \tau]$. We see easily that they hold for arbitrary intervals of the type $[t, t+\tau]$ and $[t-\tau, t]$ by choosing suitable $\tau$ which in general depend on $t$. For $t$ in a compact set, we can choose $\tau$ to be independent of $t$ such that (8.8) holds.

Our next lemma gives uniform estimates of $u^\varepsilon$.

LEMMA 8.3. *There exist positive constants $\varepsilon_0(\omega,t)$, $C(\omega,t)$, such that*

$$|u_\omega^\varepsilon(x,t)| \leq C(\omega,t) \tag{8.31}$$

*for $x \in [0,1]$, $0 < \varepsilon \leq \varepsilon_0(\omega,t)$. Furthermore, $\varepsilon_0(\,\cdot\,,t)$ and $C(\,\cdot\,,t)$ are stationary random processes in $t$.*

*Proof.* The basic idea is to use the fact that for $\varepsilon$ small, the functional integral is concentrated near minimizers whose velocities are estimated in Lemma B.1. We will prove (8.31) for $t = 0$ by working on the time interval $[-\tau, 0]$ where $\tau$ is as defined in Lemma 8.1. It will be clear that the proof works with little change for arbitrary $t$. □

Let $N = \frac{C_0(\omega)}{\|F\|_\tau}$, where $C_0(\omega)$ is as defined in Lemma 8.1. Denote $A^*(x,z,\tau) = \inf_{m \in \mathbb{N}} A(x, 0, z+m, -\tau)$ for $x, z \in [0,1]$. It is easy to see that for $\tau \ll 1$,

$$A(x, 0, z+m, -\tau) - A^*(x, z, \tau) \geq \frac{1}{3}\frac{(|m|+1)^2}{\tau} - C(\omega) \tag{8.32}$$

for $|m| > N$.

Again for simplicity of notation, we will denote $K_m = K_\varepsilon(x, 0, y+m, -\tau)$, $\mu(dy) = e^{-\frac{1}{\varepsilon}h^\varepsilon(y,0)}dy$. Using Lemma 8.2 and (8.7), we have

$$|u^\varepsilon(x,0)| \leq I_1 + \sqrt{\varepsilon}\|F\|_\tau I_2$$

where

$$I_1 = \frac{\sum_m \int_0^1 v_m K_m \mu(dy)}{\sum_m \int_0^1 K_m \mu(dy)}, \tag{8.33}$$

$$I_2 = \frac{\sum_m \int_0^1 g_m K_m \mu(dy)}{\sum_m \int_0^1 K_m \mu(dy)}.$$

For $|m| > N$, we can use (8.32) and Lemma 8.1 to get

$$K_m \leq e^{-\frac{1}{\varepsilon}A(x,0,y+m,-\tau)}g_m$$

$$\leq C(\omega)e^{-\frac{1}{\varepsilon}A^*(x,y,\tau)}e^{-\frac{1}{3\varepsilon}\frac{(|m|+1)^2}{\tau} + \frac{8(|m|+1)^2}{\varepsilon}\|F\|_\tau^2\tau}$$

$$\leq C(\omega)^{-\frac{1}{\varepsilon}A^*(x,y,\tau)}e^{-\frac{1}{4\varepsilon}\frac{(|m|+1)^2}{\tau}}$$

if $\tau$ is small enough. Hence we get, using the fact that $|v_m| \leq \frac{C(\omega)(|m|+1)}{\tau}$,

$$|I_1| \leq \frac{\sum\limits_{|m|<N} \int_0^1 |v_m| K_m \mu(dy)}{\sum\limits_m \int_0^1 K_m \mu(dy)} + \frac{\sum\limits_{|m|>N} \int_0^1 |v_m| K_m \mu(dy)}{\sum\limits_m \int_0^1 K_m \mu(dy)}$$

$$\leq C(\omega) N \cdot \frac{1}{\tau} + C(\omega) \frac{\left(\sum\limits_{|m|>N} \frac{C(\omega)(|m|+1)}{\tau} e^{-\frac{1}{4\varepsilon} \frac{(|m|+1)^2}{\tau}}\right) \int_0^1 e^{-\frac{1}{\varepsilon} A^*(x,y,\tau)} \mu(dy)}{\int_0^1 e^{-\frac{1}{\varepsilon} A^*(x,y,\tau)} \mu(dy)}$$

$$\leq C(\omega).$$

In the last step we used the fact that $\tau$ depends only on $\omega$. Similarly

$$|I_2| \leq C(\omega) + \frac{\sum\limits_{|m|>N} \int_0^1 g_m K_m \mu(dy)}{\sum\limits_m \int_0^1 K_m \mu(dy)}$$

$$\leq C(\omega)$$

where we used

$$K_m g_m \leq C(\omega) e^{-\frac{1}{\varepsilon} A^*(x,y,\tau)} e^{-\frac{1}{4\varepsilon} \frac{(|m|+1)^2}{\tau}}$$

for small enough $\tau$. This completes the proof of Lemma 8.3. $\square$

Define for $C > 0$

(8.34) $$Q_C = \Big\{ h \in \text{Lip}[0,1], \quad \text{such that}$$

$$h(y) = \int_0^y u(z) dz, \ |u| \leq C, \ \int_0^1 u(z) dz = 0 \Big\}.$$

Take $x \in I(\omega)$. Denote the unique minimizer that passes through $(x,0)$ by $\xi^*$. For $h \in Q_C$, $T < 0$, define the modified action as

(8.35) $$\mathcal{A}_{T,0}^h(\xi) = \mathcal{A}_{T,0}(\xi) + h(\xi(T))$$

and denote by $\xi_h^{**}$ the minimizer of $\mathcal{A}_{T,0}^h$. Obviously $\xi_h^{**}$ in general depends on $h$ and $T$.

LEMMA 8.4. *Fix a constant $C > 0$. For any $\delta > 0$, there exists $T^* = T^*(\delta) < 0$, such that*

(8.36) $$|\dot{\xi}_h^{**}(0) - \dot{\xi}^*(0)| < \frac{\delta}{2}$$

*for $T < T^*$ and all $h \in Q_C$ ($T^*$ in general depends on $C$).*

*Proof.* Assume to the contrary that there exists a sequence $T_j \to -\infty$, $h_j \in Q_C$, such that

$$|\dot{\xi}_{h_j}^{**}(0) - \dot{\xi}^*(0)| \geq \frac{\delta}{2}.$$

Then from Lemma 3.3, the $\{\dot{\xi}^{**}_{h_j}(0)\}$'s are uniformly bounded and we can choose a subsequence, still denoted by $\{\xi^{**}_{h_j}\}$, such that $\xi^{**}_{h_j}$ converges (uniformly on compact sets of $(-\infty, 0]$ and $\dot{\xi}^{**}_{h_j}(0) \to \dot{\tilde{\xi}}(0)$ to a limiting path $\tilde{\xi}$ defined on $(-\infty, 0]$. From Lemma 3.6, $\tilde{\xi}$ is also a one-sided minimizer. Since $\tilde{\xi}(0) = x$, and
$$|\dot{\tilde{\xi}}(0) - \dot{\xi}^*(0)| \geq \frac{\delta}{2},$$
this violates the assumption that there exists a unique one-sided minimizer passing through $(x, 0)$. □

LEMMA 8.5. *Fix a constant $C > 0$. There exists a function $\alpha(\cdot)$ defined on $(0, +\infty)$, $\alpha > 0$, with the following properties: For any $\delta > 0$, one can find a $T^* = T^*(\delta) < 0$, such that for any path $\xi$ defined on $[T^*, 0]$, with $\xi(0) = x$, the inequality $|\dot{\xi}(0) - \dot{\xi}^*(0)| > \delta$ implies*

$$|\mathcal{A}^h_{T^*,0}(\xi) - \mathcal{A}^h_{T^*,0}(\xi^{**}_h)| > \alpha(\delta) \tag{8.37}$$

*for all $h \in Q_C$.*

*Proof.* Assume to the contrary that there exist a $\delta > 0$, and a sequence $\{T_j\}$, $T_j \to -\infty$, $h_j \in Q_C$, and $\xi_j$ defined on $[-T_j, 0]$, such that $|\dot{\xi}_j(0) - \dot{\xi}^*(0)| > \delta$, and

$$\left|\mathcal{A}^{h_j}_{T^*_j,0}(\xi_j) - \mathcal{A}^{h_j}_{T^*_j,0}(\xi^{**}_{h_j})\right| < \frac{1}{j}. \tag{8.38}$$

From the estimates proved in Section 3, $\{\dot{\xi}_j(0)\}$ are uniformly bounded. Therefore we can choose a subsequence, still denoted by $\{\xi_j\}$, such that $\xi_j$ converges (uniformly on compact sets of $(-\infty, 0]$ and $\dot{\xi}_j(0) \to \dot{\tilde{\xi}}(0)$) to $\tilde{\xi}$ defined on $(-\infty, 0]$. From (8.38), $\tilde{\xi}$ is also a one-sided minimizer. Since $\tilde{\xi}(0) = x$, $|\dot{\tilde{\xi}}(0) - \dot{\xi}^*(0)| > \delta$, we arrive at a contradiction to the assumption that there exists a unique one-sided minimizer passing through $(x, 0)$. □

Now we are ready to prove Theorem 8.1.

*Proof of Theorem* 8.1. Fix an $x \in I(\omega)$. Denote by $\xi^*$ the unique one-sided minimizer passing through $(x, 0)$. Take $\delta > 0$. From Lemmas 8.4 and 8.5 we can find a $T^* < 0$, such that (8.36) and (8.37) hold.

Let $n$ be a sufficiently large integer (depending only on $\omega$ and $T^*$), such that the estimates in Lemma 8.1 hold on the intervals $[(k+1)s, ks]$ where $s = \frac{T^*}{n}$, $k = 0, 1, \ldots, n-1$. Using Lemma 8.2, we have $\tau = -s$, for,

$$u^\varepsilon(x, 0) = -\frac{\varepsilon \int_0^1 \frac{\partial}{\partial x} M(x, 0, y, T^*) \mu(dy)}{\int_0^1 M(x, 0, y, T^*) \mu(dy)}$$

where
$$\mu(dy) = -\frac{1}{\varepsilon}h(y),\ h(y) = \int_0^y u^\varepsilon(z,T^*)dz\,.$$

Hence

$$u^\varepsilon(x,0) = \frac{-\varepsilon\int_0^1\int_0^1 \frac{\partial}{\partial x}M(x,0,z_1,s)M(z_1,s,y,T^*)dz_1\mu(dy)}{\int_0^1\int_0^1 M(x,0,z_1,s)M(z_1,s,y,T^*)dz_1\mu(dy)}$$

$$= \frac{\int_0^1\int_0^1 \Big(\sum_m v(x,0,z_1+m,s)K_\varepsilon(x,0,z_1+m,s)\Big)M(z_1,s,y,T^*)dz_1\mu(dy)}{\int_0^1\int_0^1 \Big(\sum_m K_\varepsilon(x,0,z_1+m,s)\Big)M(z_1,s,y,T^*)dz_1\mu(dy)}$$

$$+ O(1)\sqrt{\varepsilon}\|F\|_\tau \frac{\int_0^1\int_0^1 \Big(\sum_m g_m K_\varepsilon(x,0,z+m,s)\Big)M(z_1,s,y,T^*)dz_1\mu(dy)}{\int_0^1\int_0^1 \Big(\sum_m K_\varepsilon(x,0,z+m,s)\Big)M(z_1,s,y,T^*)dz_1\mu(dy)}$$

$$= I_3 + O(1)\sqrt{\varepsilon}\|F\|_\tau I_4$$

where $O(1)$ denotes a uniformly bounded quantity. As in the proof of Lemma 8.3, we can show:

$$|I_4| \leq C(\omega)N = \frac{C(\omega)C_0(\omega)}{\|F\|_\tau}\,.$$

For $z_k, z_{k+1} \in [0,1]$ we have, using Lemma 8.1,

$$M(z_k, ks, z_{k+1}, (k+1)s)$$
$$= \sum_{|m|<N} K_\varepsilon(z_k, ks, z_{k+1}+m, (k+1)s)$$
$$+ \sum_{|m|>N} K_\varepsilon(z_k, ks, z_{k+1}+m, (k+1)s)$$
$$\leq e^{-\frac{1}{\varepsilon}A^*(z_k,ks,z_{k+1},(k+1)s)}\Big(C(\omega)N + \sum_{|m|>N} e^{-\frac{(|m|+1)^2}{3\varepsilon}+\frac{8(|m|+1)^2}{\varepsilon}\|F\|_\tau^2\tau}\Big)$$
$$\leq C(\omega)Ne^{-\frac{1}{\varepsilon}A^*(z_k,ks,z_{k+1},(k+1)s)}\,.$$

Letting $x = z_0$, $y = z_n$, we obtain for fixed $\{z_0, z_1, z_2, \ldots, z_n\}$,

$$\prod_{k=1}^{n-1} M(z_k, ks, z_{k+1}, (k+1)s) \leq (C(\omega)N)^n e^{-\frac{1}{\varepsilon}\sum_{k=0}^{n-1} A^*(z_k,ks,z_{k+1},(k+1)s)}\,.$$

Denote by $\int' dz_1 \sum'_m$ and $\int'' dz_1 \sum''_m$ summation and integration over the sets of $(z_1, m)$ such that $|v_m - \dot{\xi}^*(0)| > \delta$ and $|v_m - \dot{\xi}^*(0)| < \delta$ respectively, where $v_m = v(x, 0, z_1 + m, s)$. From Lemma 8.4, the second sum and integral $\int'' dz_1 \sum''_m$ cover the set $|v_m - \dot{\xi}_h^{**}(0)| < \frac{\delta}{2}$. We have

$$I_5 = \int' dz_1 \sum_m{}' K_\varepsilon(x,0,z_1+m,s) M(z_1,s,y,T^*)\mu(dy)$$

$$= \int' dz_1 \sum_m{}' K_\varepsilon(x,0,z_1+m,s)$$

$$\times \int_{[0,1]^{n-1}} \prod_{k=1}^{n-1} M(z_k,ks,z_{k+1},(k+1)s) dz_2 \ldots dz_n \mu(dy)$$

$$\leq (C(\omega)N)^n \int' dz_1 \sum_m{}'$$

$$\times \int_{[0,1]^{n-1}} e^{-\frac{1}{\varepsilon}(A(x,0,z_1+m,s) + \sum_{k=1}^{n-1} A^*(z_k,ks,z_{k+1},(k+1)s) + h(z_n))} dz_2 \ldots dz_n$$

$$\leq (C(\omega)N)^n \int' dz_1 \sum_m{}' \int_{[0,1]^{n-1}} e^{-\frac{1}{\varepsilon} \sum_{k=0}^{n-1} A^*(z_k,ks,z_{k+1},(k+1)s) + h(z_n)} dz_2 dz_3 \ldots dz_n$$

$$\leq (C(\omega)N)^n e^{-\frac{1}{\varepsilon}(\mathcal{A}_{T^*,0}^h(\xi_h^{**}) + \alpha(\delta))}.$$

In the last step, we used (8.39). On the other hand, there exists a $\delta_2 > 0$, such that if $|z - \xi_h^{**}(s)| < \delta_2$, then

$$|v(x,0,z,s) - \dot{\xi}_h^{**}(0)| < \frac{\delta}{2}.$$

Choose a $\delta_1 > 0$, such that $\delta_1 < \delta_2$, $\delta_1 < \frac{\alpha(\delta)}{2nC_3(\omega)}$, with $C_3(\omega)$ to be defined later. Using Lemma 8.4 we get

$$I_6 = \int'' dz_1 \sum_m{}'' K_\varepsilon(x,0,z_1+m,s) M(z_1,s,y,T^*)\mu(dy)$$

$$\geq \int_{|z_k - \xi_h^{**}(ks)| < \delta_1} \prod_{k=1}^{n} K_\varepsilon(z_{k-1},(k-1)s,z_k,ks) e^{-\frac{1}{\varepsilon}h(z_n)} dz_1 \ldots dz_n$$

$$\geq \frac{1}{C(\omega)^n} \int_{|z_k - \xi_h^*(ks)| < \delta_1} e^{-\frac{1}{\varepsilon}\left(\sum_{k=1}^{n} A^*(z_{k-1},(k-1)s,z_k,ks) + h(z_n)\right)} dz_1 \ldots dz_n.$$

It is easy to see that if $|z_k - \xi_h^*(ks)| < \delta_1$, $|z_{k-1} - \xi_h^*((k-1)s)| < \delta_1$, then

$$|A^*(z_{k-1},(k-1)s,z_k,ks) - \mathcal{A}_{ks,(k-1)s}(\xi_h^{**})| \leq C_3(\omega)\delta_1,$$

and if $|z_n - \xi_h^{**}(T^*)| < \delta_1$, then

$$|h(z_n) - h(\xi_h^{**}(T^*))| \leq C_3(\omega)\delta_1,$$

with $C_3(\omega)$ defined by the above estimates. Hence we have

$$I_6 \geq \frac{1}{C(\omega)^n} \delta_1^n e^{-\frac{1}{\varepsilon}(\mathcal{A}_{T^*,0}^h(\xi_h^{**}) + nC_3(\omega)\delta_1)}.$$

Therefore

$$\left|\frac{I_5}{I_6}\right| \leq C(\omega)^{2n} \cdot N^n \delta_1^{-n} e^{-\frac{1}{2\varepsilon}\alpha(\delta)}.$$

Similarly, if we define

$$I_7 = \int' dz_1 \sum_m{}' v(x,0,z_1+m,s) K_\varepsilon(x,0,z_1+m,s) M(z_1,s,y,T^*) \mu(dy),$$

$$I_8 = \int'' dz_1 \sum_m{}'' v(x,0,z_1+m,s) K_\varepsilon(x,0,z_1+m,s) M(z_1,s,y,T^*) \mu(dy),$$

then we can also get

$$\left|\frac{I_7}{I_6}\right| \leq C(\omega)^{2n} N^n \delta_1^{-n} e^{-\frac{1}{2\varepsilon}\alpha(\delta)}.$$

Finally we obtain

$$\begin{aligned}
|u^\varepsilon(x,0) - \dot\xi^*(0)| &\leq \left|\frac{I_7 + I_8 - \xi^*(0)(I_5 + I_6)}{I_5 + I_6}\right| + \sqrt{\varepsilon} C(\omega) \|F\|_\tau \\
&\leq \left|\frac{I_8 - \xi^*(0) I_6}{I_6}\right| + \frac{|I_7|}{I_6} + |\dot\xi(0)|\frac{I_5}{I_6} + \sqrt{\varepsilon} \|F\|_\tau C(\omega) \\
&\leq \delta + C(\omega)^{2n} N^n \delta_1^{-n} e^{-\frac{1}{2\varepsilon}\alpha(\delta)} + \sqrt{\varepsilon} \|F\|_\tau C(\omega) \\
&\leq \delta + \delta = 2\delta
\end{aligned}$$

if we choose $\varepsilon$ sufficiently small. This completes the proof of Theorem 8.1. □

## Appendix A. Proof of Lemma 2.1 for the random case

In this appendix we comment on the proof of Lemma 2.1 for the random case. Let $F^\delta(x,t) = \sum_{k=1}^\infty F_k(x) \dot B_k^\delta(t)$, where $B_k^\delta$ is the standard mollification of $B_k$. Denote by $u^\delta(x,t)$ the unique entropy solution of

$$(\text{A.1}) \qquad \frac{\partial u}{\partial t} + \frac{\partial}{\partial x}\left(\frac{u^2}{2}\right) = \frac{\partial F^\delta}{\partial x}$$

with the initial data $u(x,t_0) = u_0(x)$. We will assume that $\|u_0\|_{L^\infty} \leq$ Const, $\int_0^1 u_0(z) dz = 0$. From classical results [Li] we know that $u^\delta(x,t)$ is given by

$$(\text{A.2}) \qquad u^\delta(x,t) = \frac{\partial}{\partial x} \inf_{\xi:\xi(t)=x} \left\{\mathcal{A}_{t_0,t}^\delta(\xi) + \int_0^{\xi(t_0)} u_0(z) dz\right\},$$

where

(A.3)
$$\mathcal{A}^\delta_{t_0,t}(\xi) = \int_{t_0}^t \frac{1}{2}\dot\xi(s)^2 ds + \int_{t_0}^t \sum_{k=1}^\infty F_k(\xi(s))dB_k^\delta(s)$$
$$= \sum_{k=1}^\infty F_k(\xi(t))(B_k^\delta(t) - B_k^\delta(t_0))$$
$$+ \int_{t_0}^t \left(\frac{1}{2}\dot\xi(s)^2 - \sum_{k=1}^\infty f_k(\xi(s))\dot\xi(s)(B_k^\delta(s) - B_k^\delta(t_0))\right) ds.$$

It is easy to see that the boundary terms at the right hand side of (A.3) resulted from integration by parts do not affect the variational formula (A.2) and can be neglected. Denote

(A.4)
$$U(x,t) = \inf_{\xi:\xi(t)=x}\left\{\mathcal{A}_{t_0,t}(\xi) + \int_0^{\xi(t_0)} u_0(z)dz\right\}$$

and

(A.5)
$$U^\delta(x,t) = \inf_{\xi:\xi(t)=x}\left\{\mathcal{A}^\delta_{t_0,t}(\xi) + \int_0^{\xi(t_0)} u_0(z)dz\right\}.$$

It is clear that $U(x,t)$ is well-defined; i.e., the variational problem in (A.4) does have a solution. We will show that $u^\delta(x,t) \to u(x,t) = \frac{\partial U}{\partial x}(x,t)$ in $L^1_{\text{loc}}(\mathbb{S}^1 \times [t_0,\infty))$ as $\delta \to 0$ and consequently $u(x,t)$ is a weak entropy solution of (1.1). This follows from:

LEMMA A.1. *For almost all $\omega$, there exist $c_1$, $c_2$, $c_3$, depending only on $\omega$, $t$ and $t_0$, such that*

(A.6) $$\|U^\delta(\,\cdot\,,t) - U(\,\cdot\,,t)\|_{L^\infty(\mathbb{S}^1)} \leq c_1(\omega,t,t_0)\delta^{1/3},$$

(A.7) $$\|u^\delta(\,\cdot\,,t)\|_{L^\infty(\mathbb{S}^1)} \leq c_2(\omega,t,t_0),$$

*and*

(A.8) $$\|u^\delta(\,\cdot\,,t)\|_{\text{BV}(\mathbb{S}^1)} \leq c_3(\omega,t,t_0)$$

*where $\text{BV}(\mathbb{S}^1)$ is the space of functions on $\mathbb{S}^1$ with bounded variation.*

*Proof.* For any $\xi \in C^1[t_0,t]$ we have

$$|\mathcal{A}^\delta_{t_0,t}(\xi) - \mathcal{A}_{t_0,t}(\xi)| \leq \sum_k \int_{t_0}^t |f_k(\xi(s))|\,|\dot\xi(s)|\,|B_k^\delta(s) - B_k(s)|ds.$$

For almost every $\omega$, $\{B_k(\,\cdot\,)\}$ is $C^{\frac{1}{3}}$ for all $k$. Hence

$$|B_k^\delta(s) - B_k(s)| \leq C(\omega)\delta^{1/3}.$$

This gives

$$(A.9) \qquad |\mathcal{A}_{t_0,t}^\delta(\xi) - \mathcal{A}_{t_0,t}(\xi)| \leq \max_{t_0 \leq s \leq t} |\dot{\xi}(s)| C(\omega) \delta^{1/3} (t - t_0).$$

Denote by $\xi_\delta^*$ and $\xi^*$ the minimizers in (A.5) and (A.4) respectively. We have, using Lemma B.1,

$$U^\delta(x,t) - U(x,t) = \mathcal{A}_{t_0,t}^\delta(\xi_\delta^*) + \int_0^{\xi_\delta^*(t_0)} u_0(z)dz - \left(\mathcal{A}_{t_0,t}(\xi^*) + \int_0^{\xi^*(t_0)} u_0(z)dz\right)$$

$$\leq \mathcal{A}_{t_0,t}^\delta(\xi^*) + \int_0^{\xi^*(t_0)} u_0(z)dz - \left(\mathcal{A}_{t_0,t}(\xi^*) + \int_0^{\xi^*(t_0)} u_0(z)dz\right)$$

$$\leq \max_{t_0 \leq s \leq t} |\dot{\xi}^*(s)| C(\omega) \delta^{1/3} |t - t_0|$$

$$\leq C(\omega, t, t_0) \delta^{1/3}.$$

Similarly,

$$U(x,t) - U^\delta(x,t) \leq C(\omega, t, t_0) \delta^{1/3}.$$

To prove (A.7) we use the theory of backward characteristics (see [D]). If $(x,t)$ is a point of continuity of $u^\delta(\,\cdot\,,t)$, then there exists a unique backward characteristic $\gamma$ coming to $(x,t)$ and, for $s \in [t_0, t]$,

(A.10)
$$u^\delta(\gamma(s), s) = u_0(\gamma(t_0), t_0) + \sum_{k=1}^\infty \left(F_k(\gamma(s))B_k^\delta(s) - F_k(\gamma(t_0))B_k^\delta(t_0)\right)$$

$$- \int_{t_0}^s \sum_{k=1}^\infty f_k(\gamma(r)) u^\delta(\gamma(r), r) B_k^\delta(r) dr.$$

Hence

$$(A.11) \qquad |u^\delta(\gamma(s), s)| \leq |u_0(\gamma(t_0), t_0)| + c(\omega) + c(\omega) \int_{t_0}^s |u^\delta(\gamma(r), r)| dr$$

and $|u^\delta(\gamma(t), t)| \leq c_2(\omega, t - t_0)$. Since the points of continuity form a set of full measure in $\mathbb{S}^1$, we have (A.7).

Now consider two points of continuity for $u^\delta(\,\cdot\,,t)$, $x_1$ and $x_2$, and let $\gamma_1(s)$ and $\gamma_2(s)$ be the characteristics coming to $(x_1,t)$ and $(x_2,t)$ respectively. For

$i = 1, 2$ denote $u_i(s) = u^\delta(\gamma_i(s), s)$. Then

(A.12)
$$\frac{d}{ds}\left(\frac{u_1 - u_2}{\gamma_1 - \gamma_2}\right) = -\left(\frac{u_1 - u_2}{\gamma_1 - \gamma_2}\right)^2 + \frac{1}{\gamma_1 - \gamma_2}\frac{d}{ds}(u_1 - u_2)$$

$$= -\left(\frac{u_1 - u_2}{\gamma_1 - \gamma_2}\right)^2 + \sum_{k=1}^{\infty}\left(\int_0^1 f_k(\gamma_1 + r(\gamma_2 - \gamma_1))dr\right)dB_k^\delta(r)$$

by the mean value theorem and (2.3). This implies $\dfrac{u_1 - u_2}{\gamma_1 - \gamma_2} \leq c(\omega)$ since it solves an equation of the form $\dot{y} = -y^2 + C$. This, together with (A.7), gives (A.8). $\square$

LEMMA A.2. *For almost every $\omega$, the sequence $u^\delta$ converges in $L^1_{\text{loc}}(\mathbb{S}^1 \times [t_0, \infty))$ to a limit $u$ as $\delta \to 0$. Moreover $u(x, t) = \frac{\partial}{\partial x}U(x, t)$ and $u$ is an entropy-weak solution of (1.4).*

*Proof.* Integrating (A.1) on $\mathbb{S}^1 \times [t, t+\tau]$, we get

$$\int_0^1 dx|u^\delta(x, t+\tau) - u^\delta(x, t)| \leq \frac{1}{2}\int_t^{t+\tau} ds\|(u^\delta)^2\|_{\text{BV}(\mathbb{S}^1)}$$
$$+ \sum_k |B_k(t+\tau) - B_k(t)|\int_0^1 |f_k(x)|dx$$
$$\leq C_1(\omega)\tau + C_2(\omega)\tau^{1/3}.$$

In the last step, we used Lemma A.1 and Hölder continuity of the Wiener process. Hence $u^\delta$ is uniformly continuous in $t$, viewed as a function of $t$ in $L^1(\mathbb{S}^1)$. Therefore there exists a subsequence, still denoted by $u^\delta$, and $u \in L^\infty_{\text{loc}}([t_0, \infty), \text{BV}(\mathbb{S}^1)) \cap C([t_0, \infty), L^1(\mathbb{S}^1))$ such that

$$u^\delta \to u \quad \text{in} \quad L^1_{\text{loc}}(\mathbb{S}^1 \times [t_0, \infty)),$$

as $\delta \to 0$. From (A.6), we have

$$u = \frac{\partial U}{\partial x}.$$

From (A.8), the convergence also takes place in $L^p_{\text{loc}}(\mathbb{S}^1 \times [t_0, \infty))$ for $p < +\infty$. Hence $u$ is an entropy weak solution of (1.4). $\square$

Finally, observe that the solution operator for the mollified problem is order-preserving; i.e., $u_1^\delta(\,\cdot\,, t_0) \leq u_2^\delta(\,\cdot\,, t_0)$ implies $u_1^\delta(\,\cdot\, t) \leq u_2^\delta(\,\cdot\,, t)$ for $t \geq t_0$. Therefore the limiting solution, as $\delta \to 0$, is also order-preserving. Together with the conservation properties, we see that the solution operator is contractive in $L^1(\mathbb{S}^1)$ by the Crandall-Tartar lemma [CM]

(A.13) $\qquad \|u_1(\,\cdot\,, t) - u_2(\,\cdot\,, t)\|_{L^1(\mathbb{S}^1)} \leq \|u_1(\,\cdot\,, t_0) - u_2(\,\cdot\,, t_0)\|_{L^1(\mathbb{S}^1)}.$

This implies uniqueness of order-preserving weak solutions. In particular, since the solutions obtained in the zero-viscosity limit of (1.1) is also order-preserving, as a consequence of the comparison principle, we conclude that $u = \frac{\partial}{\partial x} U$ is the viscosity limit.

## Appendix B. Some technical estimates

Denote by $(x(t; x_0, v_0), v(t; x_0, v_0))$ the solution of (2.3) with initial data $x(0; x_0, v_0) = x_0$, $v(0; x_0, v_0) = v_0$. Sometimes we will also use the abbreviation $(x(t), v(t))$. Consider the stochastic flow $G_t^\omega$ defined by

$$(B.1) \qquad G_t^\omega(x_0, v_0) = (x(t; x_0, v_0), v(t; x_0, v_0)).$$

Since $f_k \in C^r$ the stochastic flow $G_t^\omega$ is $C^r$ smooth with probability 1. For $\tau > 0$ and $\omega \in \Omega$, define $\Gamma_\tau$ to be the set of $\tau$-minimizers,

$$(B.2) \qquad \Gamma_\tau = \{\gamma \in C^1[-\tau, 0]; \gamma(0), \gamma(-\tau) \in \mathbb{S}^1, \mathcal{A}_{-\tau,0}(\gamma) = \min_{\substack{\xi(0)=\gamma(0) \\ \xi(-\tau)=\gamma(-\tau)}} \mathcal{A}_{-\tau,0}(\xi)\}.$$

We shall also consider the case when endpoints belong to the universal cover $\mathbb{R}^1$, rather than $\mathbb{S}^1$. Denote

$$(B.3) \qquad \begin{aligned} \Gamma_{\tau,m} = \{&\gamma \in C^1[-\tau, 0]; 0 \leq \gamma(0) \leq 1, m \leq \gamma(-\tau) \leq (m+1), \\ & \mathcal{A}_{-\tau,0}(\gamma) = \min_{\substack{\xi(0)=\gamma(0) \\ \xi(-\tau)=\gamma(-\tau)}} \mathcal{A}_{-\tau,0}(\xi)\}. \end{aligned}$$

Of course $\Gamma_\tau, \Gamma_{\tau,m}$ depend on $\omega$. Let

$$V_\tau(\omega) = \sup_{\gamma \in \Gamma_\tau} |\dot\gamma(0)|,$$

$$V_{\tau,m}(\omega) = \sup_{\gamma \in \Gamma_{\tau,m}} |\dot\gamma(0)|,$$

$$\bar V_\tau(\omega) = \sup_{\gamma \in \Gamma_\tau} \max_{-\tau \leq s \leq 0} |\dot\gamma(s)|,$$

$$\bar V_{\tau,m}(\omega) = \sup_{\gamma \in \Gamma_{\tau,m}} \max_{-\tau \leq s \leq 0} |\dot\gamma(s)|.$$

In Lemma 3.3 it was shown that $V_\tau(\omega) \leq C(\omega)$ for $\tau \geq T(\omega)$. We consider now the case of small $\tau$.

LEMMA B.1. *There exists a constant $\tau(\omega)$ such that for $0 < \tau < \tau(\omega)$*

$$(B.4) \qquad V_\tau(\omega) \leq \bar V_\tau(\omega) \leq \frac{2}{\tau}.$$

*Furthermore, for any $m \in \mathbb{Z}^1$*

$$\text{(B.5)} \qquad V_{\tau,m}(\omega) \leq \bar{V}_{\tau,m}(\omega) \leq \frac{2(|m|+1)}{\tau}.$$

*Proof.* The proof is similar to the proof of Lemma 3.3. Let

$$\text{(B.6)} \qquad \|F\|_\tau = \max_{-\tau \leq s \leq 0} \sum_k \|F_k(x)\|_{C^3(\mathbb{S}^1)} |B_k(s) - B_k(0)|.$$

For arbitrary solutions of (2.3), $(x(t), v(t))$, $t \in [-\tau, 0]$, if we denote $\gamma(t) = x(t)$, $v_0 = |\dot{\gamma}(0)|$, $v = \max_{-\tau \leq t \leq 0} |\dot{\gamma}(t)|$, then we have

$$\text{(B.7)} \qquad v \leq v_0 + \|F\|_\tau + \|F\|_\tau \tau v$$

or

$$\text{(B.8)} \qquad v \leq \frac{v_0 + \|F\|_\tau}{1 - \|F\|_\tau \tau}$$

provided that $\|F\|_\tau \tau < 1$. Now

$$\text{(B.9)} \qquad |\dot{\gamma}(t) - \dot{\gamma}(0)| \leq \|F\|_\tau + \|F\|_\tau v\tau \leq \frac{\|F\|_\tau (1 + v_0 \tau)}{1 - \|F\|_\tau \tau}.$$

Assume that $\tau \leq 1$ is small enough so that $\|F\|_\tau \leq \epsilon$ (to be chosen later). Then we have

$$\text{(B.10)} \qquad v \leq \frac{v_0 + \epsilon}{1 - \epsilon},$$

$$\text{(B.11)} \qquad |\dot{\gamma}(t) - \dot{\gamma}(0)| \leq \frac{\epsilon(v_0 + 1)}{1 - \epsilon}.$$

Thus, provided that $v_0 \geq 1$, we have

$$\text{(B.12)} \qquad |\dot{\gamma}(t)| \geq \frac{1 - 3\epsilon}{1 - \epsilon} v_0$$

for $t \in [-\tau, 0]$. From (B.12), we get

$$\mathcal{A}_{-\tau,0}(\gamma) \geq \frac{1}{2}\left(\frac{1-3\epsilon}{1-\epsilon}\right)^2 v_0^2 \tau - \|F\|_\tau \tau v - \|F\|_\tau$$

$$\geq \frac{1}{2}\left(\frac{1-3\epsilon}{1-\epsilon}\right)^2 v_0^2 \tau - \frac{\epsilon}{1-\epsilon}(v_0 + 1).$$

Let $\gamma_1$ be the straight line such that $\gamma_1(0) = \gamma(0)$, $\gamma_1(-\tau) = \gamma(-\tau)$. Then

$$\mathcal{A}_{-\tau,0}(\gamma_1) \leq \frac{l^2}{2\tau} + l\|F\|_\tau + \|F\|_\tau \leq \frac{l^2}{2\tau} + (l+1)\epsilon,$$

where $l = |\gamma_1(0) - \gamma_1(-\tau)|$. Since $\mathcal{A}_{-\tau,0}(\gamma_1) \geq \mathcal{A}_{-\tau,0}(\gamma)$, one can easily show that $v_0 \leq \frac{3l}{2\tau}$ if $l \geq 1$ and $\epsilon \leq \frac{1}{40}$. If $l < 1$, then

$$\mathcal{A}_{-\tau,0}(\gamma_1) \leq \frac{1}{2\tau} + 2\|F\|_\tau \leq \frac{1}{2\tau} + 2\epsilon,$$

which together with $\mathcal{A}_{-\tau,0}(\gamma_1) \geq \mathcal{A}_{-\tau,0}(\gamma)$ gives $v_0 \leq \frac{3}{2\tau}$ if $\epsilon \leq \frac{1}{40}$ and $\tau \leq 1$. It follows that $V_\tau(\omega) \leq \frac{3}{2\tau}, \bar{V}_\tau(\omega) \leq \frac{2}{\tau}$ and $V_{\tau,m}(\omega) \leq \frac{3(|m|+1)}{2}, \bar{V}_{\tau,m}(\omega) \leq \frac{2(|m|+1)}{\tau}$.

In summary $\tau(\omega) \leq 1$ can be chosen such that $\|F\|_{\tau(\omega)} \leq \frac{1}{40}$. □

LEMMA B.2. *For any $K > 1$, there exists $\bar{\tau}(\omega) > 0$, such that for all $0 < \tau \leq \bar{\tau}(\omega)$*

(B.13) $$\|D^i G_t^\omega(x,v)\| \leq K$$

*for $1 \leq i \leq r$, $|v| \leq V_\tau(\omega) + 1$, $x \in \mathbb{S}^1$, $t \in [-\tau, 0]$.*

*Proof.* We will prove Lemma B.2 for $i = 1$. For $2 \leq i \leq r$ the proof is similar.

Consider the Jacobi matrix
$$DG_t^\omega = \begin{pmatrix} J_{11}(t), & J_{12}(t) \\ J_{21}(t), & J_{22}(t) \end{pmatrix}$$

where
$$J_{11}(t) = \frac{\partial x(t)}{\partial x_0}, \ J_{12}(t) = \frac{\partial x(t)}{\partial v_0}, \ J_{21}(t) = \frac{\partial v(t)}{\partial x_0}, \ J_{22}(t) = \frac{\partial v(t)}{\partial v_0}.$$

Obviously $(J_{11}, J_{12})$ and $(J_{21}, J_{22})$ satisfy

(B.14) $$\begin{cases} \dot{J}_{11}(t) = J_{21}(t) \\ \dot{J}_{21}(t) = \sum_k f_k(x(t)) J_{11}(t) \dot{B}_k(t), \end{cases}$$

$J_{11}(0) = 1$, $J_{21}(0) = 0$; and

(B.15) $$\begin{cases} \dot{J}_{12}(t) = J_{22}(t) \\ \dot{J}_{22}(t) = \sum_k f_k(x(t)) J_{12}(t) \dot{B}_k(t), \end{cases}$$

$J_{12}(0) = 0$, $J_{22}(0) = 1$.

Consider first (B.14). Let $J(\tau) = \max_{-\tau \leq s \leq 0} |J_{21}(s)|$. Then

(B.16) $$|J_{11}(s)| \leq 1 + J(\tau)\tau$$

for $s \in [-\tau, 0]$.

(B.17)
$$J_{21}(s) = \int_s^0 \sum_k f_k(x(t)) J_{11}(t) dB_k(t)$$
$$= -\sum_k f_k(x(s)) J_{11}(s) B_k(s) - \int_s^0 \sum_k f'_k(x(t)) v(t) J_{11}(t) B_k(t) dt$$
$$- \int_s^0 \sum_k f_k(x(t)) J_{21}(t) B_k(t) dt.$$

Using (B.8), (B.4) we have

$$\text{(B.18)} \qquad \max_{-\tau \leq s \leq 0} |v(s)| \leq \frac{v_0 + \|F\|_\tau}{1 - \|F\|_\tau \tau} \leq \frac{3}{\tau}$$

if $\tau$ is sufficiently small. Therefore

$$\text{(B.19)} \qquad J(\tau) \leq \|F\|_\tau \Big(1 + J(\tau)\tau\Big) + 3\|F\|_\tau \Big(1 + J(\tau)\tau\Big) + \|F\|_\tau \tau J(\tau).$$

It follows from (B.19) that

$$J(\tau) \leq \frac{4\|F\|_\tau}{1 - 5\|F\|_\tau \tau} \to 0$$

and

$$|J_{11}(s) - 1| \leq J(\tau)\tau \to 0$$

as $\tau \to 0$. Hence we have

$$|J_{11}(s)|, \ |J_{21}(s)| \leq K$$

for $s \in [-\tau, 0]$, if $\tau$ is sufficiently small. In the same way, we can prove

$$|J_{12}(s)|, \ |J_{22}(s)| \leq K$$

for $s \in [-\tau, 0]$ if $\tau$ is sufficiently small. This completes the proof of Lemma B.2. □

Denote $B_\tau = \{(x, v), x \in \mathbb{S}^1, |v| \leq V_\tau(\omega) + 1\}$, and $B_\tau(t) = G_t^\omega B_\tau$, for $t \in [-\tau, 0]$. Then, similar to Lemma B.2, we have:

LEMMA B.3. *For any $K > 1$, there exists $\tilde{\tau}(\omega) > 0$, such that for all $0 < \tau \leq \tilde{\tau}(\omega)$*

$$\text{(B.20)} \qquad \|D^i(G_t^\omega)^{-1}(x, v)\| \leq K$$

*for $1 \leq i \leq r$, $(x, v) \in B_\tau(t)$, $t \in [-\tau, 0]$.*

*Proof.* Lemma B.3 follows immediately from Lemma B.2 together with the estimate

$$\text{(B.21)} \qquad \|DG_t^\omega(x, v) - I\| \leq \max\left(\frac{5\|F\|_\tau}{1 - 5\|F\|_\tau \tau}, \frac{\tau}{1 - 5\|F\|_\tau \tau}\right)$$

for $t \in [-\bar{\tau}(\omega), 0]$, $(x, v) \in B_\tau$. (B.21) can be proved in the same way as (B.19). □

We will denote by $\tau_{r,K}(\omega)$ the maximum value of $\bar{\tau}(\omega)$ such that both (B.13) and (B.20) hold for all $\tau \leq \bar{\tau}(\omega)$.

Let $O_T(\omega) = \{(x_0, v_0), |v(t; x_0, v_0)| \leq \sup_{\gamma \in \Gamma_\tau} |v_\gamma(t)| + 1, \text{ for } t \in [-T, 0]\}$ and $O_T(t, \omega) = G_t^\omega O_T(\omega)$. Define

$$(B.22) \qquad D_{T,r}(\omega) = \sup_{-T \leq t \leq 0} \sup_{(x,v) \in O_T(\omega)} \max_{1 \leq i \leq r} \log^+ \|D^i G_t^\omega(x, v)\|,$$

$$\bar{D}_{T,r}(\omega) = \sup_{-T \leq t \leq 0} \sup_{(x,v) \in O_T(t,\omega)} \max_{1 \leq i \leq r} \log^+ \|D^i (G_t^\omega)^{-1}(x, v)\|.$$

LEMMA B.4. *For any positive integer $m$,*

$$(B.23) \qquad \int D_{T,r}(\omega)^m P(d\omega) < \infty, \quad \int \bar{D}_{T,r}(\omega)^m P(d\omega) < \infty.$$

*Proof.* The proof is similar to the proofs of similar statements in [Bax], [K2]. Define a sequence of stopping times $\tau_n$:

$$\tau_0 = 0, \quad \tau_1 = \tau_{r,K}(\omega), \quad \tau_{i+1} = \tau_{r,K}(\theta^{t_i} \omega),$$

where $t_i$ is defined by $t_0 = 0$, $t_1 = -\tau_1$, $t_i = -\sum_{j=0}^{i} \tau_j$. Choose $n$ such that $|t_n(\omega)| \leq T < |t_{n+1}(\omega)|$. Then

$$(B.24) \qquad G_{-T}^\omega(x, v) = f_{n+1} \circ f_n \circ \ldots \circ f_1(x, v)$$

where $f_1 = G_{-\tau_1}^\omega$, $f_2 = G_{-\tau_2}^{\theta^{t_1}\omega}, \ldots, f_i = G_{-\tau_i}^{\theta^{t_{i-1}}\omega}$ for $1 \leq i \leq n$, and $f_{n+1} = G_{t_n(\omega)-T}^{\theta^{t_n}\omega}$. Notice that for all $(x, v) \in O_T(\omega)$,

$$(B.25) \qquad \sup_{1 \leq i \leq r} \|D^i f_j((f_{j-1} \circ \ldots \circ f_1)(x, v))\| \leq K,$$

for $1 \leq j \leq n+1$, since $f_{j-1} \circ \ldots \circ f_1(x, v) \in B_{\tau_j}$. We now use the following fact. Consider $f \circ g(x, v)$. Assume that

$$\sup_{1 \leq i \leq r} \|D^i g(x, v)\| \leq M_1 \text{ and } \sup_{1 \leq i \leq r} \|D^i f(g(x, v))\| \leq M_2.$$

Then there exists a constant $C_r$ depending only on $r$ such that

$$\sup_{1 \leq i \leq r} \|D^i (f \circ g)(x, v)\| \leq C_r M_2 M_1^r.$$

Using (B.24) and (B.25), we obtain

$$(B.26) \qquad \sup_{(x,v) \in O_T(\omega)} \|D^i G_{-T}^\omega(x, v)\| \leq K(C_r K^r)^n$$

provided that $|t_n(\omega)| \leq T < |t_{n+1}(\omega)|$.

Since for any $t \in [-T, 0]$, we can write

$$(B.27) \qquad G_t^\omega(x, v) = \bar{f}_l \circ f_{l-1} \circ \ldots \circ f_1(x, v)$$

for some $l$, $1 \leq l \leq n+1$ and (B.25) holds, we get

(B.28) $$\sup_{-T \leq t \leq 0, (x,v) \in O_T(\omega)} \|D^i G_t^\omega(x,v)\| \leq K(C_r K^r)^n.$$

We now have

$$\int D_r(\omega)^m dP \leq \sum_{n=0}^\infty \log^m(K(C_r K^r)^n) P\{|t_n(\omega)| \leq T < |t_{n+1}(\omega)|\}.$$

Let $q = P\{\tau_{r,K}(\omega) \leq T\}$. Obviously $q < 1$ since there exists a set of $\omega$'s with positive probability such that $\tau_{r,K}(\omega) > 1$. Using the strong Markov property we have

$$P\{|t_n(\omega)| \leq T < |t_{n+1}(\omega)|\} \leq P(\tau_j \leq T, 1 \leq j \leq n)$$
$$= \prod_{j=1}^n P(\tau_j \leq T) = q^n.$$

Thus,

$$\int D_r(\omega)^m dP \leq \sum_{n=0}^\infty (\log K + n \log(C_r K^r))^m q^n < +\infty.$$

The other estimate can be proved in the same way. □

A stronger estimate holds for the first derivative $DG_t^\omega$. For $T > 0$, define

(B.29) $$d_T(\omega) = \sup_{-T \leq t \leq 0} \sup_{(x,v) \in B_T} \log^+ \|DG_t^\omega(x,v)\|,$$

(B.30) $$\bar{d}_T(\omega) = \sup_{0 \leq t \leq T} \sup_{(x,v) \in B_T} \log^+ \|DG_t^\omega(x,v)\|.$$

LEMMA B.5. *Let $m$ be a positive integer, then*

(B.31) $$\int_\Omega (d_T(\omega))^m dP < +\infty, \quad \int_\Omega (\bar{d}_T(\omega))^m dP < +\infty.$$

*Remark.* Lemma B.5 is stronger than Lemma B.4 for $r=1$ since $B_T \supset O_T$.

*Proof.* Let $x(t) = x(t; x_0, v_0)$, $v(t) = v(t; x_0, v_0)$. Now,

$$v(t) = v_0 - \int_t^0 \sum_k f_k(x(s)) dB_k(t)$$
$$= v_0 - \sum_k f_k(x(t)) B_k(t) + \int_t^0 \sum_k f_k'(x(s)) B_k(s) v(s) ds.$$

Let $M^\omega(t) = \max_{t \leq s \leq 0} |v(s)|$; then for $t \in [-T, 0]$

$$M^\omega(t) \leq |v_0| + \|F\|_T + \int_t^0 \|F\|_T M^\omega(s) ds.$$

This implies, for $t \in [-T, 0]$,
$$M^\omega(t) \leq (|v_0| + \|F\|_T)e^{\|F\|_T |t|}.$$

Consider next (B.14) and (B.15). Let $J^\omega(t) = \max_{t \leq s \leq 0} |J_{21}^\omega(s)|$, $\Delta_{t,\tau,k} = \sup_{t-\tau \leq s \leq t} |B_k(s) - B_k(t)|$, for $\tau > 0$, and $\|f_k\| = \sup_{0 \leq x \leq 1} |f_k(x)|$. Since
$$J_{21}^\omega(s) = J_{21}^\omega(t) - \int_s^t \sum_k f_k(x(u)) J_{11}^\omega(u) dB_k(u),$$
we have (without loss of generality, we can assume $B_k(t) = 0$)
$$J^\omega(t - \tau) \leq J^\omega(t) + \sup_{t-\tau \leq s \leq t} \left| \sum_k f_k(x(s)) B_k(s) J_{11}^\omega(s) \right|$$
$$+ \int_{t-\tau}^t \left| \sum_k f_k'(x(u)) B_k(u) v(u) J_{11}^\omega(u) \right| du$$
$$+ \int_{t-\tau}^t \left| \sum_k f_k(x(u)) B_k(u) J_{21}^\omega(u) \right| du$$
$$\leq J^\omega(t) + \left( \sum_k \|f_k\| \Delta_{t,\tau,k} \right)(1 + J^\omega(t - \tau)T)$$
$$+ (1 + J^\omega(t-\tau)T) M^\omega(T) 2\|F\|_T \tau + J^\omega(t-\tau) 2\|F\|_T \tau.$$

Choosing $\tau$ small enough so that

(B.32) $\quad \sum_k \|f_k\| \Delta_{t,\tau,k} \leq \frac{1}{6T},\ 2T M^\omega(T) \|F\|_T \tau < \frac{1}{6},\ 2\|F\|_T \tau \leq \frac{1}{6},$

we get

(B.33) $\quad J^\omega(t - \tau) \leq 2 \left( J^\omega(t) + \frac{1}{3T} \right).$

Define a sequence of stopping times $\bar{\tau}_i$, $i \geq 1$, by
$$\bar{\tau}_1 = \inf\left\{\tau: \sum_k \|f_k\| \Delta_{0,\tau,k} = \frac{1}{6T}\right\},$$
$$\bar{\tau}_{i+1} = \inf\left\{\tau: \sum_k \|f_k\| \Delta_{t_i,\tau,k} = \frac{1}{6T}\right\},$$
where $t_i = -\sum_{j=1}^i \bar{\tau}_j$. Assume that $|t_{k-1}| \leq T < |t_k|$. We can divide $[t_{i+1}, t_i]$, $0 \leq i \leq k - 1$, into subintervals such that (B.32) holds on each subinterval. The total number of these subintervals can be estimated from above by
$$R(k, T) = k + 12T\|F\|_T (1 + T M^\omega(T)).$$

From (B.33) we get

(B.34) $$J^\omega(T) \leq 2^{R(k,T)}\left(J^\omega(0) + \frac{2}{3T}\right) = \frac{2}{3T}2^{R(k,T)}.$$

We now show that

(B.35) $$\int R(k,T)^m dP < +\infty.$$

Denote $q = P(\bar{\tau}_1 \leq T) < 1$. Using the strong Markov property, we have

(B.36) $$P\{|t_{k-1}| \leq T \leq |t_k|\} \leq q^{k-1}.$$

Hence

(B.37) $$\int k^m P(d\omega) \leq \sum_{k=1}^{\infty} k^m P\{|t_{k-1}| \leq T \leq |t_k|\} \leq \sum_{k=1}^{\infty} k^m q^{k-1} < \infty.$$

On the other hand, since $M^\omega(T) \leq (|v_0| + \|F\|_T)e^{\|F\|_T T}$,

$$\int (M^\omega(T)\|F\|_T)^m P(d\omega) \leq \int (|v_0| + \|F\|_T)^{2m} e^{m\|F\|_T T} P(d\omega)$$
$$\leq \int (V_T(\omega) + 1 + \|F\|_T)^{2m} e^{m\|F\|_T T} P(d\omega).$$

Recall that $\|F\|_T = \sum_k \|F_k\|_{C^3} \max_{-T \leq t \leq 0} |B_k(t)|$. There exist constants $A, B > 0$, such that

(B.38) $$P(\|F\|_T \geq x) \leq Ae^{-Bx^2}.$$

Therefore for any positive integers $l$ and $m$

(B.39) $$\int \|F\|_T^l e^{m\|F\|_T T} P(d\omega) < \infty.$$

We also have from Lemma 3.3 and (B.38) that

$$\int V_T(\omega)^l e^{m\|F\|_T T} P(d\omega) < +\infty.$$

Hence we obtain

$$\int (M^\omega(T)\|F\|_T)^m P(d\omega) < +\infty.$$

An estimate for $J_{11}^\omega(t)$ follows from (B.16). Similar estimates can also be proved for $J_{12}^\omega(t)$ and $J_{22}^\omega(t)$. Together we obtain the first inequality in (B.31). The second inequality can be proved in the same way. □

Consider two minimizers $\gamma_1$ and $\gamma_2$ on $(-\infty, 0]$, $\gamma_1(0) = y$, $\gamma_2(0) = x$. Denote $v_1(\tau) = \dot{\gamma}_1(\tau)$, $v_2(\tau) = \dot{\gamma}_2(\tau)$.

LEMMA B.6. *Assume $y - x > 0$, $v_1(0) - v_2(0) = \tilde{L}(y - x)$ and $\tilde{L} > 4\|F\|_1(6 + 21\|F\|_1)$. Then*

$$v_1(t) - v_2(t) \geq 0 \tag{B.40}$$

*for $t \in [-\tau_0, 0]$, where $\tau_0 = \min\left(1, \frac{1}{2\|F\|_1}\right)$ and*

$$v_1(t) - v_2(t) > \frac{\tilde{L}}{4}(y - x) \tag{B.41}$$

*for $t \in [-\tau_1, 0]$, where $\tau_1 = \min\left(1, \frac{1}{4\|F\|_1}\right)$.*

*Proof.* We first prove (B.40). We shall consider $\gamma_1, \gamma_2$ as curves on the universal cover. Suppose that for some $-\tau_0 < t \leq 0$, $v_1(t) - v_2(t) = 0$. Denote by

$$t_1 = \max\{-\tau_0 \leq t \leq 0 \colon v_1(t) - v_2(t) = 0\},$$
$$t_2 = \min\{t \colon -t_1 \leq t \leq 0, v_1(t) - v_2(t) = \tilde{L}(y - x)\}.$$

Clearly $0 \leq v_1(t) - v_2(t) \leq \tilde{L}(y - x)$, $t_1 \leq t \leq t_2$. Also, since minimizers do not intersect, $0 \leq \gamma_1(t) - \gamma_2(t) \leq \gamma_1(0) - \gamma_2(0) = y - x$. Now,

$$v_1(t_1) = v_1(t_2) + \sum_k f_k(\gamma_1(t_1))B_k(t_1) - \sum_k f_k(\gamma_2(t_2))B_k(t_2)$$
$$+ \int_{t_1}^{t_2} \sum_k f'_k(\gamma_1(s))v_1(s)B_k(s)ds,$$

$$v_2(t_1) = v_2(t_2) + \sum_k f_k(\gamma_2(t_1))B_k(t_1) - \sum_k f_k(\gamma_2(t_2))B_k(t_2)$$
$$+ \int_{t_1}^{t_2} \sum_k f'_k(\gamma_2(s))v_2(s)B_k(s)ds.$$

Thus, $0 = v_1(t_1) - v_2(t_1) = v_1(t_2) - v_2(t_2) + \Delta v$, where

$$\Delta v = \sum_k (f_k(\gamma_1(t_1)) - f_k(\gamma_2(t_1)))B_k(t_1) - \sum_k (f_k(\gamma_1(t_2)) - f_k(\gamma_2(t_2)))B_k(t_2)$$
$$+ \int_{t_1}^{t_2} \sum_k (f'_k(\gamma_1(s)) - f'_k(\gamma_2(s)))v_1(s)B_k(s)ds$$
$$+ \int_{t_1}^{t_2} \sum_k f'_k(\gamma_2(s))(v_1(s) - v_2(s))B_k(s)ds.$$

Let $C_1 = \frac{1}{4} + \|F\|_1$, $C = 20C_1$. It follows from Lemma 3.3 and (B.8) that $|v_1(0)|, |v_2(0)| \leq C$ and for all $-\tau_0 \leq s \leq 0$:

$$|v_1(s)|, |v_2(s)| \leq 2(C + \|F\|_1) \leq 10 + 42\|F\|_1.$$

Thus,

$$|\Delta v| \leq 2\|F\|_1(y-x) + (t_2 - t_1)\|F\|_1(10 + 42\|F\|_1)(y-x)$$
$$+ (t_2 - t_1)\|F\|_1 \cdot \tilde{L}(y-x) \leq \left((12 + 42\|F\|_1) + \tau_0 \tilde{L}\right)\|F\|_1(y-x).$$

Since $\tau_0 \leq \frac{1}{2\|F_1\|}$, $(12 + 42\|F\|_1)\|F\|_1 < \frac{\tilde{L}}{2}$ we have from the estimate above:

$$|\Delta v| < \left(\frac{\tilde{L}}{2} + \frac{\tilde{L}}{2}\right)(y-x) < \tilde{L}(y-x),$$

which contradicts the fact that $|\Delta v| = \tilde{L}(y-x)$.

Next we prove (B.41). Suppose $-\tau_1 \leq t \leq 0$. Then $v_1(t) - v_2(t) \geq 0$. Suppose for some $-\tau_1 \leq t \leq 0$: $v_1(t) - v_2(t) = \frac{\tilde{L}}{4}(y-x)$. Denote

$$t_3 = \max\{-\tau_1 \leq t \leq 0 : v_1(t) - v_2(t) = \frac{\tilde{L}}{4}(y-x)\},$$
$$t_4 = \min\{-t_3 \leq t \leq 0 : v_1(t) - v_2(t) = \tilde{L}(y-x)\}.$$

Clearly, $\frac{\tilde{L}}{4}(y-x) \leq v_1(t) - v_2(t) \leq \tilde{L}(y-x)$, $0 \leq \gamma_1(t) - \gamma_2(t) \leq y-x$, $t_3 \leq t \leq t_4$. Using the same estimates as above, we have

$$\frac{\tilde{L}}{4}(y-x) = v_1(t_3) - v_2(t_3) = v_1(t_4) - v_2(t_4) + \overline{\Delta v},$$

where

$$|\overline{\Delta v}| \leq 2\|F\|_1(y-x) + (t_4 - t_3)\|F\|_1(10 + 42\|F\|_1)(y-x)$$
$$+ (t_4 - t_3)\|F\|_1 \tilde{L}(y-x)$$
$$\leq ((12 + 42\|F\|_1) + \tau_1 \tilde{L})\|F\|_1(y-x) < \left(\frac{\tilde{L}}{2} + \frac{\tilde{L}}{4}\right)(y-x) = \frac{3\tilde{L}}{4}(y-x).$$

On the other hand, $v_1(t_4) - v_2(t_4) = \tilde{L}(y-x)$, and $\overline{\Delta v} = -\frac{3}{4}\tilde{L}(y-x)$, which contradicts the estimate above. □

LEMMA B.7. *Let* $y - x > 0$, $v_1(0) - v_2(0) = \tilde{L}(y-x)$. *Then, with P-probability* 1, $\tilde{L} \leq \max(4\|F\|_1(6 + 21\|F\|_1), 4)$.

*Proof.* Suppose $\tilde{L} > 4\|F\|_1(6 + 21\|F\|_1)$. Then, for all $-\tau_1 \leq t \leq 0$, $v_1(t) - v_2(t) > \frac{\tilde{L}}{4}(y-x)$. Thus the two minimizers would intersect before the time $\tau_* = -\frac{(y-x)}{\frac{\tilde{L}}{4}(y-x)} = -\frac{4}{\tilde{L}}$, where $-\tau_1 \leq \tau_* \leq 0$ since

$$\tilde{L} > \max\left(4\|F\|_1(6 + 21\|F\|_1), 4\right).$$

This contradiction proves Lemma B.7. □

Denote

$$\|F\|_{-1,1} = \max_{-1 \leq s \leq 1} \sum_k \|F_k(x)\|_{C^3} |B_k(s) - B_k(0)|.$$

Let $L_0 = 4 + 24\|F\|_{-1,1} + 84\|F\|_{-1,1}^2$. Obviously,

$$L_0 > \max\left(4\|F\|_1(6 + 21\|F\|_1), 4\right).$$

Consider two minimizers at time $t = 1$: $\gamma_1(\tau), \gamma_2(\tau)$, $-\infty \leq \tau \leq 1$. Denote $y = \gamma_1(0)$, $v(y) = \dot\gamma_1(0)$, $x = \gamma_2(0)$, $v(x) = \dot\gamma_2(0)$.

LEMMA B.8. *With P-probability 1,*

$$|v(y) - v(x)| \leq L_0 |y - x|.$$

*Proof.* Suppose $y - x > 0$. Then, it follows from Lemma B.7 that: $v(y) - v(x) \leq L_0(y - x)$. Similarly we can prove an estimate from the other side. □

LEMMA B.9. *Consider the process*

$$\begin{aligned} dx &= v dt, \\ dv &= \sum_k f_k(x(t)) dB_k(t), \\ da &= b dt, \\ db &= a \sum_k f'_k(x(t)) dB_k(t), \end{aligned}$$

*and let $a(0)$, $b(0)$ satisfy $a(0)^2 + b(0)^2 = 1$. Assume that there exists a constant $\alpha_0 > 0$, such that*

$$\sum_k f'_k(x)^2 \geq \alpha_0$$

*for all $x \in [0,1]$. Then the joint probability distribution of $(a(1), b(1))$ has density $\bar{p}(a,b)$ which is uniformly bounded (with respect to $(a(0), b(0))$) on any compact domain.*

*Proof.* We will give only an outline for the proof. The generator $\mathcal{L}$ for the diffusion process can be written as

$$\mathcal{L} = \mathcal{L}_{a,b} + \mathcal{L}'$$

where

$$\mathcal{L}_{a,b} = b\frac{\partial}{\partial a} + \frac{1}{2}a^2 \left(\sum_k f'_k(x)^2\right) \frac{\partial^2}{\partial b^2},$$

$$\mathcal{L}' = v\frac{\partial}{\partial x} + a\sum_k f'_k(x)f_k(x)\frac{\partial^2}{\partial b \partial v}$$

$$+ \frac{1}{2}\sum_k f_k(x)^2 \frac{\partial^2}{\partial v^2}.$$

The operator $\mathcal{L}_{a,b}$ is hypoelliptic on $R^2 \setminus \{(0,0)\}$ for each fixed $x \in [0,1]$ (see [IK]). Therefore for each fixed $x \in [0,1]$, the solution of

$$\partial_t p^x = \mathcal{L}^*_{a,b} p^x,$$
$$p^x(a,b,0) = \delta(a - a(0), b - b(0))$$

is smooth for $t > 0$, except at $(a,b) = (0,0)$ [IK]. Since the delta function $p^x(\cdot, 0)$ is concentrated on the unit circle, we have that for $0 \leq t \leq 1$, $p^x$ is uniformly bounded (with respect to $x$ and $(a(0), b(0))$) on the circle

$$0 \leq p^x(a,b,t) \leq C^*$$

if $a^2 + b^2 = \frac{1}{4}$, and $0 \leq t \leq 1$. Using the maximum principle for the operator $\mathcal{L}_{a,b}$ on the domain $\{(a,b), a^2 + b^2 \leq \frac{1}{4}\} \times [0,1]$, we conclude that

$$p^x(a,b,t) \leq C^*$$

if $a^2 + b^2 \leq \frac{1}{4}$, and $t \leq 1$. Since $C^*$ is independent of $(x,v)$ and $(a(0), b(0))$, and since $p^x$ is smooth away from the origin, we obtain the desired result. □

## Appendix C. Hopf-Cole transformation and the Feynman-Kac formula

The Hopf-Cole transform and the Feynman-Kac formula are standard tools used in the analysis of (1.1). In the random case, some care has to be taken because of the appearance of stochastic integrals [S2].

Consider the stochastic PDE

$$(C.1) \qquad d\psi = \frac{\varepsilon}{2}\frac{\partial^2 \psi}{\partial x^2}dt + \left(-\frac{1}{\varepsilon}\sum_k F_k(x)dB_k(t) + c(x)dt\right)\psi$$

(the function $c(x)$ to be defined later). In the following, stochastic integrals will be understood in the Ito sense.

Let $v = -\varepsilon \ln \psi$. Using the Ito formula, we have

$$(C.2) \qquad dv = -\frac{\varepsilon^2}{2}\frac{1}{\psi}\frac{\partial^2 \psi}{\partial x^2}dt + \sum_k F_k(x)dB_k(t) + c(x)dt + \frac{1}{2\varepsilon}a(x)dt$$

where
$$a(x)dt = \mathbb{E}\left(\sum F_k(x)dB_k(t)\right)^2 = \left(\sum_k F_k^2(x)\right)dt.$$

Choosing $c(x) = -\frac{1}{2\varepsilon}a(x)$, we get

(C.3) $$dv = -\frac{\varepsilon^2}{2}\frac{1}{\psi}\frac{\partial^2 \psi}{\partial x^2}dt + \sum_k F_k(x)dB_k(t).$$

Let $u = -v_x$. It is straightforward to verify that $u$ satisfies

(C.4) $$du + \left(u\frac{\partial u}{\partial x} - \frac{\varepsilon}{2}\frac{\partial^2 u}{\partial x^2}\right)dt = \sum_k f_k(x)dB_k(t).$$

The Feynman-Kac formula for (C.1) takes the form

(C.5) $$\psi(x,t) = \mathbb{E}_\beta \psi\left(x + \sqrt{\varepsilon}\beta(t_0), t_0\right) e^{-\frac{1}{\varepsilon}\int_{t_0}^t \sum_k F_k\left(x+\sqrt{\varepsilon}\beta(s)\right)dB_k(s)},$$

where $\mathbb{E}_\beta$ denotes expectation with respect to the Wiener process on $[t_0, t]$ such that $\beta(t) = 0$. It is easy to verify that the extra terms that occur in the Ito formula for the exponential function in (C.5) are accounted for by the last term $c(x)\psi dt$ in (C.1).

(C.1) can also be rewritten as

(C.6) $$d\psi = \frac{\varepsilon}{2}\frac{\partial^2 \psi}{\partial x^2}dt - \frac{1}{\varepsilon}\psi \circ \sum_k F_k(x)dB_k(t)$$

where "$\circ$" denotes product in the Stratonovich sense.

## Appendix D. The basic collision lemma

This appendix is devoted to the proof and discussion of Lemma 5.2. We will use the notion of the backward Lagrangian map. It will be convenient to work with $\mathbb{R}^1$ instead of $\mathbb{S}^1$. Fix $t, s \in \mathbb{R}^1$, $t > s$, and $x \in \mathbb{R}^1$. Let $\xi_+$, $\xi_-$ be the maximal and minimal backward characteristics (see [D]) such that $\xi_+(t) = \xi_-(t) = x$. We define $Y_{s,t}^+(x) = \xi_+(s)$, $Y_{s,t}^-(x) = \xi_-(s)$.

We will study the case of
$$F(x,t) = -\frac{1}{2\pi}\cos(2\pi x)dB(t),$$
$$f(x,t) = \sin(2\pi x)dB(t).$$

It will be clear that the general situation follows from the same argument. From Lemma B.1, we can assume, without loss of generality, that $\|u(\cdot, 0)\|_{L^\infty} \leq C$ for some random constant $C$. Otherwise we change the initial time from $t = 0$ to some positive number, say $\frac{1}{16}$. It follows from Lemma B.1, that $\left\|u\left(\cdot, \frac{1}{16}\right)\right\|_{L^\infty} \leq C_1$ for some random constant $C_1$ depending on the forces on $[0, \frac{1}{16}]$. In addition, we will consider a particular case when $x_1^0 = \frac{1}{8}$, $x_2^0 = \frac{7}{8}$. It is easy to see from the proof that the argument works in the general case as well. We will use the notation $O(\delta)$ to denote quantities that are bounded in absolute value by $A\delta$, where $A$ is an absolute constant.

The basic strategy is to construct forces that are large on $[0, t_1]$ and small on $[t_1, 1]$ for some $t_1$ in order to set up approximately the following picture: At $t = t_1$, $u$ is very positive for $x \in [0, \frac{1}{2}]$ and very negative for $x \in [\frac{1}{2}, 1]$. If the forcing is small on $[t_1, 1]$, a shock must form which will absorb a sufficient amount of mass, if we imagine that there is a uniform distribution of masses on $[0, 1]$ at $t = t_1$. In order to make this intuitive picture rigorous, we must carefully control the value of $u$ when the forcing is small.

On $B$ and $t_1$, assume

$$\text{(D.1)} \qquad B(0) = 0, \max_{0 \le s \le t_1} |B(s)| \le 2B(t_1),\ 4\pi t_1 B(t_1) < \delta_0, B(t_1) > \bar{C}$$

with $\bar{C}, \delta_0$ as chosen below. We will show that if $B$ satisfies (D.1), then $x_1^0$ and $x_2^0$ merge before $t = 1$. Therefore the probability of merging is no less than the probability of the Brownian paths satisfying (D.1) which is positive.

Fix $x \in [0, 1]$. Let $\xi$ be a genuine backward characteristic emanating from $x$ at $t = t_1$; $\xi(t_1) = x$. Denote $y = \xi(0)$; then

$$\text{(D.2)} \quad u(x, t_1) = u(y, 0) + \sin(2\pi x) B(t_1) - \int_0^{t_1} 2\pi \cos(2\pi \xi(s)) \dot{\xi}(s) B(s) ds.$$

Hence $|u|_\infty = \max_{\substack{0 \le x \le 1 \\ 0 \le t \le t_1}} |u(x, t)|$ satisfies

$$|u|_\infty \le C + B(t_1) + 2\pi |u|_\infty \int_0^{t_1} |B(s)| ds.$$

Therefore

$$|u|_\infty \le M = \frac{C + B(t_1)}{1 - \delta_0}.$$

We now estimate $u(\cdot, t_1)$. The idea is that on the set where the force is bounded away from zero, $u$ is either very negative or very positive, reflected by the term involving $\delta_2$ below. We will bound $u$ on the complement of this set. For $x \in [\frac{1}{16}, \frac{1}{2} - \varepsilon]$, $0 < \varepsilon \ll 1$, $\varepsilon$ to be fixed later, we have

$$u(x, t_1) \ge -C + \delta_1 B(t_1) - 2\pi M \int_0^{t_1} |B(s)| ds$$
$$\ge -C + \delta_1 B(t_1) - \delta_0 M$$

with $\delta_1 = \sin\left(2\pi \left(\frac{1}{2} - \varepsilon\right)\right) = \sin(2\pi\varepsilon)$. For $x \in \left[\frac{1}{16}, \frac{7}{16}\right]$, with $\delta_2 = \sin\frac{\pi}{8}$, this can be improved to

$$u(x, t_1) \ge -C + \delta_2 B(t_1) - \delta_0 M.$$

The size of $\varepsilon$ is chosen such that there is a finite gap between $\delta_2$ and $\delta_1$. For $x \in \left[\frac{1}{2} - \varepsilon, \frac{1}{2}\right]$, we have

$$|u(x, t_1)| \leq C + \delta_1 B(t_1) + \delta_0 M.$$

Similarly on $\left[\frac{1}{2}, 1\right]$ we have the estimates

$$|u(x, t_1)| \leq C + \delta_1 B(t_1) + \delta_0 M,$$

for $x \in \left[\frac{1}{2}, \frac{1}{2} + \varepsilon\right]$;

$$u(x, t_1) \leq C - \delta_2 B(t_1) + \delta_0 M$$

for $x \in \left[\frac{9}{16}, \frac{15}{16}\right]$; and

$$u(x, t_1) \leq C - \delta_1 B(t_1) + \delta_0 M$$

for $x \in \left[\frac{1}{2} + \varepsilon, \frac{15}{16}\right]$.

Next on $[t_1, 1]$ we will choose $B(t)$ to be so small that

(D.3) $$\max_{t_1 \leq s \leq 1} |B(s) - B(t_1)| \leq \delta.$$

The value of $\delta$ will be chosen later. We first prove the following approximate monotonicity lemma.

LEMMA D.1. *Let $x^*$ be a point of shock at $t = 1$, $y_1 = Y^-_{t_1,1}(x^*)$, $y_2 = Y^+_{t_1,1}(x^*)$ and $y \in (y_1, y_2)$. Then*

$$\int_{y_1}^{y} (z + tu(z, t_1)) dz - x^*(y - y_1) \geq -C\delta |u|_\infty,$$

$$x^*(y_2 - y) - \int_{y}^{y_2} (z + tu(z, t_1)) dz \geq -C\delta |u|_\infty,$$

*where $t = 1 - t_1$, $|u|_\infty = \|u(\cdot, t_1)\|_{L^\infty}$.*

*Remark.* In the absence of forces, the correct statement is

$$\int_{y_1}^{y} (z + tu(z, t_1)) dz - x^*(y - y_1) \geq 0.$$

These statements were used in [ERS] as the basis for an alternative formulation of the variational principle. In the presence of force, similar statements appear to be invalid due to the presence of conjugate points. However, when the force is small, the error is also small, as claimed in Lemma D.1.

*Proof of Lemma* D.1. Define $y^*$ by:

(D.4) $$y^* + tu(y^*, t_1) = x^*.$$

Denote by $\xi_+$, $\xi_-$ the maximal and minimal backward characteristics such that $\xi_+(1) = \xi_-(1) = x^*$. Then

$$y_1 + \int_{t_1}^1 u(\xi_-(s), s) ds = x^*,$$

$$y_2 + \int_{t_1}^1 u(\xi_+(s), s) ds = x^*.$$

Furthermore, for $s \in [t_1, 1]$, we have

(D.5) $\quad u(\xi_+(s), s) = u(y_2, t_1) + \sin(2\pi\xi_+(s))(B(s) - B(t_1))$
$$- \int_{t_1}^s 2\pi \cos(2\pi\xi_+(s)) \dot{\xi}_+(s)(B(s) - B(t_1)) ds.$$

Hence

(D.6) $\quad |u(\xi_+(s), s) - u(y_2, t_1)| \leq O(\delta)|u|_\infty$

and

(D.7) $\quad y_2 + tu(y_2, t_1) = x^* + O(\delta)|u|_\infty.$

Similarly

(D.8) $\quad y_1 + tu(y_1, t_1) = x^* + O(\delta)|u|_\infty$

and

(D.9) $\quad |u(\xi_-(s), s) - u(y_1, t_1)| \leq O(\delta)|u|_\infty.$

From the action minimizing property of $\xi_-$, we get, by comparing the action of $\xi_-$ and $\xi(s) = y^* + (s - t_1)u(y^*, t_1)$,

$$\int_{t_1}^1 \frac{u^2(\xi_-(s), s)}{2} ds - \frac{1}{2\pi} \cos(2\pi x^*)(B(1) - B(t_1))$$
$$+ \int_{t_1}^1 \dot{\xi}_-(s) \sin(2\pi\xi_-(s))(B(s) - B(t_1)) ds$$
$$\leq t \frac{u^2(y^*, t_1)}{2} - \frac{1}{2\pi} \cos(2\pi x^*)(B(1) - B(t_1))$$
$$+ \int_{t_1}^1 u(y^*, t_1) \sin 2\pi(y^* + (s - t_1)u(y^*, t_1))(B(s) - B(t_1)) ds$$
$$+ \int_{y_1}^{y^*} u(y, t_1) dy.$$

This gives

$$t \frac{u^2(y_1, t_1)}{2} \leq t \frac{u^2(y^*, t_1)}{2} + \int_{y_1}^{y^*} u(y, t_1) dy + O(\delta)|u|_\infty.$$

Finally, we get, using (D.4), (D.7) and (D.8):

$$\int_{y_1}^{y^*} u(y,t_1)dy \geq \frac{t}{2}(u^2(y_1,t_1) - u^2(y^*,t_1)) + O(\delta)|u|_\infty\ ;$$

$$\int_{y_1}^{y^*} (z+tu(z,t_1))dz - (y^* - y_1)x^*$$

$$= \frac{(y^*)^2}{2} - \frac{y_1^2}{2} + t\int_{y_1}^{y^*} u(z,t_1)dz - (y^* - y_1)x^*$$

$$\geq \frac{(y^*)^2}{2} - \frac{y_1^2}{2} + \frac{t^2}{2}(u^2(y_1,t_1) - u^2(y^*,t_1)) - (y^* - y_1)x^* + O(\delta)|u|_\infty$$

$$= O(\delta)|u|_\infty\ . \qquad \square$$

*Proof of Lemma* 5.2. Let $z_1^0$ and $z_2^0$ be the Eulerian positions of $x_1^0$ and $x_2^0$ respectively at time $t_1$ following the forward characteristics defined by $u$. They are well-defined if the forward characteristics are continued properly by shocks [D]. Moreover, $|x_1^0 - z_1^0| \leq O(\delta_0), |x_2^0 - z_2^0| \leq O(\delta_0)$, since $t_1$ is small. Assume to the contrary that $z_1^0$ and $z_2^0$ do not merge until time 1. Then there exist $x_1$ and $x_2$, such that $x_1 < x_2$, and $z_1^0 \in [Y_{t_1,1}^-(x_1), Y_{t_1,1}^+(x_1)]$, $z_2^0 \in [Y_{t_1,1}^-(x_2), Y_{t_1,1}^+(x_2)]$. Let $\alpha_1 = Y_{t_1,1}^-(x_1)$, $\alpha_2 = Y_{t_1,1}^+(x_1)$, $\beta_1 = Y_{t_1,1}^-(x_2)$, $\beta_2 = Y_{t_1,1}^+(x_2)$. We have $\alpha_1 < \alpha_2 < \beta_1 < \beta_2$. Using the estimates obtained earlier on $u(\cdot, t_1)$, we have the following:

If $\alpha_2 < \frac{3}{8}$, then

$$\text{(D.10)} \qquad x_1 = \alpha_2 + (1-t_1)u(\alpha_2, t_1) + O(\delta)|u|_\infty$$

$$\geq \alpha_2 + (1-t_1)(-C + \delta_2 B(t_1) - \delta_0 M) + O(\delta)|u|_\infty\ .$$

Similarly, if $\beta_1 > \frac{5}{8}$, then

$$\text{(D.11)} \qquad x_2 = \beta_1 + (1-t_1)u(\beta_1, t_1) + O(\delta)|u|_\infty$$

$$\leq \beta_1 + (1-t_1)(C - \delta_2 B(t_1) + \delta_0 M) + O(\delta)|u|_\infty\ . \qquad \square$$

To deal with the case when either $\beta_1 < \frac{5}{8}$, or $\alpha_2 > \frac{3}{8}$, we introduce the parametrized measure $dQ_s(\cdot)$ which is the pullback of Lebesgue measure by the backward Lagrangian map from $t = t_1$ to $t = s$: $Q_s[x_1, x_2) = Y_{t_1,s}^+(x_2) - Y_{t_1,s}^-(x_1)$.

If we define $\rho = dQ_s$, $u(x,t) = \frac{1}{2}(u(x+,t) + u(x-,t))$, then it is easy to see that $(\rho, u)$ satisfies

$$\text{(D.12)} \qquad \rho_t + (\rho u)_x = 0$$

in the distributional sense.

Let $\xi_1$, $\xi_2$ be two genuine backward characteristics defined on $[t_1, 1]$, such that $\xi_1 < \xi_2$. Multiplying the above equation by $x$ and integrating over the region: $t_1 \leq s \leq 1$, $\xi_1(s) \leq x \leq \xi_2(s)$, we get

$$\text{(D.13)} \quad \int_{\xi_1(1)}^{\xi_2(1)} x dQ_1(x) - \int_{\xi_1(t_1)}^{\xi_2(t_1)} x dx = \int_{t_1}^{1} ds \int_{\xi_1(s)}^{\xi_2(s)} u(x,s) dQ_s(x).$$

LEMMA D.2. *For $s \in [t_1, 1]$,*

$$\int_{\xi_1(s)}^{\xi_2(s)} u(x,s) dQ_s(x) = \int_{\xi_1(t_1)}^{\xi_2(t_1)} u(y, t_1) dy + O(\delta)|u|_\infty.$$

*Proof.* First we assume that $x$ is a point of continuity of $u(\cdot, s)$. Then from the arguments presented earlier, we have

$$\text{(D.14)} \quad u(x,s) = u(y, t_1) + O(\delta)|u|_\infty$$

where $y = Y^+_{t_1,s}(x) = Y^-_{t_1,s}(x)$.

If $x$ is a point of discontinuity of $u(\cdot, s)$, let $y_1 = Y^-_{t_1,s}(x)$, $y_2 = Y^+_{t_1,s}(x)$. We then have

$$u(x,s) = \frac{1}{2}(u(x-,s) + u(x+,s))$$
$$= \frac{1}{2}(u(y_1, t_1) + u(y_2, t_1)) + O(\delta)|u|_\infty.$$

On the other hand, similar to the proof of Lemma D.1, we also have

$$y_1 + (s - t_1)u(y_1, t_1) = y_2 + (s - t_1)u(y_2, t_1) + O(\delta)|u|_\infty \frac{u^2(y_1, t_1)}{2}(s - t_1)$$
$$= \frac{u^2(y_2, t_1)}{2}(s - t_1) + \int_{y_1}^{y_2} u(y, t_1) dy + O(\delta)|u|_\infty.$$

Hence

$$\int_{y_1}^{y_2} u(y, t_1) dy = \frac{s - t_1}{2}(u^2(y_1, t_1) - u^2(y_2, t_1)) + O(\delta)|u|_\infty$$
$$= \frac{s - t_1}{2}(u(y_1, t_1) - u(y_2, t_1))(u(y_1, t_1) + u(y_2, t_1)) + O(\delta)|u|_\infty$$
$$= (y_2 - y_1)\frac{u(y_1, t_1) + u(y_2, t_2)}{2} + O(\delta)|u|_\infty.$$

This can be written as

(D.15) $$u(x,s)Q_s(\{x\}) = \int_{y_1}^{y_2} u(y,t_1)dy + O(\delta)|u|_\infty.$$

Now Lemma D.2 follows from (D.14) and (D.15) when we use a standard approximation argument. □

We now continue with the proof of Lemma 5.2. Using Lemma D.2, we can rewrite (D.13) as

$$\int_{\xi_1(1)}^{\xi_2(1)} x\,dP_1(x) - \int_{\xi_1(t_1)}^{\xi_2(t_1)} x\,dx = (1-t_1)\int_{\xi_1(t_1)}^{\xi_2(t_1)} u(y,t_1)dy + O(\delta)|u|_\infty.$$

Applying this to $Y_{t_1,1}^\pm(x_1)$, $Y_{t_1,1}^\pm(x_2)$, we get

$$x_1 - \frac{1}{\alpha_2 - \alpha_1}\int_{\alpha_1}^{\alpha_2}\{x + (1-t_1)u(x,t_1)\}dx = \frac{O(\delta)|u|_\infty}{\alpha_2 - \alpha_1},$$

$$x_2 - \frac{1}{\beta_2 - \beta_1}\int_{\beta_1}^{\beta_2}\{x + (1-t_1)u(x,t_1)\}dx = \frac{O(\delta)|u|_\infty}{\beta_2 - \beta_1}.$$

Assume that $\alpha_2 > \frac{3}{8}$. Then $\alpha_2 - z_1^0 > \frac{1}{4} - O(\delta_0)$. Integrating both sides of the equation $u_t + \left(\frac{u^2}{2}\right)_x = -F_x$ over the region: $0 \le t \le t_1$, $\xi_1(t) \le x \le \xi_2(t)$, where $\xi_1(t) = Y_{t,t_1}^-(x_1^0)$, $\xi_2(t) = Y_{t,t_1}^+(\alpha_2)$, we get

(D.16)
$$\int_{x_1^0}^{\alpha_2} u(x,t_1)dx - \int_{\xi_1(0)}^{\xi_2(0)} u(x,0)dx = \frac{1}{2}\int_0^{t_1}[u(\xi_1(t),t)^2 - u(\xi_2(t),t)^2]dt$$

$$- \frac{1}{2\pi}B(t_1)\Big(\cos(2\pi\xi_2(t_1)) - \cos(2\pi\xi_1(t_1))\Big)$$

$$+ \int_0^{t_1} B(t)\Big(\dot\xi_2(t)\sin(2\pi\xi_2(t)) - \dot\xi_1(t)\sin(2\pi\xi_1(t))\Big)dt.$$

This implies that

$$\int_{x_1^0}^{\alpha_2} u(x,t_1)dx \ge -\frac{1}{2\pi}B(t_1)\Big(\cos(2\pi\alpha_2) - \cos(2\pi x_1^0)\Big)$$

$$- C - t_1|u|_\infty^2 - |u|_\infty\int_0^{t_1}|B(s)|ds.$$

Hence, using Lemma D.1, we obtain

$$\text{(D.17)} \quad x_1 \geq \frac{1}{\alpha_2 - x_1^0} \int_{x_1^0}^{\alpha_2} (z + (1-t_1)u(z,t_1))dz + O(\delta)|u|_\infty$$

$$\geq -\frac{2}{2\pi}(1-t_1)B(t_1)\Big(\cos(2\pi\alpha_2) - \cos(2\pi x_1^0)\Big)$$

$$- C - 4t_1|u|_\infty^2 - 4\delta_0|u|_\infty + O(\delta)|u|_\infty.$$

Similarly, if $\beta_1 < \frac{5}{8}$,

$$\text{(D.18)} \quad x_2 \leq \frac{1}{x_2^0 - \beta_1} \int_{\beta_1}^{x_2^0} \Big(z + (1-t_1)u(z,t_1)\Big)dz + O(\delta)|u|_\infty$$

$$\leq \frac{2}{2\pi}(1-t_1)B(t_1)(\cos 2\pi x_2^0 - \cos 2\pi \beta_1)$$

$$+ C + 4t_1|u|_\infty^2 + 4\delta_0|u|_\infty + O(\delta)|u|_\infty.$$

If $\beta_1 \geq \frac{3}{8}$, we have $\cos(2\pi x_2^0) - \cos(2\pi\beta_1) \geq 0$, and if $\alpha_2 \leq \frac{5}{8}$, we have $\cos(2\pi\alpha_2) - \cos(2\pi x_1^0) \leq 0$. Otherwise, we can use (D.10) and (D.11). In any case, we always have, for some positive constant $C^*$,

$$\text{(D.19)} \quad x_1 - x_2 \geq C^*B(t_1) - C_0(1 + t_1|u|_\infty^2 + \delta_0|u|_\infty + \delta|u|_\infty)$$

$$\geq C^*B(t_1) - C_0 - C_0(4\delta_0 + \delta)|u|_\infty$$

$$\geq \left[C^* - \frac{2C_0(4\delta_0 + \delta)}{1-\delta_0}\right]B(t_1) - C_1.$$

The constants $C^*$, $C_0$, $C_1$ do not depend on $\delta_0$, $\delta$, $B(t_1)$.

If we choose $\delta_0$, $\delta$, such that

$$\text{(D.20)} \quad C^* - \frac{2C_0(4\delta_0 + \delta)}{1-\delta_0} > 0$$

we can then choose $\bar{C}$, such that

$$x_1 - x_2 > 0,$$

contradicting the assumption that $x_1 \leq x_2$. This completes the proof of Lemma 5.2.

We now estimate the location of the shock $x^*$ where $x_1^0$ and $x_2^0$ have merged at $t = 1$, assuming that the forces are chosen as in the proof of Lemma 5.2. Let $y_1 = Y_{t_1,1}^-(x^*)$, $y_2 = Y_{t_1,1}^+(x^*)$. Now,

$$x^* = \frac{1}{y_2 - y_1} \int_{y_1}^{y_2} (y + (1-t_1)u(y,t_1))dy + O(\delta)|u|_\infty.$$

Also,

(D.21)

$$\left| \frac{1}{y_2 - y_1} \int_{y_1}^{y_2} y\, dy - \frac{1}{2} \right| < \frac{1}{2}(1 - x_2^0 + x_1^0),$$

$$\left| \int_{y_1}^{y_2} u(y, t_1) dy \right| \leq (1 - x_2^0 + x_1^0)|u(\cdot\, .0)|_\infty + t_1 |u|_\infty^2$$

$$+ \frac{1}{2\pi} |(\cos 2\pi y_2 - \cos 2\pi y_1) B(t_1)| + \delta_0 |u|_\infty,$$

where we used an analog of (D.16). The factors $1 - x_2^0 + x_1^0$, $\frac{1}{2}(1 - x_2^0 + x_1^0)$ can be made arbitrarily small by choosing $x_1^0$ close to 0, and $x_2^0$ close to 1.

Notice that in (D.21) the coefficient in front of $B(t_1)$ is approximately equal to $|y_1 + 1 - y_2| \sin 2\pi \bar{y}$ for some $\bar{y} \in (y_2, y_1 + 1)$, whereas $C^*$ in (D.19) is bounded from below by $\min(\sin 2\pi x_1^0, \sin 2\pi x_2^0)$. Therefore by choosing $x_1^0$ close to 0, $x_2^0$ close to 1, and $B(t_1)$ such that (D.20) holds but $|\cos 2\pi y_2 - \cos 2\pi y_1| B(t_1)$ is small, we can make $x^*$ arbitrarily close to $\frac{1}{2}$. We have arrived at:

LEMMA D.3. *Assume that $F(x, t) = -\frac{1}{2\pi} \cos(2\pi x) dB(t)$. Fix any $\varepsilon_1$, $\varepsilon_2 > 0$. Then the following event has positive probability $p_0(\varepsilon_1, \varepsilon_2)$. There exists $x^* \in \left[ \frac{1}{2} - \varepsilon_1, \frac{1}{2} + \varepsilon_1 \right]$, such that $[\varepsilon_2, 1 - \varepsilon_2] \subset [Y_{0,1}^-(x^*), Y_{0,1}^+(x^*)]$. In other words, the interval $[\varepsilon_2, 1 - \varepsilon_2]$ is mapped to a point $x^* \in \left[ \frac{1}{2} - \varepsilon_1, \frac{1}{2} + \varepsilon_1 \right]$ by the forward Lagrangian map.*

To prove this, we just have to take an $\varepsilon_3 < \varepsilon_2$, and $x_1^0 = \varepsilon_3$, $x_2^0 = 1 - \varepsilon_3$ and use the argument outlined above. We omit the details.

*Acknowledgement.* We have discussed the results of this paper with many people. We are especially grateful to P. Baxendale, U. Frisch, Yu. Kifer, J. Mather, J. Mattingly, J. Moser, A. Polyakov, T. Spencer, S. R. S. Varadhan and V. Yakhot for their suggestions and comments. The work of Weinan E is partially supported by an NSF Presidential Faculty Fellowship. The work of Khanin and Mazel is supported in part by RFFI of Russia under grant 96-01-00377, with further support for Khanin from Leverhulme Trust Research Grant RF&G/4/9700445. The work of Sinai is supported in part by NSF grant DMS-9706794.

Courant Institute, New York University, New York, NY
*Current address*: Princeton University, Princeton, NJ
*E-mail address*: weinan@math.princeton.edu

Isaac Newton Institute, University of Cambridge, England, UK
BRIMS, Hewlett-Packard Laboratories, Bristol, UK
Herrot-Watt University, Edinburgh, Scotland
Landau Institute of Theoretical Physics, Moscow, Russia
*E-mail address*: K.Khanin@newton.cam.ac.uk

Rutgers University, Piscataway, NJ
*E-mail address*: mazel@math.rutgers.edu

Princeton University, Princeton, NJ
Landau Institute of Theoretical Physics, Moscow, Russia
*E-mail address*: sinai@math.princeton.edu